\documentclass[3p,times]{elsarticle}

\usepackage{ecrc}
\usepackage{amscd}

\volume{00}

\firstpage{1}

\journalname{}

\runauth{}

\jid{procs}

\jnltitlelogo{}

\usepackage{amssymb}

\usepackage[figuresright]{rotating}
\usepackage{mhchem}
\usepackage{esint}
\usepackage{mathtools}
\usepackage{bm}
\usepackage{amsmath}
\allowdisplaybreaks
\numberwithin{equation}{section}
\numberwithin{figure}{section}
\numberwithin{table}{section}

\usepackage{graphicx }
\usepackage{subfigure,txfonts}
\usepackage{amsthm}
\usepackage{amssymb}
\usepackage{cancel}
\usepackage{mathdots}
\usepackage{stmaryrd}
\usepackage{stackrel}
\usepackage{graphicx}
\usepackage{booktabs}
\usepackage{lineno}
\usepackage{nameref}
\usepackage{xcolor}
\usepackage[colorlinks,
linkcolor=blue,
anchorcolor=yellow,
citecolor=red,
]{hyperref}

\usepackage{version}
\usepackage{cases}
\usepackage{color}

    \def\H{\boldsymbol{H}}
	\def\L{\boldsymbol{L}}
	\def\W{\boldsymbol{W}}
    \def\X{\boldsymbol{X}}
    \def\Y{\boldsymbol{Y}}
    \def\Z{\boldsymbol{Z}}
    \def\I{\boldsymbol{I}}
    \def\P{\boldsymbol{P}}

    \def\PPi{\bm \Pi}

    \def\B{\boldsymbol{B}}
    \def\C{\bm \zeta}
    
    \def\u{\boldsymbol{u}}
    \def\v{\boldsymbol{v}}
    \def\w{\boldsymbol{w}}
    
    \def\0{\boldsymbol{0}}
    \def\x{\boldsymbol{x}}
    \def\y{\boldsymbol{y}}
    \def\n{\boldsymbol{n}}

\usepackage{lineno}
\makeatletter

\theoremstyle{plain}
\ifx\thesection\undefined
\newtheorem{thm}{\protect\theoremname}
\else

\fi
\theoremstyle{remark}
\ifx\thesection\undefined
\newtheorem{rem}{\protect\remarkname}

\else

\fi
\theoremstyle{plain}
\ifx\thesection\undefined
\newtheorem{lem}{\protect\lemmaname}
\else

\fi
\theoremstyle{plain}
\ifx\thesection\undefined
\newtheorem{Example}{\protect\lemmaname}
\else

\fi

\graphicspath{{./figures/}}

\begin{document}
\newtheorem{The}{Theorem}[section]
\newtheorem{lemma}{Lemma}[section]
\newtheorem{Remark}{Remark}[section]
\newtheorem{Assumption}{Assumption}[section]
\newtheorem{Definition}{Definition}[section]
\newtheorem{Proposition}{Proposition}[section]

\newcommand{\sign}{\mathrm{sign}}

\begin{frontmatter}



\dochead{}
\title{Optimal ${L^2}$ error estimates for 2D/3D incompressible Cahn--Hilliard--magnetohydrodynamic equations \tnoteref{t1}}

\tnotetext[t1]{This work is partly supported by the NSF of China (No. 12061076, 12171340 ) support.
}

\author[zzu]{Haiyan Su}%
\ead{shymath@126.com }
\cortext[cor1]{Corresponding author. Tel./fax number: 86 9918582482.}

\author[zzu1]{Jilu Wang\corref{cor1}}
\ead{(wangjilu03@gmail.com, wangjilu@hit.edu.cn}

\author[zzu2]{Zeyu Xia}
\ead{zeyuxia@uestc.edu.cn}

\author[zzu]{Ke Zhang}
\ead{zhangkemath@139.com, zkmath@stu.xju.edu.cn}
\address[zzu]{College of Mathematics and System Sciences,
Xinjiang University, Urumqi 830046, P.R. China}
\address[zzu1]{School of Science, Harbin Institute of Technology, Shenzhen 518055,
China}
\address[zzu2]{School of Mathematical Sciences, University of Electronic Science and Technology of China, Chengdu 611731, China.}

\begin{abstract}
This paper focuses on an optimal error analysis of a fully discrete finite element scheme for the Cahn--Hilliard--magnetohydrodynamic (CH-MHD) system. The method  use the standard inf-sup stable Taylor--Hood/MINI elements to solve the Navier--Stokes equations, Lagrange elements to solve the phase field, and particularly, the N\'ed\'elec elements for solving the magnetic induction field. Suffering from the strong coupling and high nonlinearity, the previous works just provide  suboptimal error estimates for phase field and velocity field in  $L^{2}/\L^2$-norm under the same order  elements, and the suboptimal error estimates for magnetic induction field in $\H(\rm curl)$-norm.  To this end, we utilize the Ritz, Stokes, and Maxwell quasi-projections to eliminate the low-order pollution of the phase field and magnetic induction field. In addition to the optimal $\L^2$-norm error estimates, we present the optimal convergence rates for magnetic induction field  in $\H(\rm curl)$-norm and for velocity field in $\H^1$-norm.
Moreover, the unconditional energy stability and mass conservation of the proposed scheme are preserved. Numerical examples are illustrated to validate the theoretical analysis and show the performance of the proposed scheme.
\end{abstract}


\end{frontmatter}
 
\section{Introduction}

The diffuse interface model of the magnetohydrodynamic system describes the dynamic behavior of two incompressible and immiscible conducting fluids under an external magnetic field. The governing equation consists of the Cahn--Hilliard equations, Navier--Stokes equations, and the Maxwell's equations, which are coupled  through convection, stresses, and Lorentz forces. Such a model has extensive application prospects in the fields of nuclear fusion, metallurgy, liquid metal magnetic pumps, aluminum electrolysis and so on \cite{2010An, Jean2006Mathematical, 2000Liquid}.

In this paper, we consider the following CH-MHD model \cite{2019A, 2023Error}: for $(\boldsymbol{x}, t )\in  \Omega\times (0,T]$,
\begin{align}
&\phi_{t}+\nabla \phi \cdot\u=\gamma \nabla\cdot (M \nabla \omega), \label{PDE1} \\
&-\gamma\Delta\phi+\gamma^{-1}(\phi^{3}-\phi) = \omega,  \label{PDE2} \\
&\u_{t}+(\u\cdot \nabla) \u- \nabla\cdot (\nu  \nabla\u) + \nabla p + \mu^{-1}\B\times\nabla\times\B = \lambda \omega \nabla \phi,	 	\label{PDE3}	\\
&\nabla \cdot \u=0, 	\label{PDE4} \\
&\B_{t}+ \mu^{-1}\nabla\times(\sigma^{-1}\nabla\times\B)- \nabla\times(\u\times\B) = \0, 	\label{PDE5}
\end{align}
where $\Omega$ is a bounded smooth domain in $\mathbb{R}^{d}$ with the spatial dimension $d=2, 3$, $\n$ denotes the unit outward normal on the boundary $\partial\Omega$ of domain $\Omega$, and $T$ represents the terminal time. 
This system is solved subject to the following boundary and initial conditions
\begin{align}
	& \frac{\partial\phi}{\partial \boldsymbol{n}}\Big|_{\partial\Omega}=0, 
			\quad 
		\frac{\partial w}{\partial \boldsymbol{n}}\Big|_{\partial\Omega}=0, 
			\quad 
		\u|_{\partial\Omega}=\0,
			\quad
		\B\times \boldsymbol{n}|_{\partial\Omega}=\0,	\qquad \mbox{on}\ \ \partial\Omega,	\\
	&\phi|_{t=0}=\phi_0,
			\quad 
		\u|_{t=0}=\u_0,
			\quad 
		\B|_{t=0}=\B_0,		\hspace{1.14in}			\mbox{in}	\ \ \Omega\times\{0\}.
\end{align}
In equation \eqref{PDE1}-\eqref{PDE5}, $\phi$ is the phase field used to distinguish the mixture of two incompressible immiscible fluids, $(\u, p)$ represents the velocity-pressure pair, $\omega$ denotes the chemical potential, and $\B$ stands for the magnetic induction field. In this problem, we consider the physical parameters to be positive constants, including interfacial width $\gamma$ between two phases, mobility parameter $M$, kinematic viscosity $\nu$, magnetic permeability $\mu$, electric conductivity $\sigma$, and capillary coefficient $\lambda$.
 
To solve such a two-phase MHD system, there already exists an extensive literature on the energy stability, which is an essential requirement when designing a numerical scheme. For example, the pressure correct method \cite{2022Highly, wang2024fully} and Gauge-Uzawa method \cite{ZHANG2023107477} are developed to decouple the velocity and pressure fields. In addition, there are lots of efficient algorithms to handle the strong coupling and high nonlinearity arising from the phase field, such as the invariant energy quadratization (IEQ) method \cite{2022Highly, ZHANG2023107477, 2023Energy}, the semi-implicit stabilization method \cite{2022Highly, ZHANG2023107477, chen2022unconditional},  the convex splitting method \cite{2019A, 2023Error, 2025Error}, and the scalar auxiliary variables (SAV) method \cite{wang2024fully, 2024ErrorJIA}. In this paper, we mainly focus on the convex splitting scheme, which was proposed in \cite{1998UnconditionallyE} and has been widely used in practice \cite{2014ExistenceHAN, 2010UnconditionallyWISE}. 


Additionally, convergence analysis of numerical methods for the CH-MHD system remains an active area of research, especially concerning the combination of Taylor--Hood/MINI elements for  the velocity-pressure pair \cite{2023Error, wang2024convergence,  2025Error} and  N\'{e}d\'{e}lec edge elements for magnetic induction fields \cite{2022Highly, ZHANG2023107477, 2023Energy, yang2025unconditionally}. However, since the high-order N\'{e}d\'{e}lec edge elements are quite complicated in implementation and time-consuming to compute, it is preferable to choose the lower N\'{e}d\'{e}lec elements than the Lagrange elements, which thus may bring the accuracy pollution in the convergence analysis.
In \cite{2023New, huang2023new}, the authors define the Maxwell quasi-projections for the incompressible MHD system to eliminate  the pollution of the lower-order N\'{e}d\'{e}lec edge element approximation. Then, in \cite{yang2025unconditionally} the authors propose a fully discrete finite element scheme based on the ``zero-energy-contribution'' method for the CH-MHD system, and provide the following error estimate: 
\begin{align}
&\|\phi^{k+1}-\phi_{h}^{k+1} \| + \|\u^{k+1}-\u_{h}^{k+1}\| 
		\leq C\Big(\Delta t+h^{l+1}+h^{r+2} \Big), \label{N1}\\
&\|\B^{k+1}-\B_{h}^{k+1}\| 
		\leq C\Big(\frac{\Delta t}{h}+h^{l}+h^{r+1} \Big), 
				\quad 
	\bigg(\Delta t\sum_{k=0}^{N_T-1}  \|\nabla\times(\B^{k+1}-\B_{h}^{k+1} )\|^{2} \bigg)^{\frac{1}{2}}
		\leq C\Big(\frac{\Delta t}{h}+h^{l}+h^{r+1}\Big) \label{BL2}\\
&\bigg(\Delta t\sum_{k=0}^{N_T-1}  \|\nabla (\u^{k+1}-\u_{h}^{k+1} )\|^{2} \bigg)^{\frac{1}{2}}
	\leq C\Big(\frac{\Delta t}{h}+h^{l}+h^{r+1}\Big). \label{vb}
\end{align}
In their methods, they look for the solution $(\phi_{h}^{k+1}, \omega_{h}^{k+1}, \u_{h}^{k+1}, p_{h}^{k+1}, \B_{h}^{k+1})$ in the finite element space $S_{h}^{r+1}\times S_{h}^{r+1}\times \X_{h}^{r+1}\times \mathring{S}_{h}^{r}\times \Y_{h}^{l}$ (the definitions of spaces seen in Section \ref{sec-main}).
It is observed  in \eqref{N1} that, with the help of Maxwell quasi-projections matched variable electric conductivities defined in \cite{yang2025unconditionally}, the accuracy of velocity field can be $\mathcal{O}(h^{l+1})$ when using the $l$-th order N\'{e}d\'{e}lec elements.
However, for the phase field $\phi$, to achieve the optimal error estimate, they should be solved by using the matched elements to the velocity field. Thus, it remains to consider how to eliminate the artificial pollution arising from the phase field when utilizing the lower-order elements. Moreover, it is seen that the error estimate of $\H(\rm{curl})$-norm for the magnetic induction field in \eqref{BL2} and $\H^{1}$-norm for velocity field in \eqref{vb} are $\mathcal{O}(\frac{\Delta t}{h})$, whereas the desired orders are expected to be $\mathcal{O}(\Delta t)$. 
  
$\bullet$ In this paper, our first goal is to develop optimal $L^{2}$-norm error estimates of all variables by employing the Ritz, Stokes and Maxwell quasi-projections.  We look for solutions $(\phi_{h}^{k+1}, \omega_{h}^{k+1}$, $\u_{h}^{k+1}$, $p_{h}^{k+1}$, $\B_{h}^{k+1})\in S_{h}^{r}\times S_{h}^{r}\times \X_{h}^{r+1}\times \mathring{S}_{h}^{r}\times \Y_{h}^{l}$ and then obtain the  error estimates
\begin{equation}\label{temporal}
\left\{
\begin{aligned}
\|\phi^{k+1}-\phi_{h}^{k+1} \|
	&\leq C\big(\Delta t+h^{l+1}+h^{r+1}\big),	\\
\|\u^{k+1}-\u_{h}^{k+1}\|
	&\leq C\big(\Delta t+h^{l+1}+\beta_{h}\big),\\
\|\B^{k+1}-\B_{h}^{k+1}\|
	&\leq C\big(\Delta t+h^{l}+\beta_{h}\big), 
\end{aligned}
\right.
\qquad  
\beta_{h} = 
\left\{
\begin{aligned}
	h^{r+2},\ \  &r\geq2,\\
	h^{r+1},\ \  &r=1.
\end{aligned}
	\right.
\end{equation}
Clearly, for the MINI element and lowest-order N\'{e}d\'{e}lec edge element, namely, $r=1$ and $l=1$, the error estimates above for both phase field, velocity and magnetic induction field are optimal in the sense of interpolation. Meanwhile, for the Taylor-Hood type element matching the N\'{e}d\'{e}lec edge element of the first-kind, namely, $r>2$ and $l>3$, we can also obtain the optimal error estimates. For $r=2$ and $l=2$, the error estimate for the velocity is  one order lower than the interpolation theory, which can be confirmed by numerical results.

$\bullet$ Furthermore, we enhance the accuracy of the velocity and magnetic induction field in $\H^{1}$-norm and $\H({\rm curl})$-norm, respectively, as follows:
\begin{equation}\label{temporal2}
\left\{
\begin{aligned}
\bigg(\Delta t\sum\limits_{k=0}^{N_T-1} 
	\|\nabla (\u^{k+1}-\u_{h}^{k+1} )\|^{2}  
	\bigg)^{\frac{1}{2}}
&\leq C\Big(\Delta t+h^{l+1}+\beta_{h}^{\star}\Big),	\\
\bigg(\Delta t\sum\limits_{k=0}^{N_T-1}
	\|\nabla\times(\B^{k+1}-\B_{h}^{k+1} )\|^{2} \bigg)^{\frac{1}{2}}
&\leq C\Big(\Delta t+h^{l }+\beta_{h}^{\star}\Big),
\end{aligned}
\right.
\quad  
\beta_{h}^{\star}= 
\left\{
\begin{array}{ll}
	h^{r+1},\,  &r\geq2,	\vspace{0.04in}\\
	h^{r+1},\,   &r=1 \ \ \big( (\u_{h}^{k+1},p_{h}^{k+1})\in \X_{h}^{2}\times \mathring{S}_{h}^{1} \big),	\vspace{0.04in}\\
	h^{r},\, 	& r=1 \ \ \big( (\u_{h}^{k+1},p_{h}^{k+1})\in \X_{h}^{1b}\times \mathring{S}_{h}^{1} \big).
\end{array}
\right.
\end{equation}
For $r\geq1$,  $l\geq1$, we conclude that the above error estimates for the velocity field and magnetic induction field are optimal in the sense of interpolation. Clearly, the temporal convergence order is $\mathcal{O}(\Delta t)$ in \eqref{temporal2}, which has been enhanced compared with the previous order of  $\mathcal{O}(\frac{\Delta t}{h})$ in \eqref{BL2}-\eqref{vb}. Thus, the optimal $\H^{1}$-norm error estimates for velocity fields and the optimal $\H(\rm{curl})$-norm error estimates for magnetic induction fields are obtained. The detailed rigorous analysis will be presented in Theorem \ref{theorem2-2}.
 
The rest of this work is organized as follows: In Section \ref{sec-main}, we propose the fully discrete finite element scheme as well as the main results. In Section \ref{sec-projections}, we introduce the Ritz, Stokes and Maxwell quasi-projections, and   in Section \ref{sec-proof}, we give the proof of the main Theorem \ref{theorem2-2}.  Numerical examples are conducted in Section \ref{sec-examples} to confirm our theoretical analysis and demonstrate the efficiency of the method. The concluding remarks are summarized in Section \ref{sec-conclusion}.
 
\section{Numerical discretization and main results}\label{sec-main}

In this section, we design a fully discrete convex-splitting  finite element scheme for the CH-MHD model \eqref{PDE1}-\eqref{PDE5} and then present the main convergence results.
 
\subsection{Variational formulation}
Let $W^{k,p}(\Omega)$ denote the standard Sobolev spaces equipped with the norms $\|\cdot\|_{W^{k,p}}$ for integer $k\geq0$ and $p\in[1, \infty]$. As general, we denote by $H^{k}(\Omega)=W^{k,2}(\Omega)$ and $L^{p}(\Omega)=W^{0,p}(\Omega)$, and the vector-valued spaces by $\W^{k,p}(\Omega) = [W^{k,p}(\Omega)]^d$. The $L^2/\L^2$ inner product and norm are denoted by $(\cdot, \cdot)$ and $\|\cdot\|$, respectively.
The admissible spaces are defined as follows:
\begin{align*}
&\H_{0}^{1}(\Omega) = 
	\big\{\v\in  \H^{1}(\Omega): \v|_{\partial\Omega}=\0\, \big\},
&&L_{0}^{2}(\Omega) = 
	\big\{q\in L^{2}(\Omega): (q, 1)=0\, \big\},	\\
&\H({\rm curl}, \Omega) = 
	\big\{\C\in \L^{2}(\Omega): \nabla\times \C\in \L^{2}(\Omega) \, \big\}, 
&&\H_{0}({\rm curl}, \Omega) = 
	\big\{\C\in \H({\rm curl}, \Omega): \boldsymbol{n}\times \C|_{\partial\Omega}=\0\, \big\},	
\end{align*}
Then, it is natural to obtain that the exact solution of the CH-MHD system \eqref{PDE1}-\eqref{PDE5} satisfies the following variational formulation: for any test function $(\varphi, \psi, \v, q, \C)\in (H^{1}(\Omega), H^{1}(\Omega), \H_{0}^{1}(\Omega), L_{0}^{2}(\Omega), \H_{0}( {\rm curl}, \Omega) )$ and $t\in(0, T]$, it holds
\begin{align}
	&(\phi_{t}, \varphi) + (\nabla \phi \cdot\u, \varphi) + (\nabla \omega, \nabla \varphi) =  0,  \label{weak1} \\
	&(\nabla\phi, \nabla \psi) + (\phi^{3}-\phi, \psi) - (\omega, \psi) =0,   \label{weak2}\\
	&(\u_{t}, \v)+b(\u, \u,\v) + (\nabla\u, \nabla\v) - (\nabla\cdot \v, p) 
			+ (\nabla\times\B, \v\times\B) - (\omega \nabla \phi, \v) = 0,  \label{weak3}\\
	&(\nabla \cdot \u, q)=0, \label{weak4} \\
	&(\B_{t}, \C)+ (\nabla\times\B, \nabla\times \C) - (\u\times\B, \nabla\times\C) =0,	\label{weak5}
\end{align}
where we define the trilinear form $b(\u, \v,\w)=\frac12[( \u\cdot \nabla \v, \w)- ( \u\cdot \nabla  \w, \v) ]$, which is anti-symmetric respect to the two last arguments.

\begin{lemma}[\cite{2019A}]
Let ($\phi, \omega, \u, p, \B$) be the solution of the  two-phase MHD model \eqref{PDE1}-\eqref{PDE5}. Then, for any $t\in(0,T]$, the mass is conserved
	\begin{equation*}
		\big(\phi(t), 1\big) = \big(\phi_0, 1\big)
	\end{equation*}
	and the system is energy-stable,
	\begin{equation*}
		\frac{dE(\phi, \u, \B)}{dt} =-\Big(\lambda\varepsilon\|\sqrt{M(\phi)} \nabla\omega\|^{2}+ \|\sqrt{\nu (\phi)}\nabla\u\|^{2}+ \frac{1}{\mu} \|\frac{1}{\sqrt{\sigma (\phi)}} \nabla\times \B\|^{2} \Big)\leq 0,
	\end{equation*}
	where the total energy  is given by
	\begin{equation}\label{system-energy}
		E(\phi, \u, \B)
		= \frac{\lambda\gamma}{2}\|\nabla\phi\|^{2} 
		+\frac{\lambda}{4\gamma}\|\phi^{2}-1\|^2
		+\frac{1}{2}\|\u\|^{2}+\frac{1}{2\mu}\|\B\|^{2}.
	\end{equation} 
\end{lemma}

Hereafter, for brevity, we consider the physical parameters $\gamma=M(\phi)=\nu(\phi)=\mu=\lambda=\sigma (\phi)$=1, and note that for any positive constants there are no more essential difficulties. Moreover, we denote by $C $ a generic positive constant independent of $\Delta t$ and $h$, which may take different values at different places.

\subsection{Numerical scheme}
Let $\mathfrak{T}_{h}$ be a quasi-uniform partition of domain $\Omega$ into the simplices $K_{j}$ 
with mesh size $h=\max_j {\rm diam} (K_{j})$. We define the following finite element spaces, 
\begin{equation*}
\begin{array}{l l}
S_{h}^{r}=\{\phi_{h}\in C^0(\overline{\Omega}): 
	\, \phi_{h}|_{K_{j}}\in P_{r}(K_{j}), \forall K_{j} \in \mathfrak{T}_{h}\},	
& \mathring{S}_{h}^{r}=S_{h}^{r}\cap L_{0}^{2}(\Omega),	
			\vspace{0.06in} \\
\X_{h}^{r+1}=\{\v_{h}\in \H_{0}^{1}(\Omega):
	 \v_{h}|_{K_{j}}\in \P_{r+1}(K_{j}), \forall K_{j} \in \mathfrak{T}_{h}\},
&\X_{h}^{1b} 
	= (S_{h}^{1}\oplus B_{d+1})^{d}\cap \H_{0}^{1}(\Omega),	\vspace{0.06in} \\
\Y_{h}^{l}=\{\C_{h}\in \H_{0}( {\rm{curl}}, \Omega): 
\C_{h}|_{K_{j}}\in \P_{l-1}(K_{j})\oplus D_{h}^{l}(K_{j}), \forall K_{j} \in \mathfrak{T}_{h}\}, 
& 	\vspace{0.06in} \\
\Z_{h}^{l}=\{\C_{h}\in\boldsymbol{Z}_{h}^{l}: (\C_{h}, \nabla w_{h})=0, \forall\, w_{h} \in H_{0}^{1}(\Omega)\cap P_{l}(K_{j}) \}. 
&
\end{array}
\end{equation*}
where $P_{r}(K_{j})$ is the polynomial space of total degree $r$ on $K_{j}$ with its vector-valued form $\P_{r}(K_{j}):=[P_{r}(K_{j})]^d$.  
Note that $\X_{h}^{r+1}\times\mathring{S}_{h}^{1}$ is the standard Taylor--Hood elements and $\X_{h}^{1b} \times \mathring{S}_{h}^{1}$ is the MINI element, where we denote by $B_{3}$ and $B_{4}$ the spaces of cubic bubbles and quartic bubbles, respectively. According to the classical finite element theory  \cite{1986FiniteG, 1991BF}, we have the following discrete inf-sup condition:
\begin{equation*}
	\begin{aligned}
		\inf_{0\neq q_{h}\in \, \mathring{S}_{h}^{r}/\mathring{S}_{h}^{1} } \,
		\sup_{\0\neq\v_{h}\in \, \X_{h}^{r+1}/\X_{h}^{1b} }
			\frac{(\nabla\cdot \v_{h}, q_{h})}{\|q_{h}\| \, \|\nabla \v_{h}\|}
		\geq \beta_{0},
	\end{aligned}
\end{equation*}
where $\beta_{0}$ is a positive constants depending only on $\Omega$.
Furthermore, $\Y_{h}^{l}$ is an $l$-th order N\'{e}d\'{e}lec element space of the first type, with $D_{h}^{l}(K_{j}) := \{\boldsymbol{p}(\x) \in [\tilde{P}_{l}(K_{j})]^{d}: \boldsymbol{p}(\x)\cdot \x=0, \forall \x \in K_{j} \}$, where $\tilde{P}_{l}(K_{j})$ is homogeneous polynomial subsets of $P_{l}(K_{j})$.
For the sake of presentation, we use the following notations
\begin{equation*}
\mathring{\boldsymbol \chi }_{h}^{r}:=\left\{
\begin{array}{l l}
S_{h}^{r}\times S_{h}^{r}\times \X_{h}^{r+1} \times \mathring{S}_{h}^{r}\times \Y_{Bh}^{l }, 
	\quad 		&r\geq2, \; l\geq 3, \vspace{0.05in}\\
S_{h}^{1}\times S_{h}^{1}\times \X_{h}^{1b} \times \mathring{S}_{h}^{1}\times \Y_{Bh}^{1},  
	\quad 		&r=1,\;  l=1,		\vspace{0.05in}\\
S_{h}^{1}\times S_{h}^{1}\times \X_{h}^{2} \times \mathring{S}_{h}^{1}\times \Y_{Bh}^{2}, 
	\quad 		&r=1,\;  l=2.
\end{array}
\right.
\end{equation*}

Let $\big\{t_k = k\Delta t\big\}_{k=0}^{N_T}$ denote a uniform partition of the temporal interval $[0, T]$ with the time stepsize $\Delta t={T}/{N_T}$ for any positive integer $N_T$. 
By denoting $f^{k}=f (\x, t_{k})$ and $d_{t}f^{k+1} = (f^{k+1} - f^{k})/{\Delta t}$, then we propose the fully discrete convex splitting finite element scheme for the CH-MHD system \eqref{PDE1}-\eqref{PDE5}: Find ($\phi_{h}^{k+1}, \omega_{h}^{k+1}, \u_{h}^{k+1}, p_{h}^{k+1}, \B_{h}^{k+1})\in \mathring{\boldsymbol \chi }_{h}^{r}$ such that it holds
\begin{align}
&\big(d_{t}\phi_{h}^{k+1}, \varphi_{h}\big) 
		+ \big(\nabla\phi_{h}^{k}\cdot \u_{h}^{k+1},  \varphi_{h}\big) 
		+ \big( \nabla\omega_{h}^{k+1}, \nabla\varphi_{h}\big)=0,		\label{scheme1}\\
&\big(\nabla\phi_{h}^{k+1}, \nabla\psi_{h}\big)
		+ \big((\phi_{h}^{k+1})^{3}-\phi_{h}^{k}, \psi_{h}\big)
		- \big(\omega_{h}^{k+1}, \psi_{h}\big) = 0, 			\label{scheme2}\\
&\big(d_{t}\u_{h}^{k+1},\v_{h}\big) 
		+ b \big(\u_{h}^{k}, \u_{h}^{k+1},\v_{h}\big)
		+ \big( \nabla\u_{h}^{k+1},\nabla\v_{h} \big)
		- \big( \nabla\cdot \v_{h},  p_{h}^{k+1}\big) 
		+ \big(\nabla\times \B_{h}^{k+1}, \v_{h}\times \B_{h}^{k}\big)
		- \big(\omega_{h}^{k+1} \nabla \phi_{h}^{k}, \v_{h} \big) =0,	\label{scheme3}\\
&\big(\nabla\cdot\u_{h}^{k+1}, q_{h}\big) =0,									\label{scheme4}\\
&\big(d_{t}\B_{h}^{k+1}, \C_{h}\big) 
		+ \big( \nabla\times \B_{h}^{k+1}, \nabla\times \C_{h} \big)
		- \big(\u_{h}^{k+1}\times \B_{h}^{k}, \nabla\times \C_{h}\big) = 0,	\label{scheme5}
\end{align}
for all $(\varphi_{h}, \psi_{h}, \v_{h}, q_{h},  \C_{h})\in \mathring{\boldsymbol \chi }_{h}^{r}$ and $k=0,1,\cdots, N_T-1$. 
The initial data is set as
\begin{align*}
	\phi_{h}^{0}=R_h\phi_{0},
		\qquad
	\u_{h}^{0}= \I_{h}\u_{0}, 
		\qquad
	\B_{h}^{0}=\PPi_h\B_{0},
\end{align*}
where $R_h$, $\I_{h}$, and $\PPi_h$ are the standard Ritz, $\L^{2} $, and Maxwell quasi-projection operators, respectively, and the definitions will be given in the next section.

The energy stability \eqref{scheme1}-\eqref{scheme5} has been presented in \cite{2019A}, and here we omit the proof for compactness.
\begin{lemma}[Theorem 4.1, \cite{2019A}]
The scheme \eqref{scheme1}-\eqref{scheme5} is mass-conserved
\begin{equation}\label{mass-conservetion}
	\big(\phi_h^{k+1}, 1\big) = \big(\phi_h^0, 1\big)
\end{equation}
and admits the following discrete energy decaying law
\begin{equation}\label{menergy}
E_h^{k+1} - E_h^k
	\leq -\Delta t\,\Big(\| \nabla\omega_{h}^{k+1}\|^{2}
	+ \|\nabla\u_{h}^{k+1} \|^{2}
	+ \| \nabla\times \B_{h}^{k+1}\|^{2} \Big),
\end{equation}
for $k=0,1,\cdots,N_T-1$, where the total energy is given by
\begin{equation}\label{algorithm-energy}
E_h^{k}
	= \frac{1}{2}\|\nabla\phi_{h}^{k}\|^{2}+\frac{1}{4}\|(\phi_{h}^{k})^{2}-1\|^2 +\frac{1}{2}\|\u_{h}^{k}\|^{2}+\frac{1}{2}\|\B_{h}^{k}\|^{2}.
\end{equation}
\end{lemma}


\subsection{Main results}
It is supposed that the unique solution of CH-MHD model \eqref{PDE1}-\eqref{PDE5} exists and satisfies the regularity assumption
\begin{equation}\label{regularity_ass}
	\begin{aligned}
		& \phi \in H^{2}(0,T; L^{2}(\Omega))\cap H^{1}(0,T; H^{r+1}(\Omega))\cap C(0,T; W^{2, 4}(\Omega)), \quad  \omega\in H^{1}(0,T; H^{r+1}(\Omega)),\\
		&\u\in \H^{2}(0,T; \L^{2}(\Omega))\cap \H^{1}(0,T; \H^{r+2}(\Omega)),\quad p\in L^{2}(0,T; H^{r+1}(\Omega)\cap L_{0}^{2}(\Omega)),\\
		&\B,\, \nabla\times \B \in \L^{\infty}(0,T; [\L^{\infty}(\Omega)\cap \H^{l}(\Omega)\cap \W^{1,3}(\Omega)]).
	\end{aligned}
\end{equation}
with $r\geq1$ and $l\geq1$.
Then, we obtain the main results of this work in the following theorem.

\begin{The}\label{theorem2-2}
Supposing that CH-MHD model \eqref{PDE1}-\eqref{PDE5} admits a unique solution $(\phi, \omega, \u, p, \B)$ satisfying the regularity assumptions \eqref{regularity_ass}, then the numerical solution $(\phi_{h}^{k}, \omega_{h}^{k}, \u_{h}^{k}, p_{h}^{k}, \B_{h}^{k}) \in \mathring{\boldsymbol \chi }_{h}^{r}$ of fully discrete scheme \eqref{scheme1}-\eqref{scheme5}, $k=1,\cdots, N_T$,  satisfies the following error estimates:
\begin{align*}
&\max\limits_{1\leq k\leq N_T}\|\phi^{k}-\phi_{h}^{k} \| + \bigg(\Delta t\sum_{k=1}^{N_T} 
	\|\omega^{k}-\omega_{h}^{k}\|^{2}\bigg)^{\frac{1}{2}} \leq C \Big(\Delta t+h^{l+1}+h^{r+1}\Big), \\
&\max\limits_{1\leq k\leq N_T}  \|\u^{k}-\u_{h}^{k} \|
    \leq C \Big(\Delta t+h^{l+1}+\beta_{h}\Big),\\
&\max\limits_{1\leq k\leq N_T} \|\B^{k}-\B_{h}^{k} \|
	\leq C \Big(\Delta t+h^{l }+\beta_{h}\Big),	\\
&\max\limits_{1\leq k\leq N_T}\|\nabla(\phi^{k}-\phi_{h}^{k}) \|\leq C \Big(\Delta t+h^{l+1}+h^{r}\Big),	\\
&\bigg(\Delta t\sum_{k=1}^{N_T}  \|\nabla (\u^{k}-\u_{h}^{k} )  \|^{2}  \bigg)^{\frac{1}{2}}
	\leq C \Big(\Delta t+h^{l+1}+\beta_{h}^{\star}\Big),\\
&\bigg(\Delta t\sum_{k=1}^{N_T} \|\nabla\times(\B^{k}-\B_{h}^{k} )\|^{2} \bigg)^{\frac{1}{2}}\leq C \Big(\Delta t+h^{l }+\beta_{h}^{\star}\Big),
\end{align*}
with
\begin{equation*}
\beta_{h}=\left\{
\begin{aligned}
& h^{r+2},\quad r\geq2,\\
& h^{r+1},\quad  r=1,
\end{aligned}
\right.\qquad
\beta_{h}^{\star}=\left\{
\begin{array}{ll}
  	h^{r+1},\quad  &r\geq2,	\\
 	h^{r+1},\quad  &r=1 \quad ( (\u_{h}^{k},p_{h}^{k})\in \X_{h}^{2}\times \mathring{S}_{h}^{1} ),\\
 	h^{r},\quad    &r=1 \quad ( (\u_{h}^{k},p_{h}^{k})\in \X_{h}^{1b}\times \mathring{S}_{h}^{1} ).
\end{array}
\right.
\end{equation*}
\end{The}

\section{Projections and their properties}\label{sec-projections}

We first present the classic projection operators \cite{wheeler1973priori}:

(1) The $L^{2}/\L^2$ projection $I_{h}: L^{2}(\Omega)\rightarrow S_{h}^{r}$ and $\I_{h}: \L^{2}(\Omega)\rightarrow \X_{h}^{r+1}/\X_{h}^{1b}$ are defined as follows,
\begin{equation*}
\begin{aligned}
&(v -I_{h}v, \varphi_{h})=0,  \quad		&&\forall \varphi_{h}\in S_{h}^{r},\\
&(\v-\I_{h}\v, \v_{h})=0,	\quad	&&\forall \v_{h}\in \X_{h}^{r+1}/\X_{h}^{1b}.\\
\end{aligned}
\end{equation*}
For these projections, the following estimates hold
\begin{equation*}
\begin{aligned}
& \|v -I_{h}v \| +h \|\nabla (v -I_{h}v)\|\leq C  h^{r+1}\|v\|_{H^{r+1}},	\quad	&I_hv&\in S_h^r, \\
& \|\v-\I_{h}\v\| +h \|\nabla (\v-\I_{h}\v)\|\leq C  h^{r+2}\|\v\|_{H^{r+2}},\quad &\I_{h}\v&\in \X_{h}^{r+1},  \\
& \|\v-\I_{h}\v\| +h \|\nabla (\v -\I_{h}\v)\|\leq C  h^{2}\|\v\|_{H^{2}},\quad &\I_{h}\v &\in \X_{h}^{1b}.
\end{aligned}
\end{equation*}

(2) The classic Ritz projection $R_h: H^{1}(\Omega)\rightarrow S_{h}^{r}$ is defined by,
\begin{equation*}
(\nabla (\varphi-R_h\varphi), \nabla\psi_{h})=0,	\quad	\forall\, \psi_{h}\in S_{h}^{r},
\end{equation*}
with $\int_{\Omega}(\varphi-R_h\varphi) \mathrm{d}\x$=0 for the uniqueness. The Ritz projection satisfies with following estimates:
\begin{equation*}
\begin{aligned}
&\|\varphi-R_h\varphi\|_{L^{s}}+h\|\varphi-R_h\varphi\|_{W^{1, s}}\leq C h^{r+1}\|\varphi\|_{W^{r+1, s}},\\
&\|\varphi-R_h\varphi\|_{H^{-1}}\leq C \beta_{h}\|\varphi\|_{H^{r+1}},\\
&\|d_{t}(\varphi^{k}-R_h\varphi^{k})\|+h\|d_{t}(\varphi^{k}-R_h\varphi^{k})\|_{H^{1}}\leq C h^{r+1}\|d_{t}\varphi^{k}\|_{H^{r+1}},\\
&\|d_{t}(\varphi^{k}-R_h\varphi^{k})\|_{H^{-1}}\leq C \beta_{h}\|d_{t}\varphi^{k}\|_{H^{r+1}},
\end{aligned}
\end{equation*}
for $s\in[2, \infty]$ and $k=1, 2, \cdots, N_T$.
  
However, suffering from the high coupling nonlinearity, the previous works failed to obtain optimal error estimate in $L^{2}$-norm due to the limitations of the traditional approach. Thanks to the introduction of the Ritz quasi-projection, Stokes quasi-projection in \cite{2023Optimalwang} and the Maxwell quasi-projection in \cite{2023New}, we manage to improve the theoretical results.
 
(3) The Ritz quasi-projection $\widetilde{R}_h: H^{1}(\Omega)\rightarrow S_{h}^{r}$ is defined by
\begin{equation*}
(\nabla (\omega-\widetilde{R}_h\omega), \nabla\varphi_{h})+  (\nabla (\phi-R_h\phi)\cdot \u,   \varphi_{h})=0,
\end{equation*}
for all $\varphi_{h}\in S_{h}^{r}$ and $\int_{\Omega}(\omega-\widetilde{R}_h\omega) \mathrm{d}\x$=0, and its estimates are as follows:
\begin{equation*}
\begin{aligned}
&\|\omega-\widetilde{R}_h\omega\| +h\|\nabla(\omega-\widetilde{R}_h\omega)\|\leq C h^{r+1}\big(\|\u\|_{L^{\infty}}\|\phi\|_{H^{r+1}}+ \|\omega\|_{H^{r+1}}\big),\\
&\|\omega-\widetilde{R}_h\omega\|_{H^{-1}}\leq C \beta_{h}\big(\|\u\|_{W^{1, 4}}\|\phi\|_{H^{r+1}}+ \|\omega\|_{H^{r+1}}\big),\\
&\|\nabla(d_{t}(\omega^{k+1}-\widetilde{R}_h\omega^{k+1}) )\|\leq C h^{r}\big(  \|\u^{k+1}\|_{L^{\infty}}\|d_{t}\phi^{k+1}\|_{H^{r}}+ \|d_{t}\u^{k+1}\|_{L^{\infty}}\|\phi^{k}\|_{H^{r}} +\|d_{t}\omega^{k+1}\|_{H^{r+1}}   \big),\\
&\|d_{t}(\omega^{k+1}-\widetilde{R}_h\omega^{k+1})\|_{H^{-1}}\leq C \beta_{h}\big(  \|\u^{k+1}\|_{W^{1,4}}\|d_{t}\phi^{k+1}\|_{H^{r+1}}+ \|d_{t}\u^{k+1}\|_{W^{1, 4}}\|\phi^{k}\|_{H^{r+1}} +\|d_{t}\omega^{k+1}\|_{H^{r+1}}   \big),
\end{aligned}
\end{equation*}
for $k=0, 1, 2, \cdots, N_T-1$.

(4) The Stokes quasi-projection $(\P_{h}, P_{h}): \H_{0}^{1}(\Omega)\times L_{0}^{2}(\Omega)\rightarrow \X_{h}^{r+1}\times \mathring{S}_{h}^{r}/ \X_{h}^{1b}\times \mathring{S}_{h}^{1}$ is given by
\begin{equation*}
\begin{aligned}
&\big(\nabla (\u-\P_{h}(\u, p)), \nabla \v_{h}\big) - \big(p-P_{h}(\u, p), \nabla\cdot \v_{h}  \big)=\big( \omega \nabla( \phi- R_h\phi), \v_{h} \big),\\
&\big(\nabla\cdot (\u-\P_{h}(\u, p)), q_{h}\big) = 0,
\end{aligned}
\end{equation*}
for all $(\v_{h}, q_{h})\in \X_{h}^{r+1}\times \mathring{S}_{h}^{r}/ \X_{h}^{1b}\times \mathring{S}_{h}^{1}$, and we denote $\P_{h}\u:=\P_{h}(\u, p)$ and $P_{h}p:=P_{h}(\u, p)$ for simplicity. Further, the Stokes quasi-projection has the following estimates: 
\begin{equation}\label{err-U}
\|\u-\P_{h}\u\|=\left\{
\begin{aligned}
&Ch^{r+2} \big(\|\u\|_{H^{r+2}}+  \|p\|_{H^{r+1}} +  \|\phi\|_{H^{r+1}}\|\omega\|_{H^{2}}\big), &&r\geq2,	\\
&Ch^{r+1} \big(\|\u\|_{H^{r+1}} +  \|p\|_{H^{r}} +  \|\phi\|_{H^{r+1}}\|\omega\|_{H^{2}}\big), &&r=1,
\end{aligned}
\right.
\end{equation}
\begin{equation}\label{err-UH1}
\|\nabla (\u-\P_{h}\boldsymbol{u})\|+\|p-P_{h}p\|=\left\{
\begin{aligned}
&Ch^{r+1} \big(\|\u\|_{H^{r+2}}+  \|p\|_{H^{r+1}} +  \|\phi\|_{H^{r+1}}\|\omega\|_{W^{1, 4}}\big), &&r\geq2,	\\
&Ch^{r+1} \big(\|\u\|_{H^{r+2}}+  \|p\|_{H^{r+1}} +  \|\phi\|_{H^{r+1}}\|\omega\|_{W^{1, 4}}),  &&r=1 \, (\X_{h}^{2}\times \mathring{S}_{h}^{1}\big),	\\
&Ch^{r} \big(\|\u\|_{H^{r+1}}+  \|p\|_{H^{r}} +  \|\phi\|_{H^{r}}\|\omega\|_{W^{1, 4}}\big),  &&r=1 \, (\X_{h}^{1b}\times \mathring{S}_{h}^{1}),
\end{aligned}
\right.
\end{equation}
\begin{equation*}\label{Stokes3}
\|d_{t} (\u^{k+1} -\P_{h}\u^{k+1} )\| =\left\{
\begin{aligned}
&Ch^{r+2} \big(\|d_{t}\u^{k+1} \|_{H^{r+2}}+  \|d_{t}p^{k+1} \|_{H^{r+1}} +  \|d_{t}\phi^{k+1} \|_{H^{r+1}}\|\omega^{k+1} \|_{H^{2}}+\|\phi^{k} \|_{H^{r+1}}\|d_{t}\omega^{k+1} \|_{H^{2}} \big), &&r\geq2, \\
&Ch^{r+1}\big(\|d_{t}\u^{k+1} \|_{H^{r+1}}+  \|d_{t}p^{k+1} \|_{H^{r}} +  \|d_{t}\phi^{k+1} \|_{H^{r+1}}\|\omega^{k+1} \|_{H^{2}}+\|\phi^{k} \|_{H^{r+1}}\|d_{t}\omega^{k+1} \|_{H^{2}} \big), &&r=1 ,
\end{aligned}
\right.
\end{equation*}



(5) The Maxwell quasi-projection $\PPi_h:[H^{1}(\Omega)]^{d}\rightarrow \Z_{h}^{l}$ is defined as (cf. \cite{2023New}),
\begin{equation}\label{def-Max}
\big(\nabla\times (\B-\PPi_h\B), \nabla\times \C_{h}\big) 
	+ \big(\nabla\times\C_{h}\times (\B-\PPi_h\B), \u \big)
	+ \big(7(|\u|^{2}+1)(\B-\PPi_h\B), \C_{h}\big)
	=0, \quad \forall \C_{h} \in \Z_{h}^{l}.
\end{equation}

The estimates of Maxwell quasi-projection are as follows: 
\begin{align}
& \|\B-\PPi_h\B\|_{H(\rm{curl})}   \leq C h^{l}\| \B \|_{H^{l}(\rm{curl})},	\quad
\|\B-\PPi_h\B\|_{L^{3}}   \leq C h^{l} (\| \B \|_{H^{l}(\rm{curl})}+   \| \B\|_{W^{1, 3}}  ), \label{err-curl}\\
&\|\B-\PPi_h\B\|_{(H^{1})'} +\|\nabla\times(\B-\PPi_h\B)\|_{(H^{1}_{0})'} \leq C h^{l+1}\| \B\|_{H^{l}(\rm{curl})}, \label{error-dual}\\
& \|d_{t}(\B-\PPi_h\B)\|_{H(\rm{curl})}   \leq C h^{l} (\|d_{t}\B\|_{L^{\infty}(0,T; H^{l}(\rm{curl}))} +\| \B\|_{H^{l}(\rm{curl})}  ), \label{error-dcurl} \\
& \|d_{t}(\B-\PPi_h\B)\|_{(H^{1})'} +\|\nabla\times d_{t}(\B-\PPi_h\B)\|_{(H^{1}_{0})'} \leq C h^{l+1} (\|d_{t}\B\|_{L^{\infty}(0,T; H^{l}(\rm{curl}))}    +\| \B\|_{H^{l}(\rm{curl})} ). \label{error-dH1}
\end{align}

\section{The Proof of Theorem \ref{theorem2-2} }\label{sec-proof}

In this section, we will provide the proof of Theorem \ref{theorem2-2}, and some essential inequality are introduced.

\begin{lemma}[\cite{heywood1990finite}]\label{Gronwall}
Let $\alpha_{k}, \beta_{k}, c_{k}, \gamma_{k}$ and $g_{0}$ be a sequence of nonnegative numbers for integers $k\geq 0$ such that
\begin{equation*}
\alpha_{k}+\Delta t\sum_{j=0}^{k}\beta_{j}\leq \Delta t\sum_{j=0}^{k}\gamma_{j}\alpha_{j}+\Delta t\sum_{j=0}^{k}c_{j}+g_{0}.
\end{equation*}
Assume that $\gamma_{j}\Delta t\leq 1$ for all $j$, and set $\sigma_{j}=(1-\gamma_{j}\Delta t)^{-1}$. Then, for all $k\geq 0$, we have
\begin{equation*}
\alpha_{k}+\Delta t\sum_{j=0}^{k}\beta_{j}\leq {\rm exp}\bigg(\Delta t\sum_{j=0}^{k}\sigma_{j} \gamma_{j} \bigg)\bigg(  \Delta t\sum_{j=0}^{k}c_{j}+g_{0}\bigg).
\end{equation*}
\end{lemma}

\begin{lemma}[\cite{2023New, hiptmair2002finite}]\label{lemmaita}
For any $\boldsymbol \eta \in \L^{2}(\Omega)$ and $\y_{h}\in \Z_{h}^{l}$, we have
\begin{align*}
|(\y_{h}, \boldsymbol \eta)|\leq C (h\|\boldsymbol \eta \|+ \|\boldsymbol \eta \|_{(H^{1})'} )\|\nabla\times \y_{h}\|.
\end{align*}

\end{lemma}

\begin{lemma}[\cite{Jean2006Mathematical, 2019A, 1975Sobolev, 1986On}]
	We have the following Poincar\'{e}  inequalities and embedding inequalities
	\begin{align}
		&\|\psi\|_{L^{q}}\leq C \|\psi\|_{H^{1}},
				\qquad 
			\forall \psi\in H^{1}(\Omega)\cap L^2_0(\Omega), \H_{0}^{1}(\Omega),
				\qquad 
			1\leq q \leq 6,\label{i2}\\
		& \|v_{h}\|_{W^{m,s}}\leq C h^{n-m+\frac{d}{s}-\frac{d}{q}}\|v_{h}\|_{W^{n,q}},
				\qquad 
			v_{h}\in S_{h}^{r}, \mathring{S}_{h}^{r}, \X_{h}^{r+1}, \X_{h}^{1b},
				\quad 
			0\leq n\leq m\leq1,\quad  1\leq q\leq s\leq\infty.	 \label{i4}
	\end{align}
\end{lemma}

\subsection{Error equations}
For simplicity, we introduce the following error functions
\begin{align*}
e_{\phi}^{k}:=R_h\phi^{k}-\phi^{k}_{h},	
	\quad e_{\omega}^{k}:=\widetilde{R}_h\omega^{k}-\omega^{k}_{h}, 
	\quad e_{p}^{k}:=P_{h}p^{k}-p^{k}_{h},
	\quad e_{\u}^{k}:=\P_{h}\u^{k}-\u^{k}_{h}, 
	\quad e_{\B}^{k}:=\PPi_h\B^{k}-\B^{k}_{h}.
\end{align*}
 
With the help of projection operators defined in the previous section, we subtract \eqref{weak1}-\eqref{weak5} at $t=t_{k+1}$ from \eqref{scheme1}-\eqref{scheme2} to get the following error equations for ($e_{\phi}^{k+1}, e_{\omega}^{k+1}, e_{\u}^{k+1}, e_{p}^{k+1}, e_{\B}^{k+1}$),
\begin{align}
 &\big(d_{t}e_{\phi}^{k+1}, \varphi_{h} \big) 
		+\big(\nabla e_{\omega}^{k+1}, \nabla  \varphi_{h}\big)
		=\big(d_{t}(R_h\phi^{k+1}- \phi^{k+1}  ), \varphi_{h}\big)
		 +\big(\nabla \phi_{h}^{k}\cdot\u_{h}^{k+1},   \varphi_{h}\big)
		-\big( \nabla R_h\phi^{k+1}\cdot\u^{k+1},   \varphi_{h}\big)
		+\big(R_{1}^{k+1}, \varphi_{h}\big), 	\label{error1}\\
&\big( \nabla e_{\phi}^{k+1}, \nabla\psi_{h}  \big)
		-\big(e_{\omega}^{k+1}, \psi_{h} \big)
		+\frac{1}{2} \big( Z^{k+1}(e_{\phi}^{k+1}-e_{\phi}^{k}  ),  \psi_{h} \big)  
		=  -\big(Z^{k+1}\bar{e}_{\phi}^{k+\frac{1}{2}},  \psi_{h}  \big)
		+ \big(e_{\phi}^{k},  \psi_{h} \big)
		- \big((\phi^{k+1})^{3}-(R_h \phi^{k+1})^{3} , \psi_{h} \big)	\nonumber\\
&\hspace{0.5in}
		+ \big(\phi^{k}-R_h\phi^{k}, \psi_{h} \big)
		+ \big(\omega^{k+1}-\widetilde{R}_h\omega^{k+1}, \psi_{h} \big)
		+ \big(R_{2}^{k+1} ,\psi_{h} \big), 	\label{error2}\\
&\big(d_{t}e_{\u}^{k+1}, \v_{h} \big)
		+ \big(\nabla  e_{\u}^{k+1}, \nabla\v_{h} \big)
		- \big(\nabla\cdot \v_{h}, e_{p}^{k+1} \big)
		= \Big[ b (\u_{h}^{k}, \u_{h}^{k+1}, \v_{h} ) -b (\u^{k}, \u^{k+1}, \v_{h})\Big] \nonumber  \\
&\hspace{0.5in} 
		+\Big[ \big(\nabla R_h\phi^{k+1}\cdot \v_{h}, \omega^{k+1}\big) 
		- (\nabla \phi_{h}^{k}\cdot \v_{h}, \omega_{h}^{k+1}) \Big]	
	 	+ \big(d_{t}(\P_{h}\u^{k+1}-\u^{k+1}), \v_{h}\big) \nonumber\\
&\hspace{0.5in}  
		+\big(\nabla\times \B_{h}^{k+1}, \v_{h} \times\B_{h}^{k}\big) 
		-\big(\nabla\times \B^{k+1}, \v_{h}\times\B^{k}\big)+ \big( R_{3}^{k+1} ,\v_{h}  \big),  \label{error3}\\
&\big(\nabla\cdot e_{\u}^{k+1}, q_{h} \big)=0, 	\label{error4}\\
&\big(d_{t} e_{\B}^{k+1}, \C_{h}\big)
		+ \big(\nabla\times e_{\B}^{k+1},  \nabla\times  \C_{h}  \big) 
		=\big(d_{t}( \PPi_h\B^{k+1}-\B^{k+1}), \C_{h}\big)
		+\big(\nabla\times (\PPi_h\B^{k+1}-\B^{k+1}), \nabla\times\C_{h} \big)	\nonumber \\
&\hspace{0.5in}  
		-\big( (\u_{h}^{k+1}\times \B_{h}^{k}, \nabla\times \C_{h}\big)
		-\big(\u^{k+1}\times \B^{k}, \nabla\times \C_{h})\big)+\big(R_{4}^{k+1} ,\C_{h}\big)  ,		 \label{errorBT}
\end{align}
for any $(\varphi_{h}, \psi_{h}, \v_{h}, q_{h},  \C_{h})\in \mathring{\boldsymbol \chi }_{h}^{r}$  and $k=0, 1, \cdots, N_T-1$, where we define
\begin{equation*}
\begin{aligned}
&e_{\phi}^{k+1}=\bar{e}_{\phi}^{k+\frac{1}{2}}+\frac{1}{2}(e_{\phi}^{k+1}-e_{\phi}^{k}),\quad \bar{e}_{\phi}^{k+\frac{1}{2}}:=\frac{1}{2}(e_{\phi}^{k+1}+e_{\phi}^{k} ),\\
&(R_h\phi^{k+1})^{3}-(\phi_{h}^{k+1})^{3}=3e_{\phi}^{k+1}\int_{0}^{1}\left((1-\theta)\phi_{h}^{k+1} +\theta R_h \phi^{k+1}   \right)^{2}d\theta= \frac{1}{2}(e_{\phi}^{k+1}-e_{\phi}^{k} )Z^{k+1}+\bar{e}_{\phi}^{k+\frac{1}{2}} Z^{k+1},\\
&Z^{k+1}:=3\int_{0}^{1}\left((1-\theta)\phi_{h}^{k+1} +\theta R_h \phi^{k+1}   \right)^{2}d\theta.
\end{aligned}
\end{equation*}
In addition, $R_{1}^{k+1}, R_{2}^{k+1}, R_{3}^{k+1}$, and $R_{4}^{k+1}$ are the truncation terms satisfying
\begin{equation*}
\begin{aligned}
\big(R_{1}^{k+1}, \varphi_{h}\big)= &\big(d_{t}\phi^{k+1} - \phi_t^{k+1}, \varphi_{h}\big),	\qquad
\big(R_{2}^{k+1}, \psi_{h}\big)=  \big(\phi^k, \psi_{h}\big) -\big(\phi^{k+1}, \psi_{h}\big),	\\
\big(R_{3}^{k+1}, \v_{h}\big) = &\big(d_{t}\u^{k+1} - \u_{t}^{k+1}, \v_{h}\big)  + b\big(\u^{k} - \u^{k+1}, \u^{k+1}, \v_{h}\big)  + \big(\nabla\times \B^{k+1}, \v_{h}\times(\B^k-\B^{k+1}) \big),  \\
\big(R_{4}^{k+1}, \C_{h}\big)= &\big(d_{t}\B^{k+1} - \B_{t}^{k+1}, \C_{h}\big) - \big(\u^{k+1}\times (\B^{k}-\B^{k+1}) , \nabla\times \C_{h}\big).
\end{aligned}
\end{equation*} 
By utilizing the Taylor expansion, it is easy to obtain the truncation error estimates 
\begin{equation*}
\bigg(\Delta t\sum_{k=0}^{N_T-1}(\|R_{1}^{k+1}\|^{2} + \|R_{2}^{k+1}\|^{2} + \|d_tR_{2}^{k+1}\|^{2} +\|R_{3}^{k+1}\|^{2}+\|R_{4}^{k+1}\|^{2}   ) \bigg)^{\frac{1}{2}}\leq C \Delta t.
\end{equation*}

Next, we give the error estimates of the numerical solutions in the following lemma, which will be repeatedly used in later analysis.
\begin{lemma}[(4.21) and (4.22), \cite{2023Optimalwang}]\label{lemmaphi}
By taking $\varphi_{h}=e_{\omega}^{k+1}$ and $\psi_{h}=d_{t}e_{\phi}^{k+1}$ in equations (\ref{error1})-(\ref{error2}), respectively, we can get the following estimate:
\begin{align}\label{phi}
\quad \|\nabla e_{\phi}^{k+1}\|^{2}+\Delta t\sum_{m=0}^{k}\|\nabla e_{\omega}^{m+1}\|^{2}\leq & C_{\epsilon}\Delta t\sum_{m=0}^{k}\Big(\|e_{\textbf{u}}^{m+1}\|^{2} + \|\nabla e_{\phi}^{m+1}\|^{2} + \|  e_{\phi}^{m}\|^{2}    \Big)+C_{\epsilon} \Big(\beta_{h}^{2}+\Delta t^{2}  \Big) +\varepsilon \Delta t \sum_{m=0}^{k} \|  e_{\omega}^{m+1}\|^{2}        \nonumber  \\
&   + C \Delta t \sum_{m=0}^{k}\|d_{t}Z^{m+1}\|_{L^{\frac{3}{2}}}\| e_{\phi}^{m}\|^{2}_{H^{1}}+ \big(2+\epsilon  \big)\|e_{\phi}^{k+1}\|^{2}.
\end{align}
By additionally taking $\varphi_{h}=(-\Delta_{h})^{-1}e_{\phi}^{k+1}$
and $\psi_{h}=e_{\phi}^{k+1}-\frac{1}{|\Omega|}( e_{\phi}^{k+1}, 1)$ in equations (\ref{error1})-(\ref{error2}), there exists a small positive constant $\Delta t_{1}$ such that for $\Delta t \leq\Delta t_{1}$,  the $L^{2}$-norm estimates for $e_{\phi}^{k+1}$ and $e_{\omega}^{k+1}$ are as follows:
\begin{align*}
&\quad \|e_{\omega}^{k+1}\|\leq C \Big(\|\nabla e_{\omega}^{k+1}\|+  \|\nabla e_{\phi}^{k+1}\| +\|\nabla e_{\phi}^{k}\| +\beta_{h}+\Delta t \Big),	\\
&\quad \| e_{\phi}^{k+1}\|^{2}\leq \varepsilon \|\nabla e_{\phi}^{k+1}\|^{2}+C_{\varepsilon}\Delta t\sum_{m=0}^{k}\|e_{\u}^{k+1}\|^{2}+ C_{\varepsilon} (\beta_{h}^{2}+\Delta t^{2}),
\end{align*}
where $\beta_{h}$ is given in Theorem \ref{theorem2-2}. Here, the discrete Laplacian operator $\Delta_{h}:\mathring{S}_{h}^{r}\rightarrow \mathring{S}_{h}^{r}$ is defined by
	\begin{equation*}
		(-\Delta_{h}\psi_{h}, \varphi_{h})=(\nabla\psi_{h}, \nabla\varphi_{h}),\quad \forall \psi_{h}, \varphi_{h}\in \mathring{S}_{h}^{r}.
	\end{equation*}

\end{lemma}

\subsection{Estimates for $e_{\textbf{u}}^{k+1}$. }

Taking $\v_{h}=e_{\u}^{k+1}$,  $q_{h}=e_{p}^{k+1}$ in equations (\ref{error3})-(\ref{error4}), we have
\begin{align}\label{error-u}
&\frac{1}{2}d_{t}\|e_{\u}^{k+1}\|^{2} + \frac{1}{2\Delta t}\|e_{\u}^{k+1}-e_{\u}^{k}\|^{2} + \|\nabla e_{\u}^{k+1}\|^{2} 	\nonumber	\\
&= 
	\Big[ b\big(\u_{h}^{k}, \u_{h}^{k+1}, e_{\u}^{k+1}\big)-b\big(\u^{k}, \u^{k+1}, e_{\u}^{k+1}\big) \Big]  
	+ \Big[ \big(\nabla R_h\phi^{k+1}\cdot e_{\u}^{k+1}, \omega^{k+1}\big) - \big(\nabla \phi_{h}^{k}\cdot e_{\u}^{k+1}, \omega_{h}^{k+1}\big)\Big]	 \nonumber \\
&\ \ \ \
		+ \Big[\big(d_{t}(\P_{h}\u^{k+1}-\u^{k+1}), e_{\u}^{k+1}\big) + \big( R_{3}^{k+1} ,e_{\u}^{k+1} \big) \Big]
		+ \Big[\big(\nabla\times \B_{h}^{k+1}, e_{\u}^{k+1} \times\B_{h}^{k}\big) - \big(\nabla\times \B^{k+1}, e_{\u}^{k+1}\times\B^{k}\big)\big] \nonumber	\\
&:=
	\sum_{i=1}^{4} I_{i}.
\end{align}
 
\begin{lemma}[Section 4.2, \cite{2023Optimalwang}]\label{lemmaI}
We have the following known results 
\begin{align*}
&I_{1} \leq\varepsilon\|\nabla e_{\u}^{k+1}\|^{2}+C_{\varepsilon}\Big(\beta_{h}^{2} +\|e_{\boldsymbol{u}}^{k}\|^{2} \Big),\\
&I_{2} \leq C_{\varepsilon}\Big(\|e_{\u}^{k+1}\|^{2} + \|\nabla e_{\phi}^{k}\|^{2} + \beta_{h}^{2}+\Delta t^{2} \Big)+  \varepsilon\Big( \| e_{\omega}^{k+1}\|^{2}+  \|\nabla e_{\u}^{k+1}\|^{2}\Big),\\
&I_{3} \leq C \Big( \beta_{h}^{2}+ \| e_{\u}^{k+1}\|^{2}+ \|R_{3}^{k+1}\|^{2}    \Big).
\end{align*}
\end{lemma}

Next, we can estimate $I_{4}$ as
\begin{align}\label{II}
I_{4}
=&\, \big(\nabla\times \B_{h}^{k+1}, e_{\u}^{k+1} \times\B_{h}^{k}\big) - \big(\nabla\times \B^{k+1}, e_{\u}^{k+1}\times\B^{k}\big)	\nonumber	\\
=&\,  \big(\nabla\times (\B_{h}^{k+1}-\PPi_h\B^{k+1}), e_{\u}^{k+1} \times\B_{h}^{k}\big) 
	+ \big(\nabla\times (\PPi_h\B^{k+1} - \B^{k+1}), e_{\u}^{k+1} \times\B_{h}^{k}\big) \nonumber \\
& 
 	+ \big(\nabla\times \B^{k+1}, e_{\u}^{k+1} \times(\B_{h}^{k}-\PPi_h\B^{k})\big)	
	+ \big(\nabla\times \B^{k+1}, e_{\u}^{k+1} \times(\PPi_h\B^{k}-\B^k)\big)		\nonumber \\
\leq&\, \Big[\big(\nabla\times  e_{\B}^{k+1}, e_{\u}^{k+1}\times (\B^{k}-\B_{h}^{k}) \big) 
	+ C\|\nabla\times e_{\B}^{k+1}\|\,\|e_{\u}^{k+1}\|\,\|\B^k\|_{L^\infty}\Big]	\nonumber	\\
&
	+ \Big[Ch^{l-\frac{d}{6}}\,\|\nabla e_{\u}^{k+1}\|\,\|e_{\B}^k\|_{L^2}
	+ Ch^{2l} \,\|\nabla e_{\u}^{k+1}\| + Ch^{l+1}\|\nabla e_{\u}^{k+1}\|\Big] \nonumber	\\
&
	+ \|\nabla\times \B^{k+1}\|_{L^\infty}\| e_{\u}^{k+1}\|\,\|e_{\B}^k\|
	+ \|\nabla\times \B^{k+1}\|_{L^\infty}\|e_{\u}^{k+1}\|_{H^1}\|\PPi_h\B^{k}-\B^k\|_{(H_0^1)'}	\nonumber	\\
\leq&
	\,\big(\nabla\times  e_{\B}^{k+1}, e_{\u}^{k+1}\times (\B^{k}-\B_{h}^{k}) \big) 
	+ \varepsilon (\|\nabla e_{\u}^{k+1}\|^{2} + \|\nabla\times e_{\B}^{k+1}\|^{2}) + C_{\varepsilon} \big( \|e_{\u}^{k+1}\|^2 + \| e_{\B}^{k}\|^2 + h^{2l+2} \big),
\end{align} 
for some sufficiently small $h$, where we utilize
\begin{align*}
&\big(\nabla\times (\B_{h}^{k+1}-\PPi_h\B^{k+1}), e_{\u}^{k+1} \times\B_{h}^{k}\big) \\
&= 
	\big(\nabla\times(\B_{h}^{k+1}-\PPi_h\B^{k+1}), e_{\u}^{k+1}\times(\B_{h}^{k}-\B^k)\big)	
	+  \big(\nabla\times (\B_{h}^{k+1}-\PPi_h\B^{k+1}), e_{\u}^{k+1} \times\B^{k}\big)	\\
&\leq 
	\big(\nabla\times  e_{\B}^{k+1}, e_{\u}^{k+1}\times (\B^{k}-\B_{h}^{k}) \big) 
	+ C\|\nabla\times e_{\B}^{k+1}\|\,\|e_{\u}^{k+1}\|\,\|\B^k\|_{L^\infty},
\end{align*}
and
\begin{align*}
&\big(\nabla\times (\PPi_h\B^{k+1} - \B^{k+1}), e_{\u}^{k+1} \times\B_{h}^{k}\big)	\\
&=
	\big(\nabla\times (\PPi_h\B^{k+1} - \B^{k+1}), e_{\u}^{k+1} \times(\B_h^{k}-\PPi_h\B^k)\big)
	+ \big(\nabla\times (\PPi_h\B^{k+1} - \B^{k+1}), e_{\u}^{k+1} \times(\PPi_h\B^k-\B^k)\big) \\
&\quad\: \big(\nabla\times (\PPi_h\B^{k+1} - \B^{k+1}), e_{\u}^{k+1} \times\B^k\big)	\\
&\leq 
	\|\nabla\times (\PPi_h\B^{k+1} - \B^{k+1})\|\,\|e_{\u}^{k+1}\|_{L^6}\|e_{\B}^k\|_{L^3}
	+ \||\nabla\times (\PPi_h\B^{k+1} - \B^{k+1})\|\,\|e_{\u}^{k+1}\|_{L^6}\|\PPi_h\B^k-\B^k\|_{L^3} \\
&\quad +\|\nabla\times (\PPi_h\B^{k+1} - \B^{k+1})\|_{(H_0^1)'}\|e_{\u}^{k+1}\|_{H^1}\|\B^k\|_{L^\infty}	\\
&\leq 
	Ch^l \,\|\nabla e_{\u}^{k+1}\|\, h^{-\frac{d}{6}}\|e_{\B}^k\|_{L^2}
	+ Ch^l \,\|\nabla e_{\u}^{k+1}\|\, h^l + Ch^{l+1}\|\nabla e_{\u}^{k+1}\|.
\end{align*}
Combining Lemma \ref{lemmaI} and choosing a sufficiently small $\varepsilon$, equation (\ref{error-u}) reduces to
\begin{align*}
d_{t}\|e_{\u}^{k+1}\|^{2} +\|\nabla e_{\u}^{k+1}\|^{2} \leq 
& C_{\varepsilon}\big(\|e_{\u}^{k+1}\|^{2} + \|e_{\u}^{k}\|^{2} + \|\nabla e_{\phi}^{k}\|^{2}
   +\beta_{h}^{2} +h^{2(l+1)}+\Delta t^{2} +\|e_{\B}^{k}\|^{2}\big) \\
& +\varepsilon \big(\|e_{\omega}^{k+1}\|^{2}+\|\nabla\times  e_{\B}^{k}\|^{2}\big) + \big(\nabla\times  e_{\B}^{k+1}, e_{\u}^{k+1}\times (\B^{k}-\B_{h}^{k}) \big).
\end{align*}
Consequently, summing up the above estimate from time step $t_{0}$ to $t_{k}$ leads to
\begin{align}\label{Nu}
\|e_{\u}^{k+1}\|^{2} + \Delta t\sum_{m=0}^{k}\|\nabla e_{\u}^{m+1}\|^{2}\leq 
&C_{\varepsilon}\Delta t \sum_{m=0}^{k} 
  \big(\|e_{\u}^{m+1}\|^{2} + \|\nabla e_{\phi}^{m+1}\|^{2} \big)
  + C_{\varepsilon} \big(\beta_{h}^{2}+h^{2(l+1)} +\Delta t^{2}\big)  \nonumber \\
&+\varepsilon \Delta t \sum_{m=0}^{k} \big(\|e_{\omega}^{m+1}\|^{2} + \|\nabla\times  e_{\B}^{m+1}\|^{2}\big) + \Delta t \sum_{m=0}^{k}  \big(\nabla\times e_{\B}^{m+1},   e_{\u}^{m+1}\times (\B^{m}-\B_{h}^{m}) \big).
\end{align}
 
\subsection{Estimates for $e_{\textbf{B}}^{k+1}$.}
 
Taking $\C_{h}= e_{\B}^{k+1}$ in equation (\ref{errorBT}), we have
\begin{align}\label{error-BB}
&\frac{1}{2}d_{t}\| e_{\B}^{k+1}\|^{2} + \frac{1}{2\Delta t}\|e_{\B}^{k+1}-e_{\B}^{k}\|^{2} + \|\nabla \times e_{\B}^{k+1}\|^{2} 
		= \big(d_{t}(  \PPi_h\B^{k+1}-\B^{k+1}), e_{\B}^{k+1} \big)  
		+ \big(R_{4}^{k+1} ,e_{\B}^{k+1}\big) 	 \nonumber \\
&\ \ \ \ 
	+ \Big[ \big(\nabla\times (\PPi_h\B^{k+1}-\B^{k+1}), \nabla\times e_{\B}^{k+1}\big) 
		 - (\u_{h}^{k+1}\times \B_{h}^{k}, \nabla\times e_{\B}^{k+1}) +  (\u^{k+1}\times \B^{k}, \nabla\times e_{\B}^{k+1}) \Big] 	\nonumber	\\
&:= \sum_{i=1}^{3}Q_i.
\end{align}
According to Lemma \ref{lemmaita} and \eqref{error-dcurl}-\eqref{error-dH1}, we obtain the estimates of $Q_1$
\begin{align}\label{ets_Q1}
	Q_{1} 
	&\leq  \Big(h\|d_{t}(\PPi_h\B^{k+1}-\B^{k+1})\| + \|d_{t}(  \PPi_h\B^{k+1}-\B^{k+1})\|_{(H^{1})'} \Big) \|\nabla\times e_{\B}^{k+1}\|  	\nonumber	\\
	&\leq C_{\varepsilon}\Big(h^{2(l+1)}+  \Delta t^{2} \Big) + \varepsilon\|\nabla\times e_{\B}^{k+1}\|  ^2,
\end{align}
and by using the Taylor expansion, we can estimate $Q_2$ as 
\begin{align}\label{ets_Q2}
	Q_2 \leq \|R_{4}^{k+1}\|\, \|\nabla\times e_{\B}^{k+1}\|  
			\leq C_{\varepsilon}\Delta t^{2} + \varepsilon\|\nabla\times e_{\B}^{k+1}\|^{2}.
\end{align}
Employing the definition of the Maxwell quasi-projection \eqref{def-Max}, we know
\begin{align*}\label{QQ3}
Q_{3}
= &
	\big(\nabla\times e_{\B}^{k+1}\times(\B^{k+1}-\PPi_h\B^{k+1}), \u^{k+1}\big)
	+ \big(7(|\u^{k+1}|^2+1)(\B^{k+1}-\PPi_h\B^{k+1}), e_{\B}^{k+1}\big)	\nonumber	\\
& 
	+ \big(\nabla\times e_{\B}^{k+1}\times\B_h^{k}, \u_h^{k+1}\big) 
	- \big(\nabla\times e_{\B}^{k+1}\times\B^k, \u^{k+1}\big)	\nonumber\\
=&
	\big(\nabla\times e_{\B}^{k+1}\times[(\B^{k+1}-\PPi_h\B^{k+1}) - (\B^k-\PPi_h\B^k)], \u^{k+1}\big)
	+ \big(\nabla\times e_{\B}^{k+1}\times(\B^k-\PPi_h\B^k), \u^{k+1}\big)	\nonumber	\\
&
	+ \big(7(|\u^{k+1}|^2+1)(\B^{k+1}-\PPi_h\B^{k+1}), e_{\B}^{k+1}\big)	
	+ \big(\nabla\times e_{\B}^{k+1}\times(\B_h^{k}-\B^k), \u_h^{k+1}\big)		\nonumber	\\
&
	+ \big(\nabla\times e_{\B}^{k+1}\times\B^k, \u_h^{k+1}-P_h\u^{k+1}\big)	
	+ \big(\nabla\times e_{\B}^{k+1}\times\B^k, P_h\u^{k+1}-\u^{k+1}\big)			\\
:=&
	\sum_{i=1}^{6}Q_{3,i}
\end{align*}
Applying the estimates of quasi-projections and Lemma \ref{lemmaita},  we can estimate $Q_{3,1}-Q_{3,6}$ as
\begin{align*}
	Q_{3,1} &\leq \|\nabla\times e_{\B}^{k+1}\|\, \Delta t\|d_t(\B^{k+1}-\PPi_h\B^{k+1})\|\,\|\u^{k+1}\|_{L^\infty}
	\leq C_{\varepsilon}\big(\Delta t^2 + h^{4l}\big) + \varepsilon\|\nabla\times e_{\B}^{k+1}\|^2,		\\
	Q_{3,2} + Q_{3,4} 
		&= \big(\nabla\times e_{\B}^{k+1}\times(\B^k-\PPi_h\B^k), \u^{k+1}\big)
			+ \big(\nabla\times e_{\B}^{k+1}\times(\B_h^{k}-\B^k), \u_h^{k+1}\big)		\\
		& = \big(\nabla\times e_{\B}^{k+1}\times(\B^k-\B_h^k), \u^{k+1}\big)
			+ \big(\nabla\times e_{\B}^{k+1}\times(\B_h^k-\PPi_h\B^k), \u^{k+1}\big) 	\\
		&\quad
			+ \big(\nabla\times e_{\B}^{k+1}\times(\B_h^{k}-\B^k), \u_h^{k+1}-\P_h\u^{k+1}\big)	
			+ \big(\nabla\times e_{\B}^{k+1}\times(\B_h^{k}-\B^k), \P_h\u^{k+1}-\u^{k+1}\big)	\\
		&\quad 
			+ \big(\nabla\times e_{\B}^{k+1}\times(\B_h^{k}-\B^k), \u^{k+1}\big)		\\
		&= -\big(\nabla\times e_{\B}^{k+1}\times e_{\B}^k, \u^{k+1}\big) 
			- \big(\nabla\times e_{\B}^{k+1}\times(\B_h^{k}-\B^k), e_{\u}^{k+1}\big)	 \\
		&\quad
			+ \big(\nabla\times e_{\B}^{k+1}\times(\B_h^{k}-\PPi_h\B^k), \P_h\u^{k+1}-\u^{k+1}\big)
			+ \big(\nabla\times e_{\B}^{k+1}\times(\PPi_h\B^k-\B^k), P_h\u^{k+1}-\u^{k+1}\big)	\\
		&\leq 
			\|\nabla\times e_{\B}^{k+1}\|\,\|e_{\B}^k\|\,\|\u^{k+1}\|_{L^\infty}	
			+ \big(\nabla\times e_{\B}^{k+1}\times(\B_h^{k}-\B^k), e_{\u}^{k+1}\big) \\
		&\quad
			+ \|\nabla\times e_{\B}^{k+1}\|\,\|e_{\B}^k\|_{L^3}\|\P_h\u^{k+1}-\u^{k+1}\|_{L^6}	
			+ \|\nabla\times e_{\B}^{k+1}\|\,\|\PPi_h\B^k-\B^k\|_{L^3}\|\P_h\u^{k+1}-\u^{k+1}\|_{L^6} \\
		&\leq 
			C \|\nabla\times e_{\B}^{k+1}\|\,\|e_{\B}^k\| 
			+ \big(\nabla\times e_{\B}^{k+1}\times(\B_h^{k}-\B^k), e_{\u}^{k+1}\big) 	
			+ C \|\nabla\times e_{\B}^{k+1}\|\, h^{-\frac{d}{6}}\|e_{\B}^k\| h^r	
			+ C  \|\nabla\times e_{\B}^{k+1}\| h^{l+r}	\\
		&\leq
			C_{\varepsilon}\big(\|e_{\B}^k\|^2 + h^{2l+2r}\big) + \varepsilon\|\nabla\times e_{\B}^{k+1}\|^2 + \big(\nabla\times e_{\B}^{k+1}\times(\B_h^{k}-\B^k), e_{\u}^{k+1}\big), 	\\
	Q_{3,3} &\leq \|7(|\u^{k+1}|^2+1)\|_{L^{\infty}}\big|\big((\B^{k+1}-\PPi_h\B^{k+1}), e_{\B}^{k+1}\big)\big|	\\
		&\leq C\big(h\|\B^{k+1}-\PPi_h\B^{k+1}\| + \|\B^{k+1}-\PPi_h\B^{k+1}\|_{(H^1)'}\big) \|\nabla\times e_{\B}^{k+1}\|	\\
		&\leq C_{\varepsilon}h^{2l+2} + \varepsilon\|\nabla\times e_{\B}^{k+1}\|^2,	\\
	Q_{3,5} &\leq \|\nabla\times e_{\B}^{k+1}\|\,\|\B^k\|_{L^\infty} \|e_{\u}^{k+1}\|
		\leq C_{\varepsilon}\|e_{\u}^{k+1}\|^2 + \varepsilon\|\nabla\times e_{\B}^{k+1}\|^2,	\\
	Q_{3,6} &\leq \|\nabla\times e_{\B}^{k+1}\|\,\|\B^k\|_{L^\infty}\|P_h\u^{k+1}-\u^{k+1}\|	
		\leq C_{\varepsilon}\beta_h^2 + \varepsilon\|\nabla\times e_{\B}^{k+1}\|^2,	
\end{align*}
Combining the above estimates, we have
\begin{align*}
Q_{3} \leq C_{\varepsilon}\Big(\Delta t^{2} + h^{2l+2} + \beta_{h}^{2}+ \|e_{\u}^{k+1} \| ^{2}+\|e_{\B}^{k }\|^{2} \Big) 
          + \varepsilon \|\nabla\times e_{\B}^{k+1}\|^{2}
          + \big(\nabla\times e_{\B}^{k+1}\times (\B^{k }- \B_{h}^{k }), e_{\u}^{k+1}  \big).
\end{align*}
Thus, using the estimate $Q_1$, $Q_2$ and $Q_3$, \eqref{error-BB} can be rewritten as
\begin{align*}
d_{t}\|e_{\B}^{k+1}\|^{2} +\|\nabla \times e_{\B}^{k+1}\|^{2}
	\leq C_{\varepsilon}\Big(\Delta t^{2}+ h^{2(l+1)}+ \beta_{h}^{2}+ \|e_{\u}^{k+1} \| ^{2}+\|e_{\B}^{k }\|^{2} \Big)  
	+ \big(\nabla\times e_{\B}^{k+1}\times (\B^{k}- \B_{h}^{k}), e_{\u}^{k+1}  \big),
\end{align*}
with a sufficiently small $\varepsilon$. Given $e_{\B}^{0}=\0$ and summing up the above inequality from time step $t_{0}$ to $t_{k}$, we obtain
\begin{align*}
  \|e_{\B}^{k+1}\|^{2} +\Delta t\sum_{m=0}^{k}\|\nabla \times e_{\B}^{m+1}\|^{2}\leq 
&  C_{\varepsilon}\Big(\Delta t^{2}+ h^{2(l+1)}+ \beta_{h}^{2}\Big)
  + C_{\varepsilon}\Delta t\sum_{m=0}^{k}\Big(\|e_{\u}^{m+1} \| ^{2}+\| e_{\B}^{m }\|^{2} \Big) \\
&+  \Delta t\sum_{m=0}^{k} \Big(\nabla\times e_{\B}^{m+1}\times (\B^{m}- \B_{h}^{m}), e_{\u}^{m+1}  \Big).
\end{align*}
which together with Lemma \ref{lemmaphi} and \eqref{Nu} finally leads to
\begin{align*}
&\|\nabla e_{\phi}^{k+1}\|^{2}+\|e_{\u}^{k+1}\|^{2}+\|e_{\B}^{k+1}\|^{2}+\Delta t\sum_{m=0}^{k}\Big(\|\nabla e_{\omega}^{m+1}\|^{2}+ \|\nabla e_{\u}^{m+1}\|^{2}+\|\nabla \times e_{\B}^{m+1}\|^{2}   \Big)\nonumber \\
\leq 
& C_{\epsilon}\Delta t\sum_{m=0}^{k}\Big(\|e_{\u}^{m+1}\|^{2} + \|\nabla e_{\phi}^{m+1}\|^{2}  +\| e_{\B}^{m }\|^{2}  \Big)
 + C_{\epsilon} \Big(\beta_{h}^{2}+h^{2(l+1)}+\Delta t^{2}  \Big)
 + C \Delta t \sum_{m=0}^{k}\|d_{t}Z^{m+1}\|_{L^{\frac{3}{2}}}\| e_{\phi}^{m}\|^{2}_{H^{1}}.
\end{align*} 
By using the discrete Gronwall's inequality in Lemma \ref{Gronwall}, there exists a positive constant $\Delta t_{2}$ such that for $\Delta t\leq \min\{\Delta t_{1}, \Delta t_{2}\}$, 
\begin{align}
\|\nabla e_{\phi}^{k+1}\|^{2}+\|e_{\u}^{k+1}\|^{2}+\|e_{\B}^{k+1}\|^{2}+\Delta t\sum_{m=0}^{k}\Big(\|\nabla e_{\omega}^{m+1}\|^{2}+ \|\nabla e_{\u}^{m+1}\|^{2}+\|\nabla \times e_{\B}^{m+1}\|^{2}   \Big)\leq C_{\epsilon} \Big(\beta_{h}^{2}+h^{2(l+1)}+\Delta t^{2}  \Big).
\end{align}
Combining the estimates of projection operators in Section \ref{sec-projections} and the triangle inequality, the proof of convergence results in Theorem \ref{theorem2-2} is completed.

\section{Numerical examples}\label{sec-examples}

In this section, we conduct several 2D/3D numerical examples to verify the theoretical analysis and demonstrate the performance of the scheme \eqref{scheme1}-\eqref{scheme5} based on the software FreeFem \cite{Hecht+2012+251+266}. In addition,  $l=1, 2$ mean the first-type  zero-order and first-order N\'{e}d\'{e}lec edge elements, which are marked by  $RT0Ortho$ and $RT1Ortho$. 
 The finite element spaces are chosen as follows:
\begin{table}[hpt]
\caption{Selection of finite element spaces}\label{case}
\centering
\begin{tabular}{ccccccccccccccccccccccccccccccc}
\toprule
 & $ \phi $ &    $ \u $   & $p  $ & $\B$ &	\\ 	 \hline
 case I  &  $P_{1}$   &   $P_{1}^{b}$    &  $P_{1}$   &   $RT0Ortho$	   & \\
 case II  & $P_{1}$   &  $P_{2}$    &  $P_{1}$   &   $RT1Ortho$    & \\
\bottomrule
\end{tabular}
\end{table}

\subsection{Convergence test}

Consider that in a $d$-dimension domain $\Omega=[0,1]^{d}$, CH-MHD system \eqref{PDE1}-\eqref{PDE5} admits the following smooth exact solutions:
\begin{equation*}
	\begin{aligned}
\mbox{in 2D case:} \,
&\left\{
\begin{aligned}
\phi	&=  \cos(t)\cos^2(\pi x)\cos^2(\pi y),\\
\u	&= \cos(t)(\pi\sin(2\pi y)\sin^2(\pi x) , - \pi \sin(2\pi x) \sin^2(\pi y) )^\top,\\
p	&= \cos(t)(2x-2) (2y -1),\\
\B	&= \cos(t)(\sin(\pi x)\cos(\pi y),  -\sin(\pi y)\cos(\pi x))^\top, \\
\end{aligned}
\right.			\\
	\mbox{and} \qquad \mbox{in 3D case:} \,
&\left\{
\begin{aligned}
\phi&= {\rm exp}(-2t)\sin^2(\pi x)\sin^2(\pi y)\sin^2(\pi z),\\
\u	&= {\rm exp}(t)(y(1-y)z(1-z) ,\ x(1-x)z(1-z),\ x(1-x)y(1-y) )^\top,\\
p	&= {\rm exp}(t)(2x-1) (2y -1) (2z-1),\\
\B	&= {\rm exp}(t)(\sin(\pi y)\sin(\pi z), \sin(\pi x)\sin(\pi z), \sin(\pi x)\sin(\pi y) )^\top. \\
\end{aligned}
\right.	
		\end{aligned}
\end{equation*}
We test the time and  space convergence rates at the terminal time $T=1$ with $\Delta t=\mathcal{O}(h^{2})$ for cases I-II in Table \ref{case}. The numerical results for case I in 2D and 3D are given in Table \ref{caseI-2D} and  Table \ref{caseI-3D}.   Table \ref{caseII-2D} shows the results of case II in 2D case. It is seen that the convergence orders are consistent with the theoretical results in Theorem \ref{theorem2-2}.
  
\begin{table}[htbp]
\caption{Convergence results  with case I in 2D.}\label{caseI-2D}
\centering
\begin{tabular}{cccccccccccc}
\toprule
$h$ & $ \|\phi^{N_T}-\phi_{h}^{N_T}\| $ & rate & $\|\nabla (\phi^{N_T}-\phi_{h}^{N_T})\|$ & rate & $\|\u^{N_T}-\u_{h}^{N_T}\|$ & rate & $\|\nabla (\u^{N_T}-\u_{h}^{N_T})\|$ & rate &\\ \hline
 1/8 &2.66e-01 & 1.48  &3.62e-01  &  1.20   &1.30e-01 & 1.66 &3.02e-01 & 0.92  & \\
 1/16 &7.48e-02 & 1.83 &1.57e-01  & 1.21  &3.44e-02 & 1.92   &1.52e-01 & 0.99  & \\
 1/32 &1.93e-02  &1.95  &7.36e-02 & 1.09  &8.67e-03 & 1.99 &7.58e-02 & 1.00  &\\
 1/64 &4.86e-03 & 1.99 &3.61e-02 & 1.03  &2.17e-03 & 2.00  &3.78e-02  & 1.00  &\\
 \bottomrule	
 \toprule
$h$ &$\|\B^{N_T}-\B_{h}^{N_T}\|$ & rate & $\| \B^{N_T}-\B_{h}^{N_T} \|_{H(\rm {curl})}$ & rate & $\|p^{N_T}-p_{h}^{N_T}\|$ & rate & \\ 
\hline
  1/8 &1.14e-01 & 1.03  &2.26e-01  &  0.99  & 3.28e-00 &1.45  & \\
  1/16 &5.68e-02 & 1.01   &1.13e-01 & 1.00  & 1.04e-00 & 1.65  & \\
  1/32 &2.84e-02  & 1.00   &5.67e-02 &  1.00 &3.36e-01 & 1.64 &\\
  1/64 &1.42e-02 &1.00 &2.83e-02 & 1.00  &1.13e-01   & 1.57  & \\
\bottomrule
\end{tabular}
\end{table}

\begin{table}[htbp]
\caption{Convergence results  with case II  in 2D.}\label{caseII-2D}
\centering
\begin{tabular}{ccccccccccccccccccccccccccccccc}
\toprule
$h$ & $\|\phi^{N_T}-\phi_{h}^{N_T}\|$ & rate & $\|\nabla (\phi^{N_T}-\phi_{h}^{N_T})\|$ & rate & $\|\u^{N_T}-\u_{h}^{N_T}\|$ & rate & $\|\nabla (\u^{N_T}-\u_{h}^{N_T})\|$ & rate &\\ \hline
 1/8 &2.66e-01 & 1.48  &3.62e-01 &  1.20   &5.99e-03 & 3.11 &4.47e-02 & 1.88  & \\
 1/16 &7.48e-02 & 1.83 &1.57e-01 & 1.21  &9.48e-04 & 2.66   &1.15e-02 & 1.96  & \\
 1/32 &1.93e-02 &1.95  &7.36e-02 & 1.09  &2.03e-04 & 2.22 &2.89e-03 & 1.99  &\\
 1/64 &4.86e-03 & 1.99 &3.61e-02 & 1.03  & 4.87e-05  & 2.01  &7.23e-04 & 2.00  &\\
\bottomrule
\toprule
$h$ &$\|\B^{N_T}-\B_{h}^{N_T}\|$ & rate & $\| \B^{N_T}-\B_{h}^{N_T} \|_{H(\rm{curl})}$ & rate & $\|p^{N_T}-p_{h}^{N_T}\|$ & rate & \\ 
\hline
  1/8 &6.68e-03 & 2.08  &1.71e-02 &  2.00  & 1.73e-00 &1.20  & \\
  1/16 &1.61e-03 & 2.05   &4.28e-03 & 2.00   &5.15e-01 & 1.75  & \\
  1/32 &3.99e-04 & 2.02  &1.07e-03 &  2.00 &1.35e-01 & 1.93 &\\
  1/64 & 9.95e-05   & 2.00 &2.68e-04 & 2.00  &3.41e-02  & 1.98  & \\
\bottomrule
\end{tabular}
\end{table}
 
\begin{table}[h]
\caption{Convergence results  with case I in 3D.}\label{caseI-3D}
\centering
\begin{tabular}{ccccccccccccccccccccccccccccccc}
\toprule
$h$ & $\|\phi^{N_T}-\phi_{h}^{N_T}\|$ & rate & $\|\nabla (\phi^{N_T}-\phi_{h}^{N_T})\|$ & rate & $\|\u^{N_T}-\u_{h}^{N_T}\|$ & rate & $\|\nabla (\u^{N_T}-\u_{h}^{N_T})\|$ & rate &\\ \hline
1/4 &1.09e-00  &    &9.59e-01  &      &1.29e-01 &   &4.15e-01  &    & \\
 1/8 &4.37e-01 & 1.32 &4.54e-01 & 1.08  &3.49e-02 & 1.88   &2.05e-01 & 1.01  & \\
 1/12 &2.17e-01 &1.73  &2.74e-01 & 1.25  &1.64e-02 & 1.86 &1.36e-01 & 1.01  &\\
 1/16 &1.27e-01 & 1.86 &1.93e-01  & 1.21  &9.49e-03 & 1.91  &1.02e-01 & 1.02  &\\
\bottomrule
\toprule
$h$ &$\|\B^{N_T}-\B_{h}^{N_T}\|$ & rate & $\|\B^{N_T}-\B_{h}^{N_T} \|_{H(\rm{curl})}$ & rate & $\|p^{N_T}-p_{h}^{N_T}\|$ & rate & \\ \hline
  1/4 &4.24e-01 &    &8.98e-00 &     &6.13e-01 &  & \\
  1/8 &1.90e-01 & 1.16   &4.55e-00  & 0.98  &1.54e-01 & 1.99  & \\
  1/12 &1.22e-01  & 1.10   &3.04e-00 &  1.00 &7.06e-02 & 1.93 &\\
  1/16 &8.99e-02 & 1.10 &2.28e-00 &  1.00  &4.02e-02 & 1.95  & \\
\bottomrule
\end{tabular}
\end{table}

\subsection{Spinodal decomposition}
The spinodal decomposition is a phase separation phenomenon that occurs in binary or multi-component alloys, polymer blends and liquid crystals \cite{2023Energy, shi2024structure}.
The computational domain is $\Omega=[0,1]^{2}$.  The initial values read as  
\begin{equation}\label{fai0}
\left\{
\begin{aligned} 
\phi_{0}&=-0.05+0.001{ \rm rand}(x),\\
\u_{0}&= \B_0 =\0, \quad p_{0}=0,
\end{aligned}
\right.
\end{equation}
where $\mathrm{rand}(x)$ is a uniformly distributed random function in $[-1,1]$ with zero mean. We select finite element pairs   cases I  to test the spinodal decomposition phenomenon. The parameters are given as
\begin{equation*}
\gamma=1/100,\quad  M=1,\quad \nu=1,\quad \mu=1,\quad \lambda=1,\quad \sigma=1.
\end{equation*}

\begin{figure}[h]
	\centering
\subfigure[t=0.0001]{
		\begin{minipage}[t]{0.2\linewidth}
			\centering
			\includegraphics[width=\textwidth]{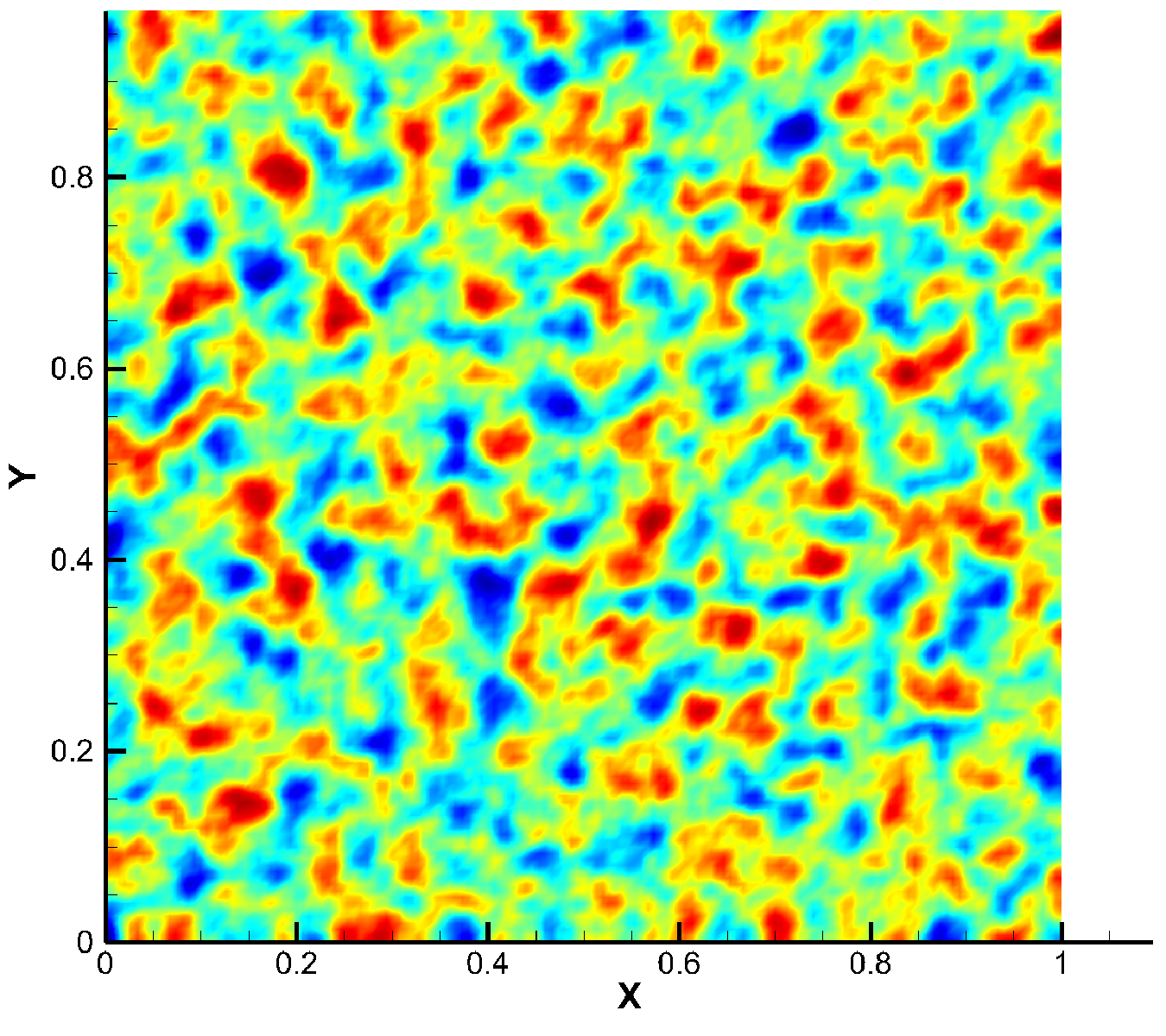}
		\end{minipage}
	}%
	\subfigure[t=0.05]{
		\begin{minipage}[t]{0.2\linewidth}
			\centering
			\includegraphics[width=\textwidth]{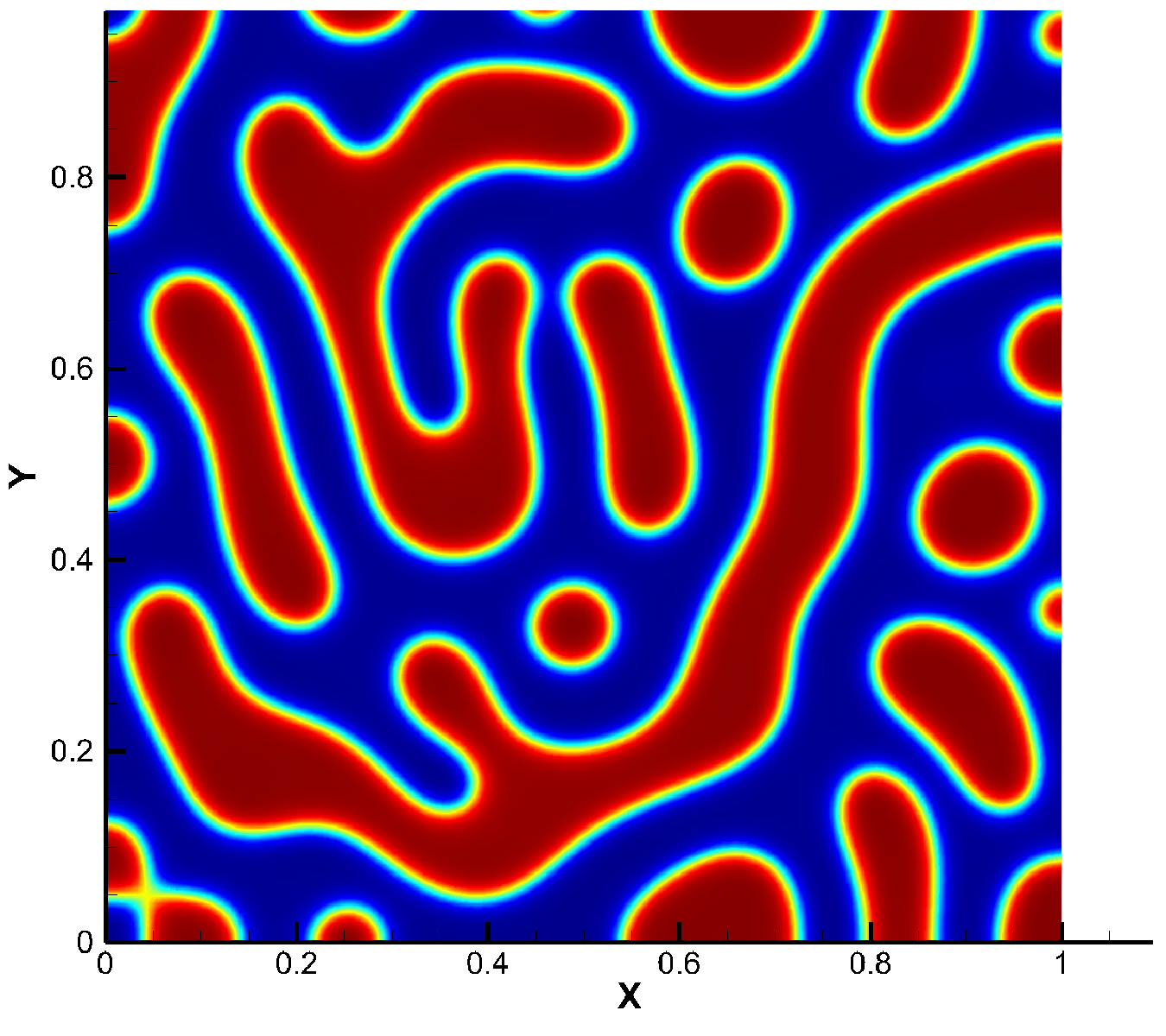}
		\end{minipage}
	}%
	\subfigure[t=0.5]{
		\begin{minipage}[t]{0.2\linewidth}
			\centering
			\includegraphics[width=\textwidth]{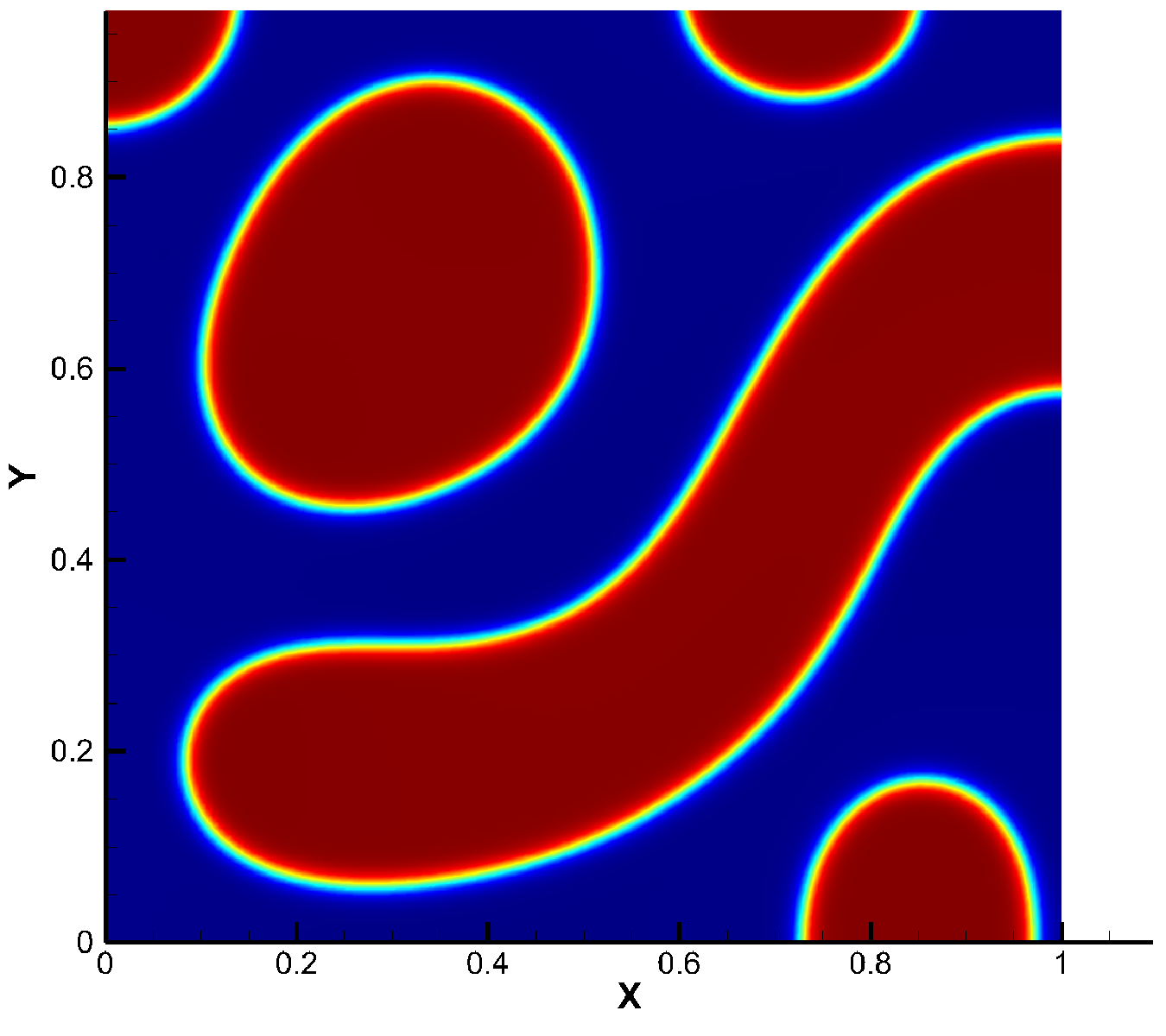}
		\end{minipage}
	}%
	\subfigure[t=2.5]{
		\begin{minipage}[t]{0.2\linewidth}
			\centering
			\includegraphics[width=\textwidth]{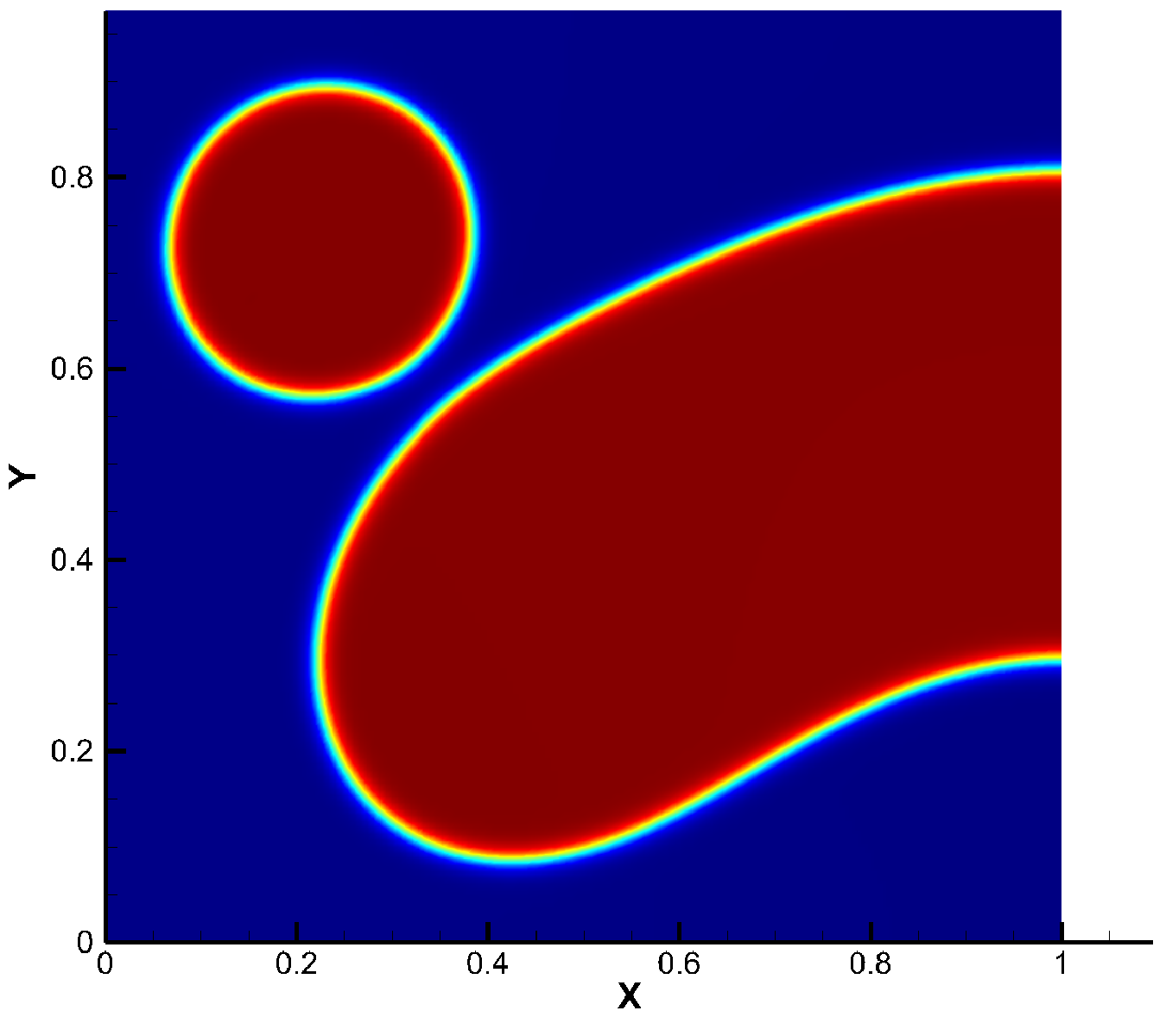}
		\end{minipage}
	}%
   \subfigure[t=4]{
		\begin{minipage}[t]{0.2\linewidth}
			\centering
			\includegraphics[width=\textwidth]{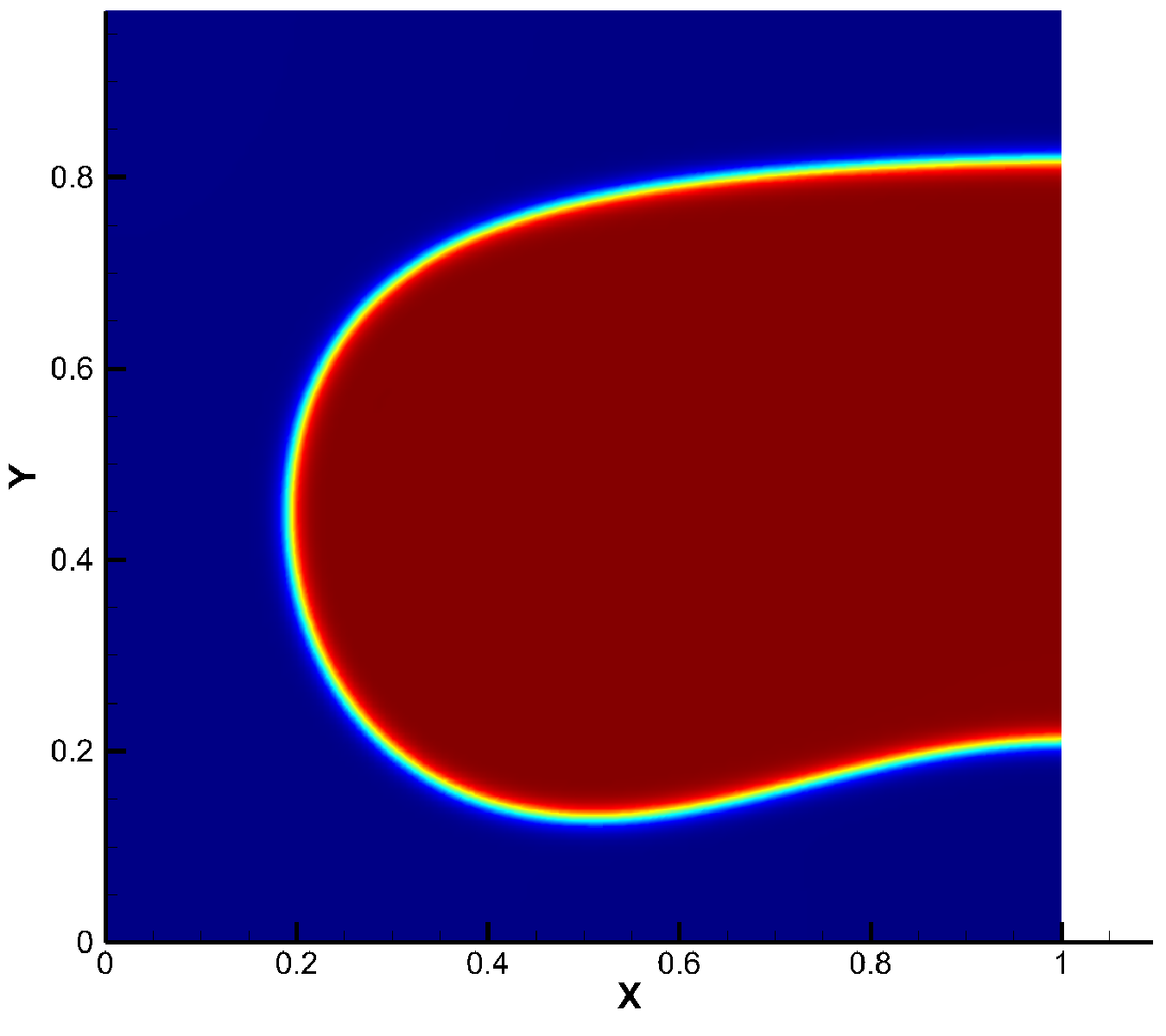}
		\end{minipage}
	}%
	\centering
	\caption{Snapshots of phase field dynamical evolution for spinodal decomposition for case I.}
	\label{spinodal-1}
\end{figure}

We apply the homogeneous Dirichlet boundary conditions to the velocity and magnetic induction fields, and enforce the homogeneous Neumann boundary conditions for the phase field and chemical potential. The time step size $\Delta t$=1/1000 and the mesh size $h$=1/150 are selected to investigate the evolution of the phase field based on the case I. In Figure \ref{spinodal-1}, we find that over time, the phase field gradually coarsens.

\begin{figure}[h]
	\centering
\subfigure[system energy]{
		\begin{minipage}[t]{0.34\linewidth}
			\centering
			\includegraphics[width=\textwidth]{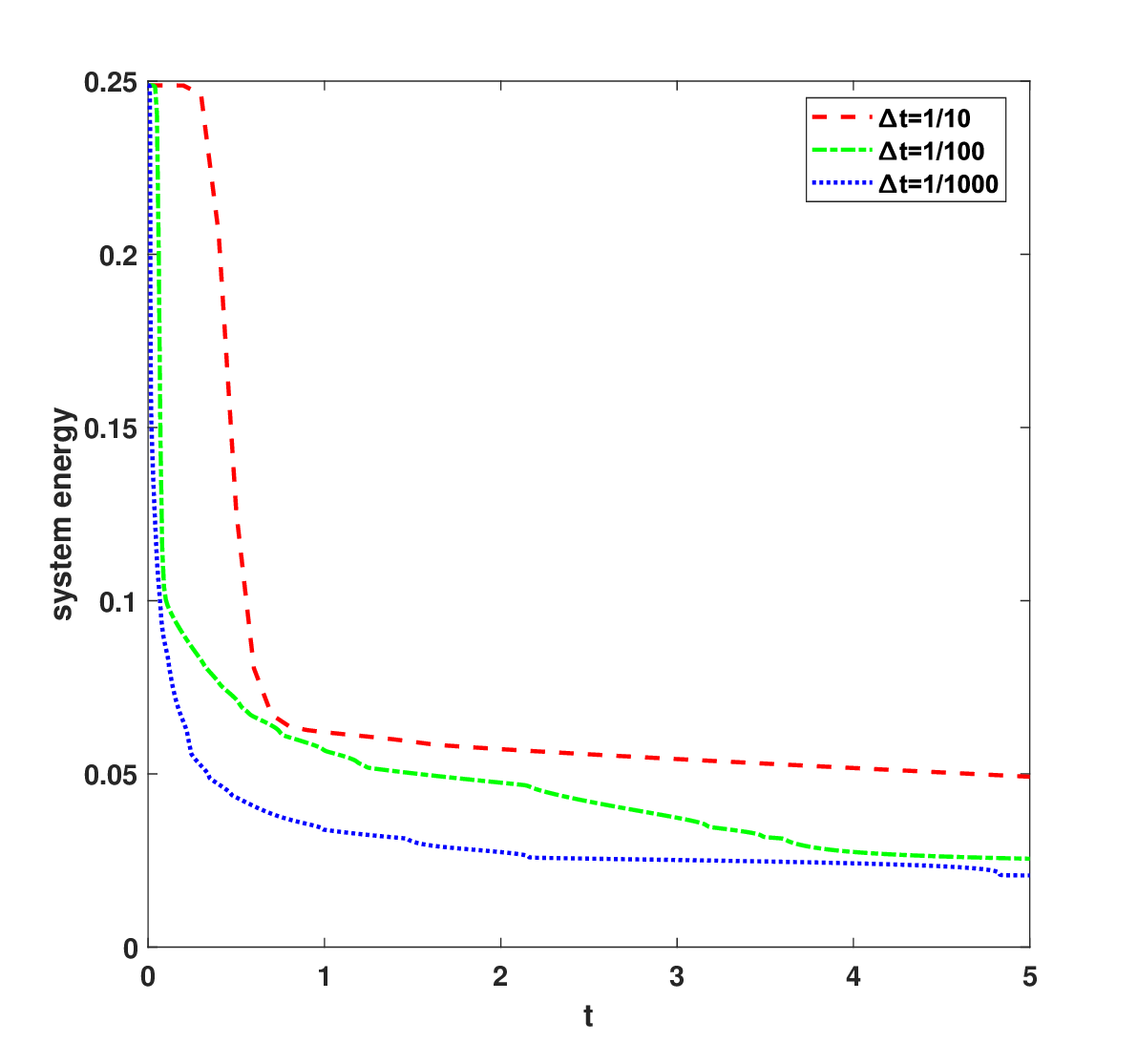}
		\end{minipage}
	}%
	\subfigure[algorithm energy]{
		\begin{minipage}[t]{0.34\linewidth}
			\centering
			\includegraphics[width=\textwidth]{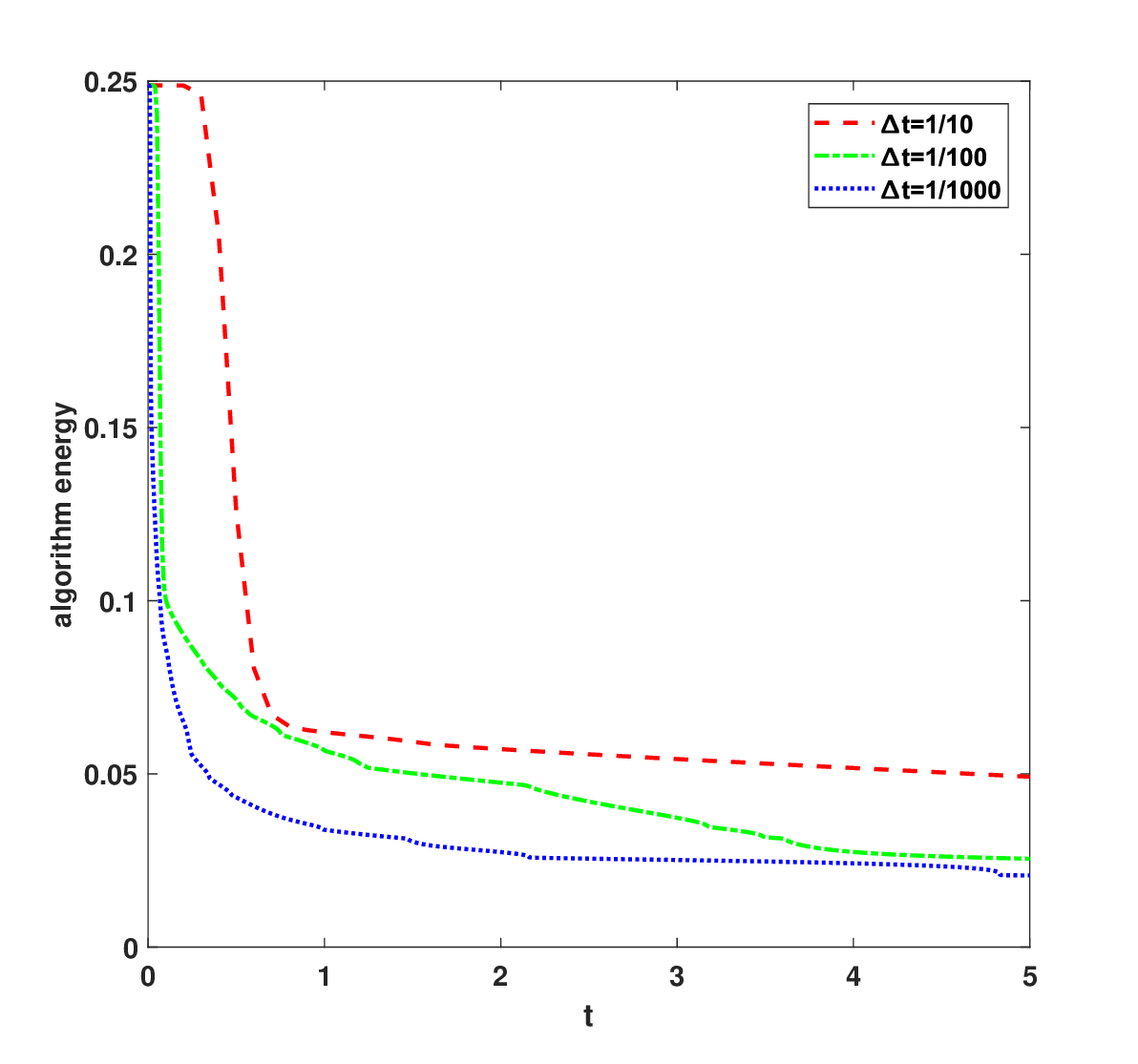}
		\end{minipage}
	}%
	\subfigure[discrete mass]{
		\begin{minipage}[t]{0.34\linewidth}
			\centering
			\includegraphics[width=\textwidth]{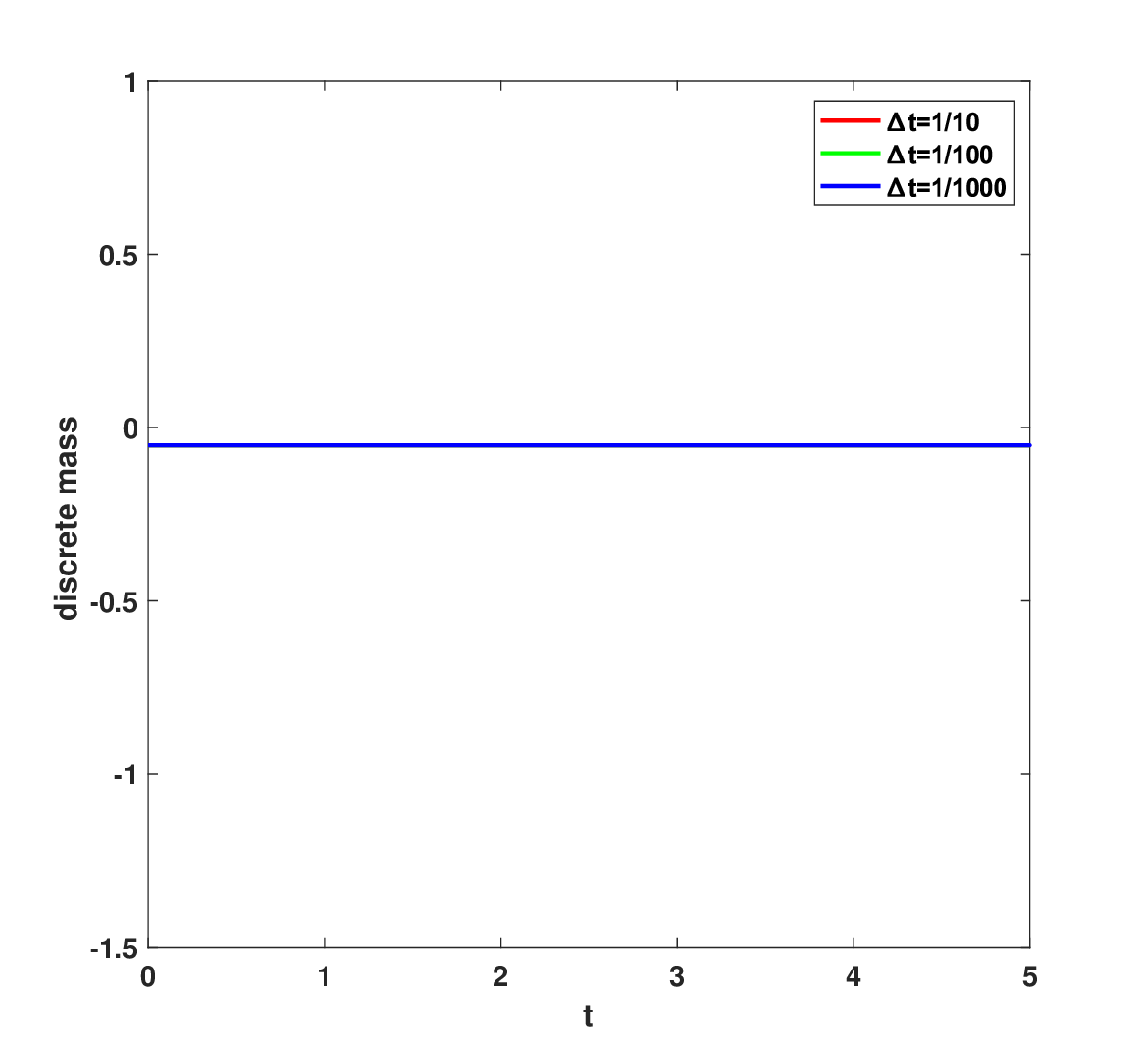}
		\end{minipage}
	}%
	\centering
	\caption{The system energy (left), algorithm energy (middle) and the discrete mass (right) for case I.}
	\label{energy-mass}
\end{figure}

Then we conduct the system energy test (\ref{system-energy}), the algorithm energy test (\ref{algorithm-energy}), and the discrete mass conservation  test (\ref{mass-conservetion}).  We fix the mesh size $h$=1/120, and set the time step size $\Delta t= 1/10, 1/100$, and $1/1000$ respectively.
The initial values are set according to equations (\ref{fai0}). The parameters are chosen as
\begin{equation*}
\gamma=1/100,\quad M=1,\quad \nu=1,\quad \mu=1,\quad \lambda=1/100,\quad \sigma=1.
\end{equation*}
In Figure \ref{energy-mass}(a),  (b), and  (c), the comparisons of system energy, algorithm energy, and discrete mass at different time steps are plotted for case I. As the time step is refined, the energy curves gradually become flat, and the discrete masses are always conserved. This indicates good numerical consistency, i.e., a smaller time step leads to more stable results in Figure \ref{energy-mass}(a) and (b). Without specific needs, we choose case I in the following contents. 

\subsection{Lid-driven cavity flow}

\begin{figure}[h]
	\centering
\subfigure[t=0.001]{
		\begin{minipage}[t]{0.2\linewidth}
			\centering
			\includegraphics[width=\textwidth]{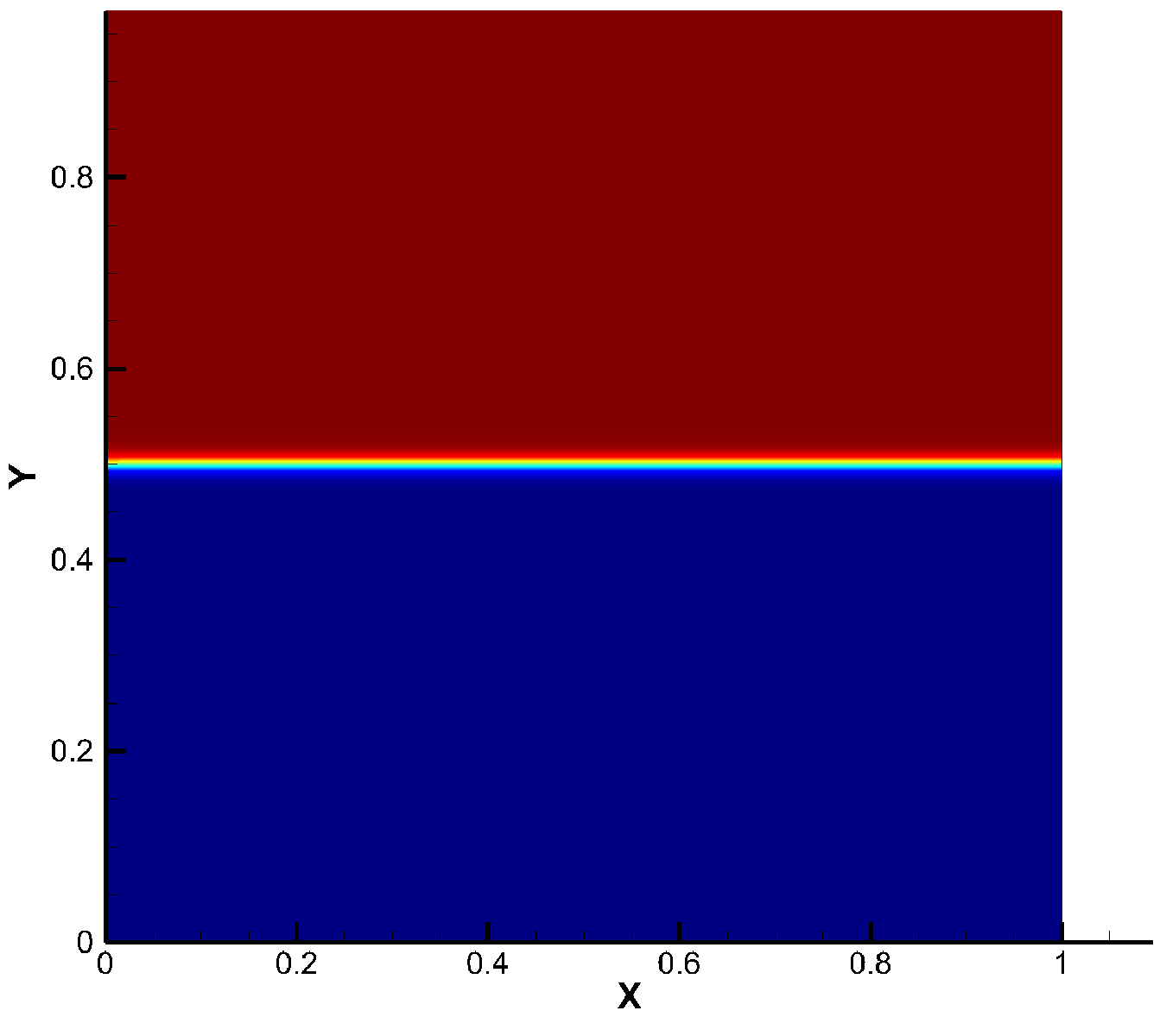}
		\end{minipage}
	}%
	\subfigure[t=2]{
		\begin{minipage}[t]{0.2\linewidth}
			\centering
			\includegraphics[width=\textwidth]{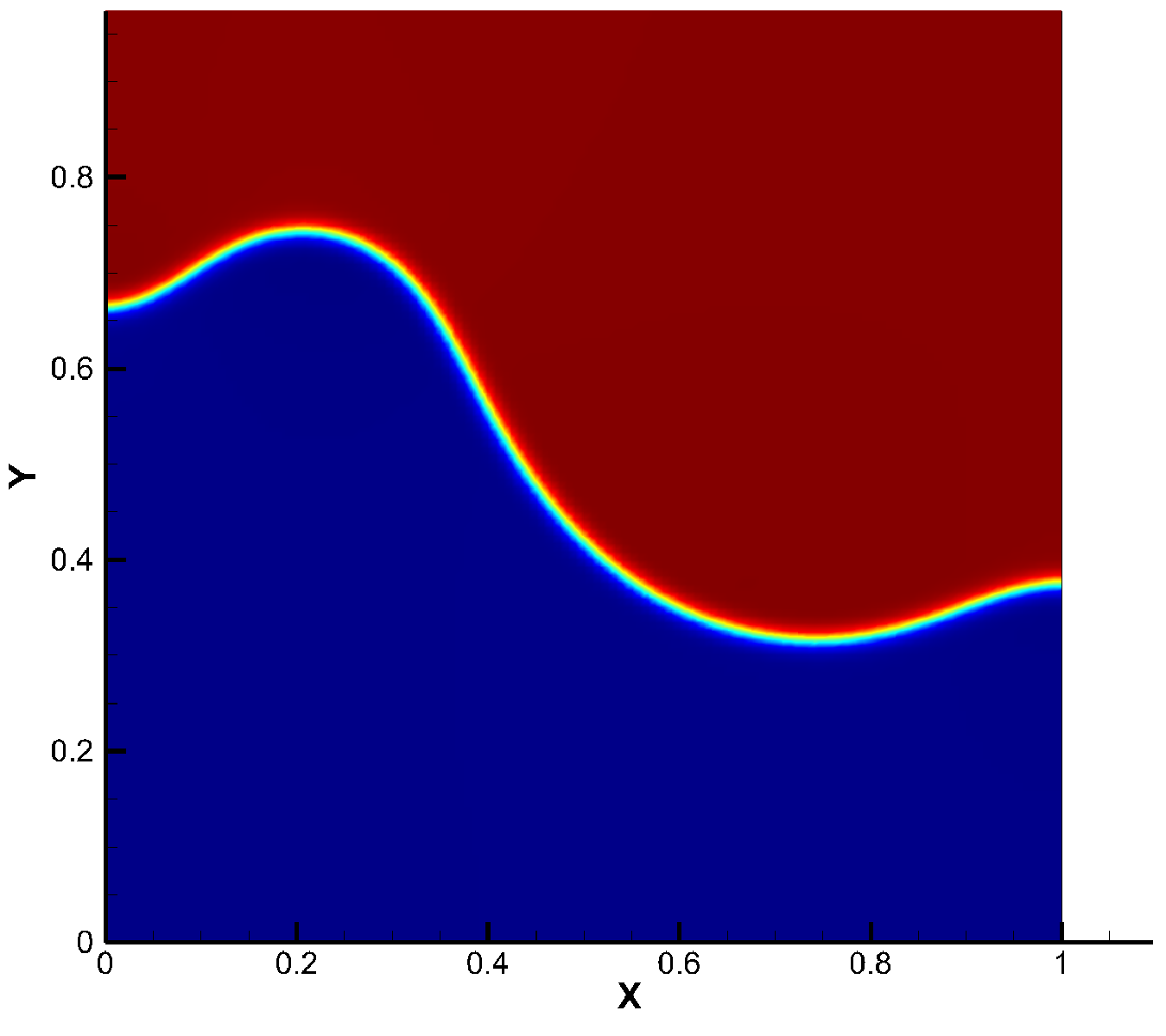}
		\end{minipage}
	}%
	\subfigure[t=3.1]{
		\begin{minipage}[t]{0.2\linewidth}
			\centering
			\includegraphics[width=\textwidth]{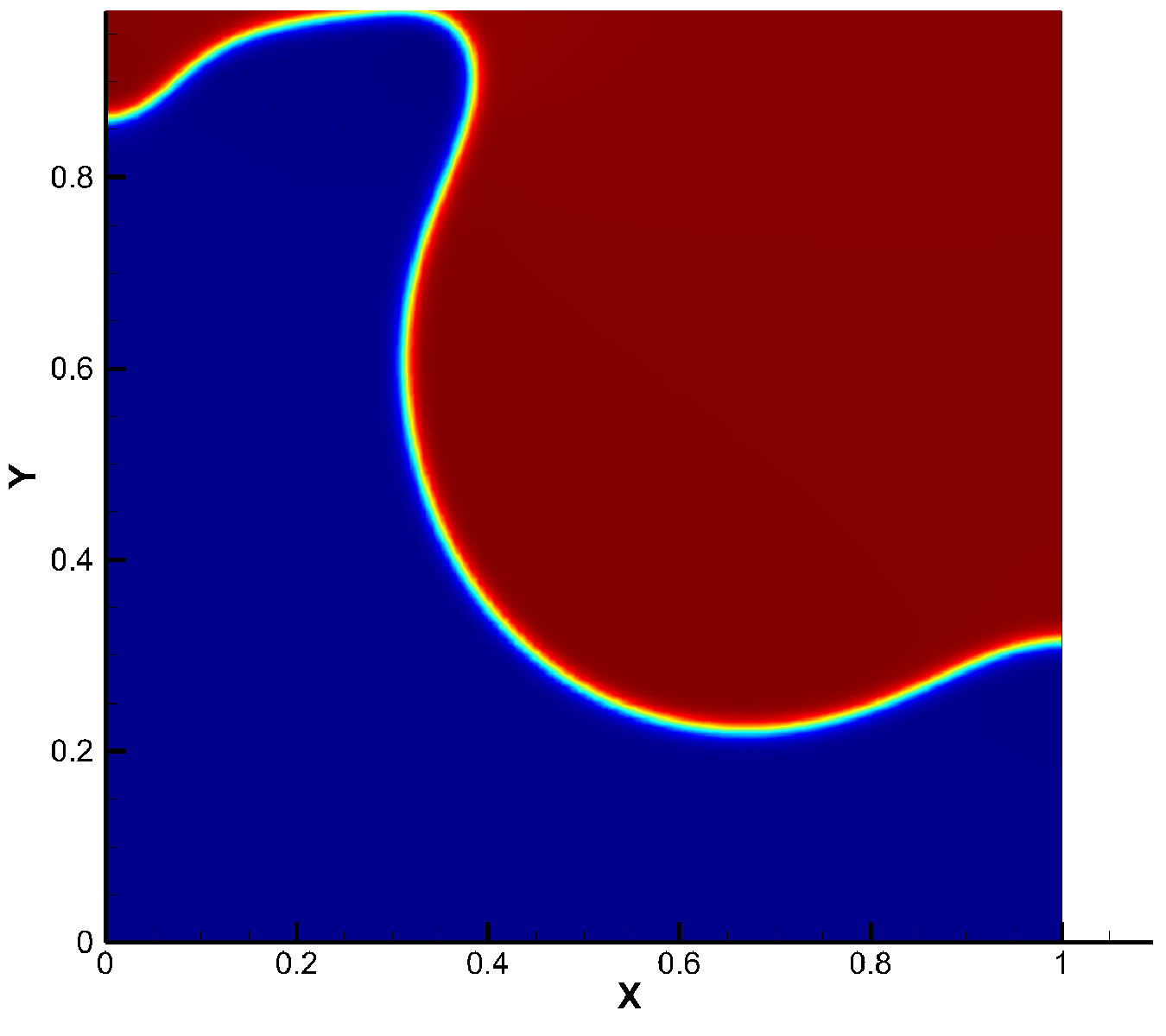}
		\end{minipage}
	}%
	\subfigure[t=5]{
		\begin{minipage}[t]{0.2\linewidth}
			\centering
			\includegraphics[width=\textwidth]{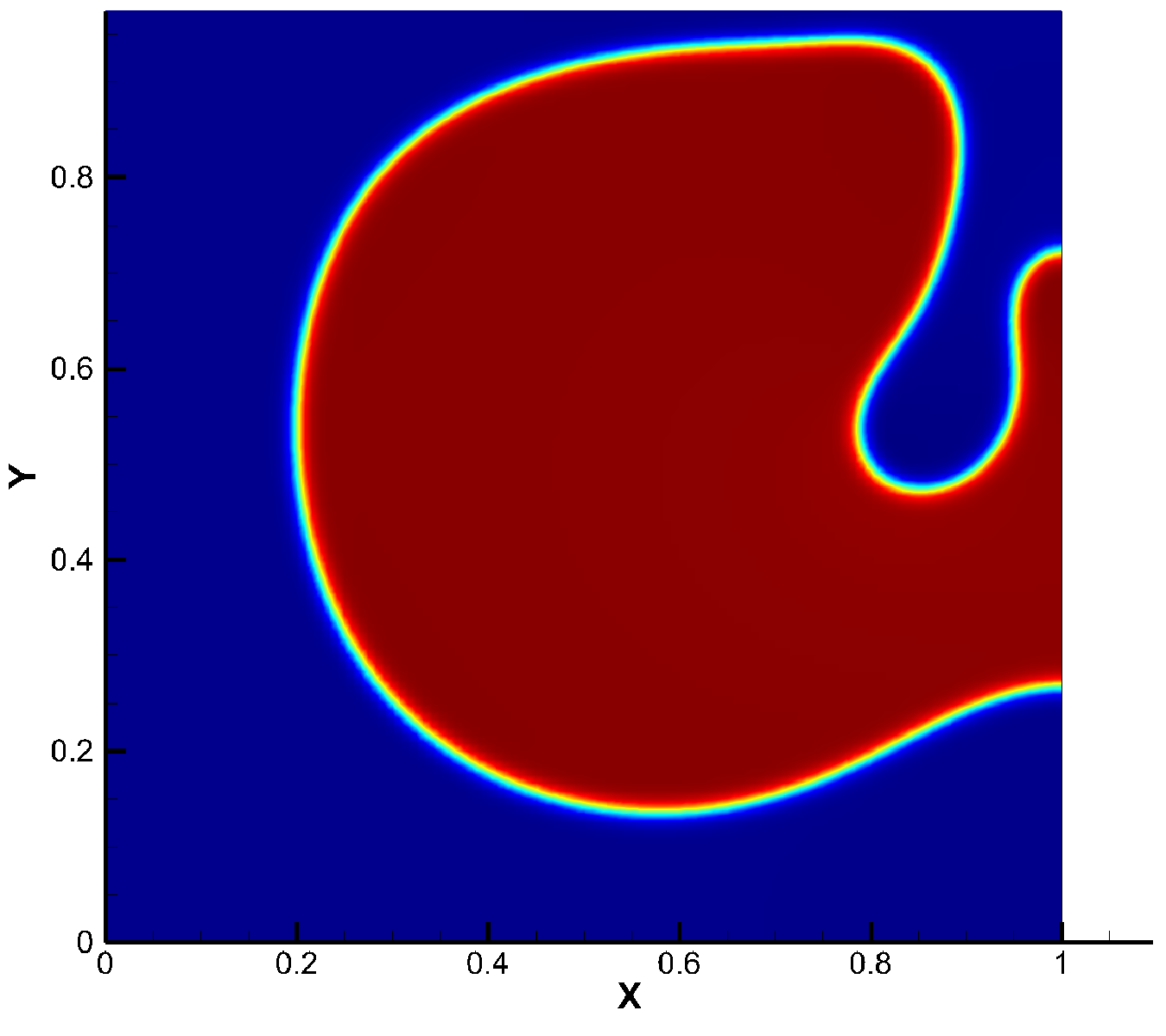}
		\end{minipage}
	}%
   \subfigure[t=9.1]{
		\begin{minipage}[t]{0.2\linewidth}
			\centering
			\includegraphics[width=\textwidth]{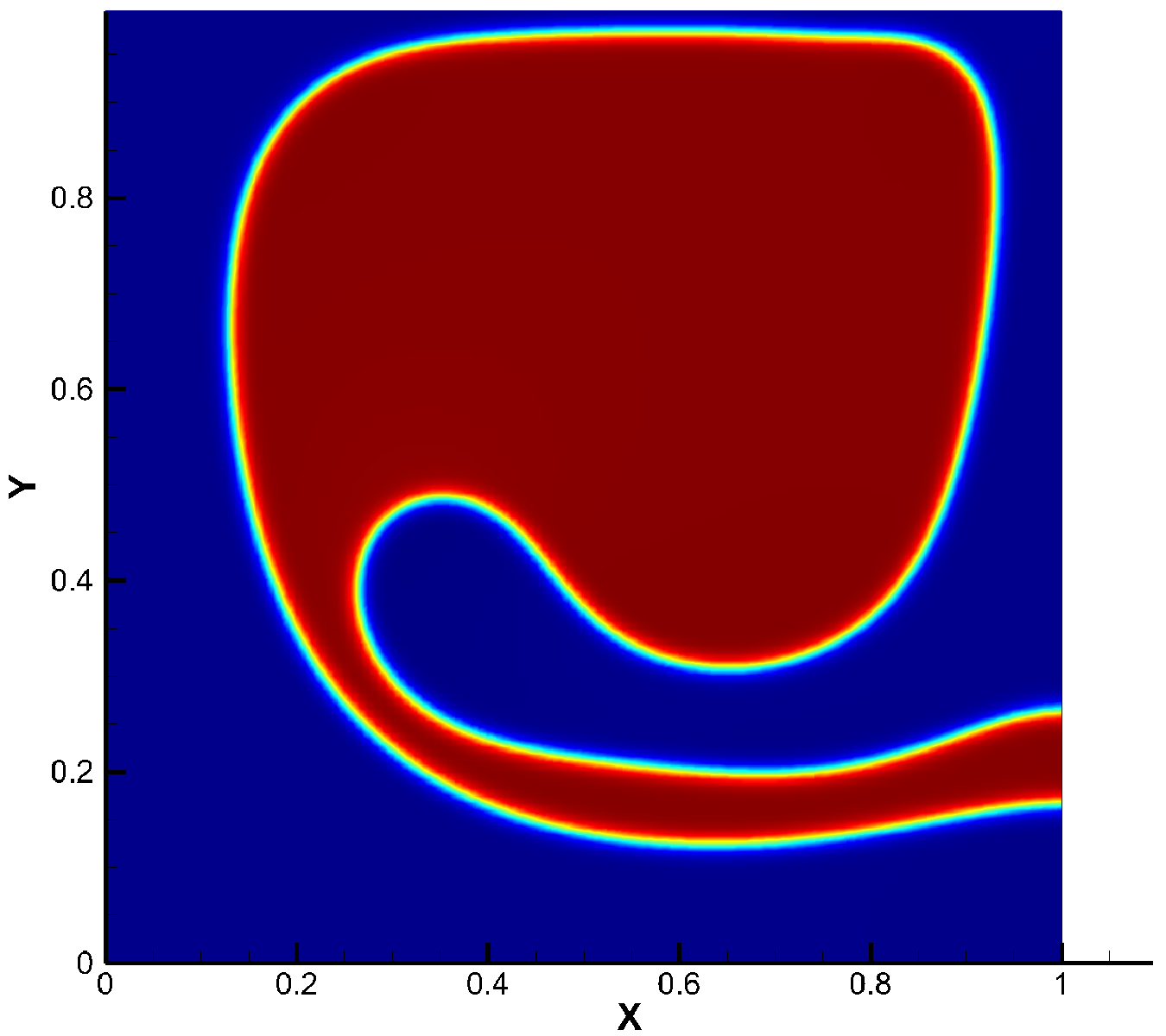}
		\end{minipage}
	}\vspace{5pt}
\subfigure[t=0.001]{
		\begin{minipage}[t]{0.2\linewidth}
			\centering
			\includegraphics[width=\textwidth]{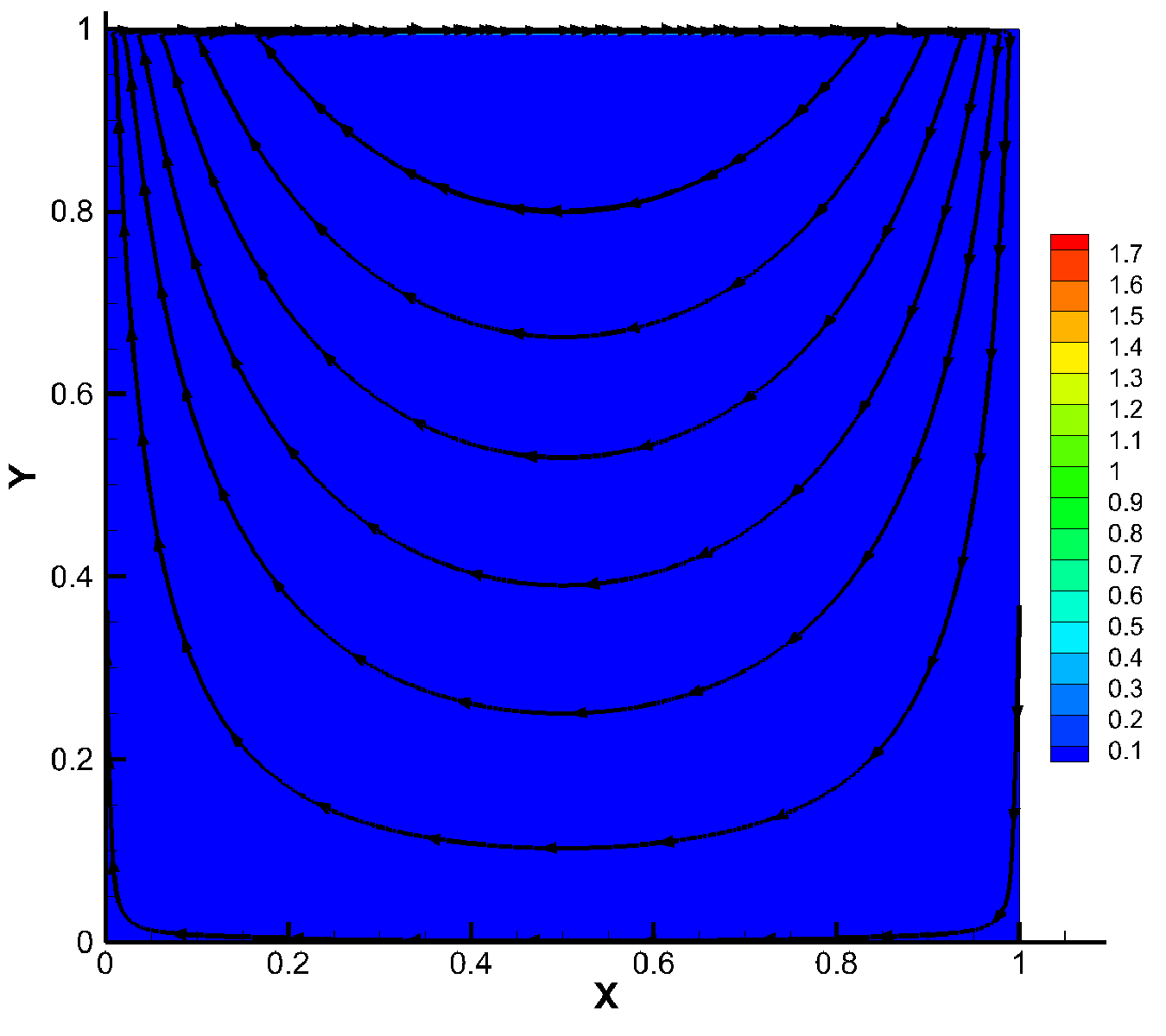}
		\end{minipage}
	}%
	\subfigure[t=2]{
		\begin{minipage}[t]{0.2\linewidth}
			\centering
			\includegraphics[width=\textwidth]{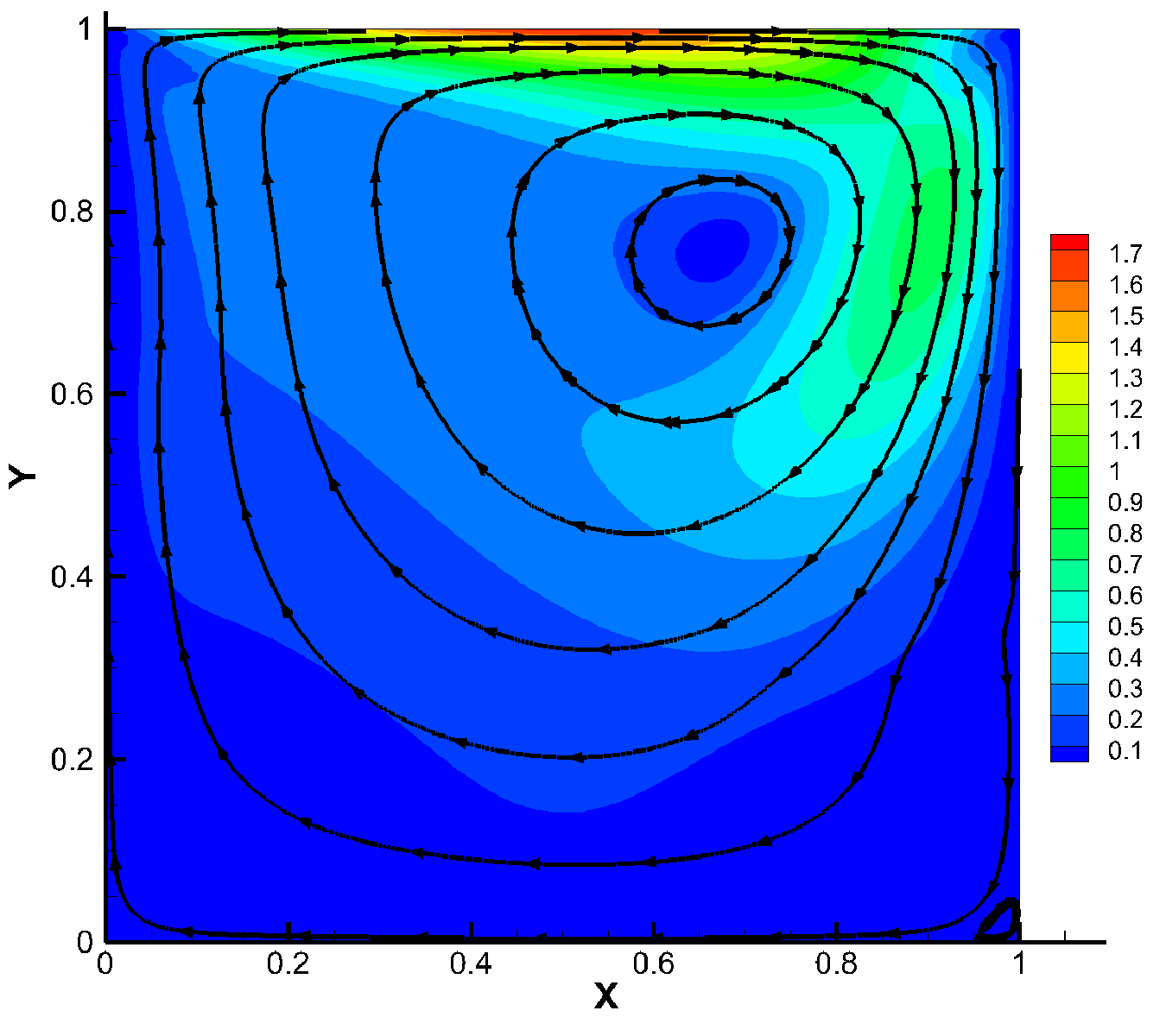}
		\end{minipage}
	}%
	\subfigure[t=3.1]{
		\begin{minipage}[t]{0.2\linewidth}
			\centering
			\includegraphics[width=\textwidth]{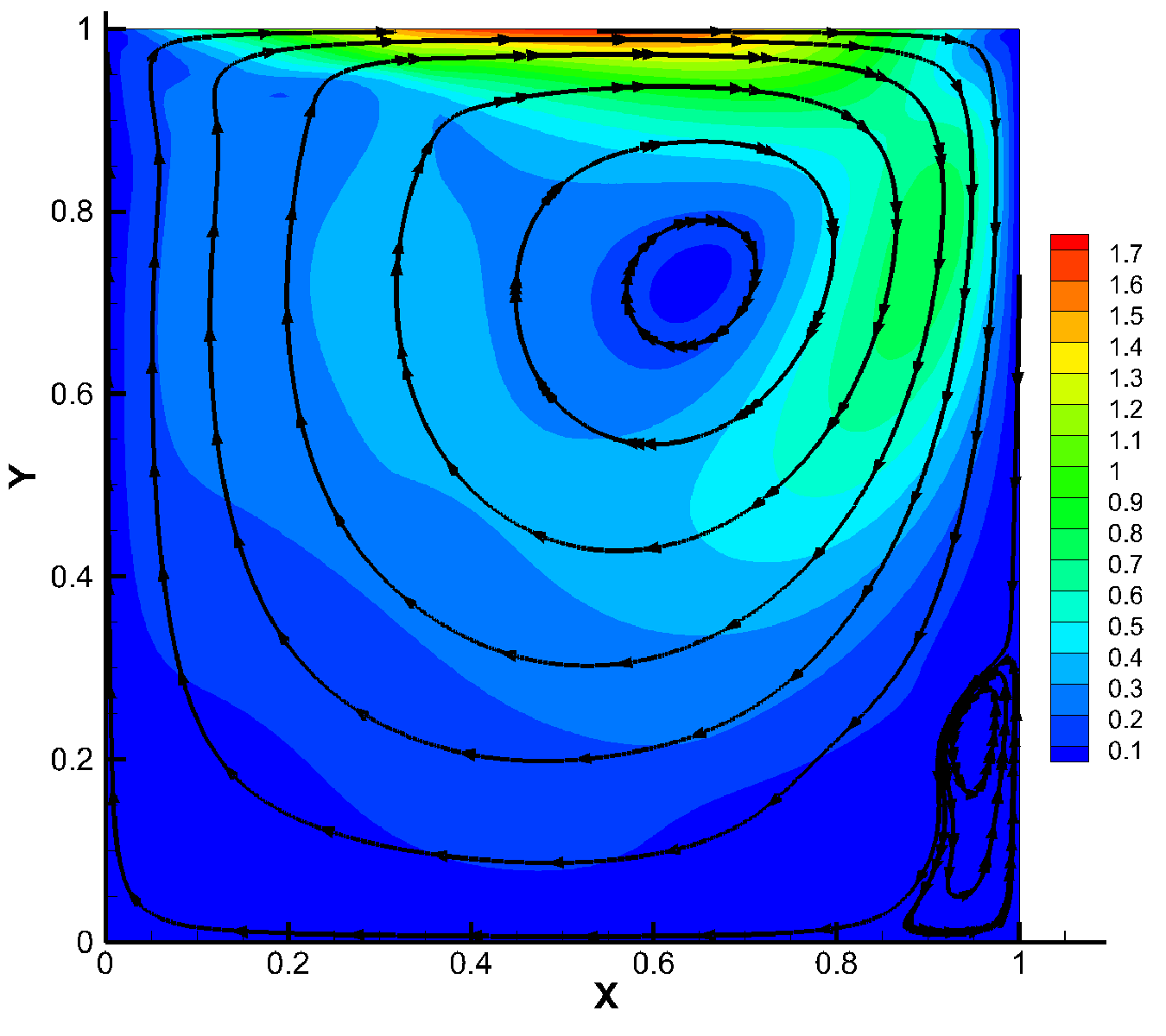}
		\end{minipage}
	}%
	\subfigure[t=5]{
		\begin{minipage}[t]{0.2\linewidth}
			\centering
			\includegraphics[width=\textwidth]{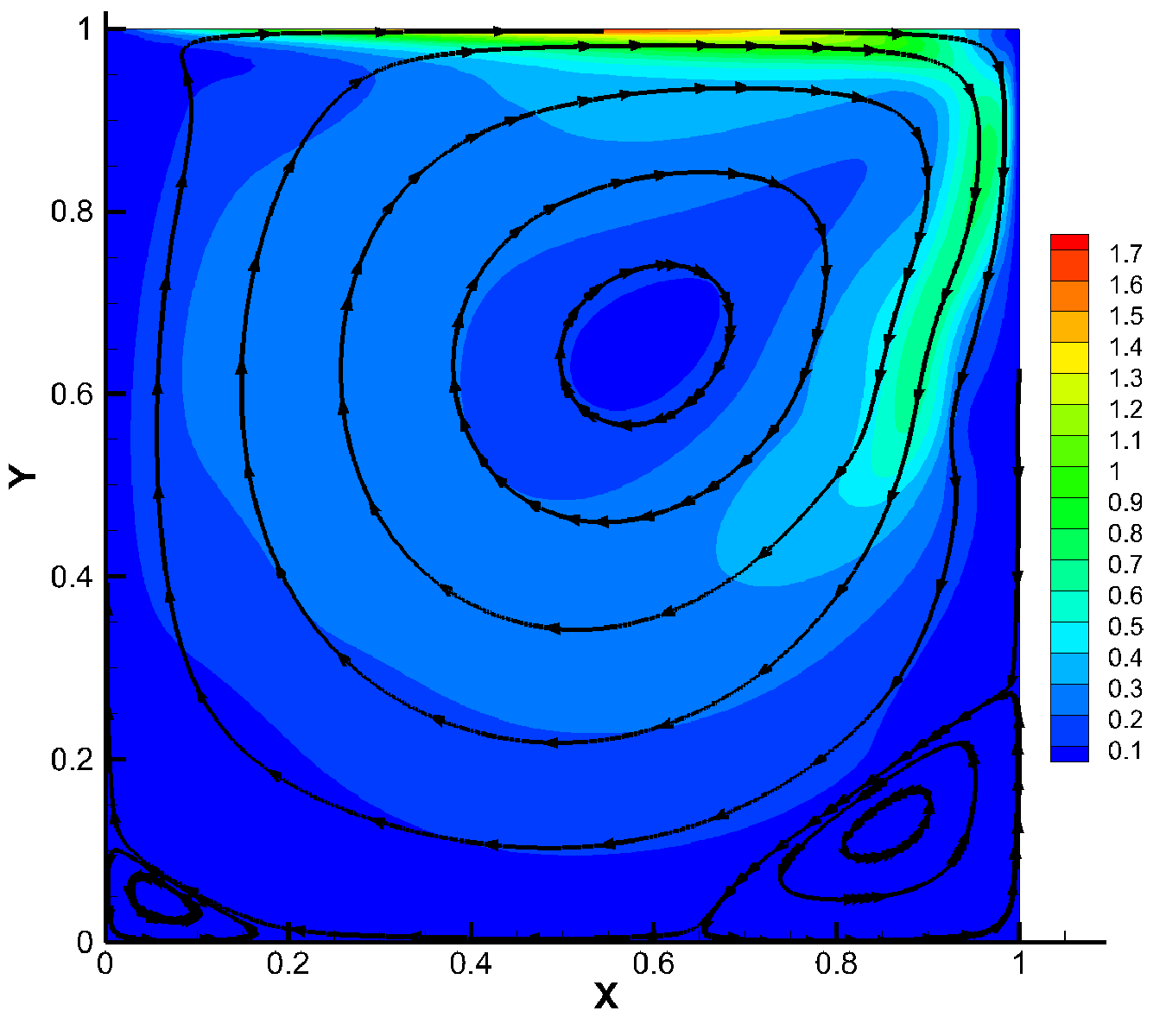}
		\end{minipage}
	}%
   \subfigure[t=9.1]{
		\begin{minipage}[t]{0.2\linewidth}
			\centering
			\includegraphics[width=\textwidth]{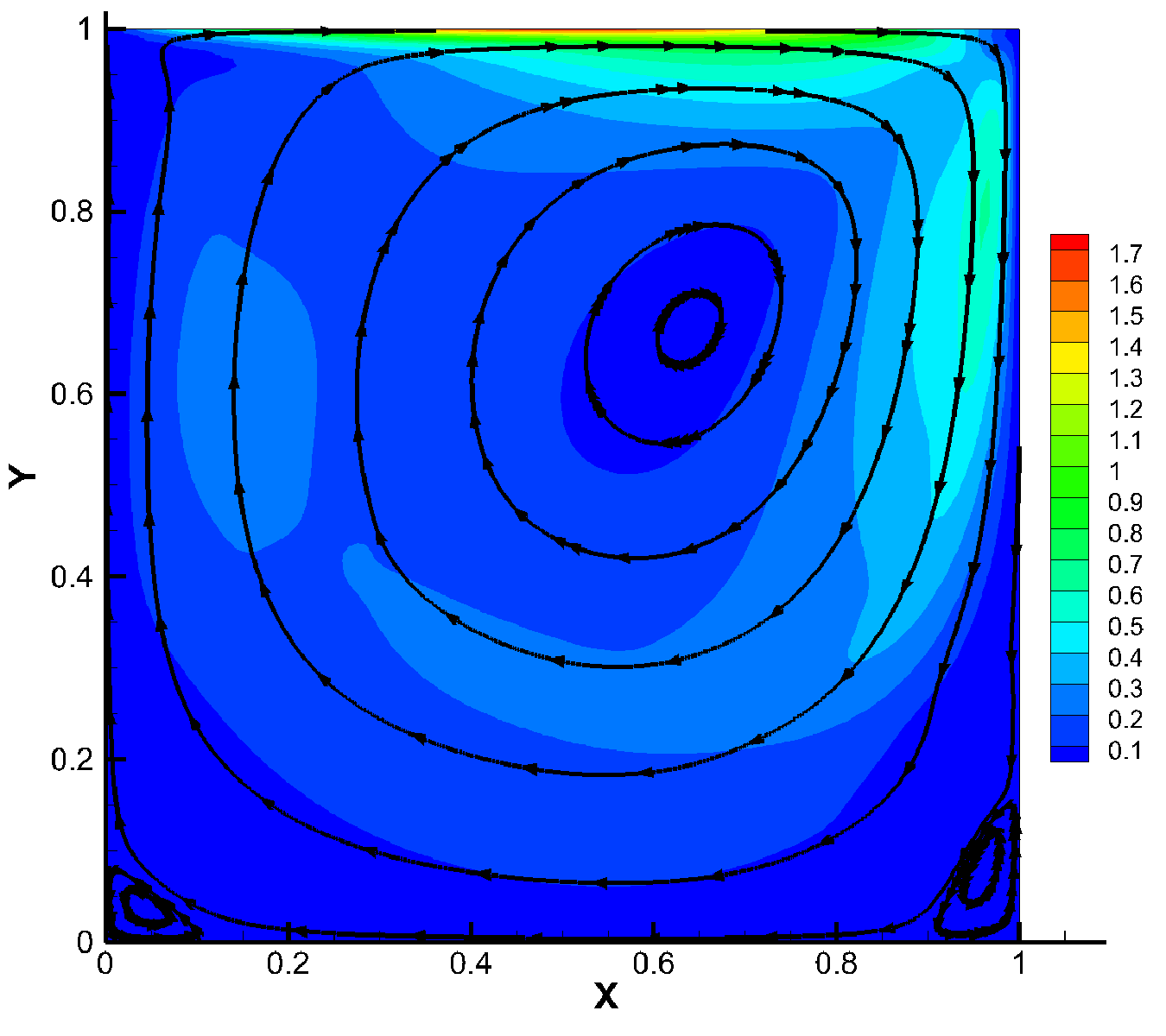}
		\end{minipage}
	}%
	\centering
	\caption{Snapshots of phase field (upper), velocity field (lower) dynamical evolution for lid driven cavity flow with $\mu$=2.}
	\label{lid-mu-2}
\end{figure}

\begin{figure}[h]
	\centering
\subfigure[t=0.001]{
		\begin{minipage}[t]{0.2\linewidth}
			\centering
			\includegraphics[width=\textwidth]{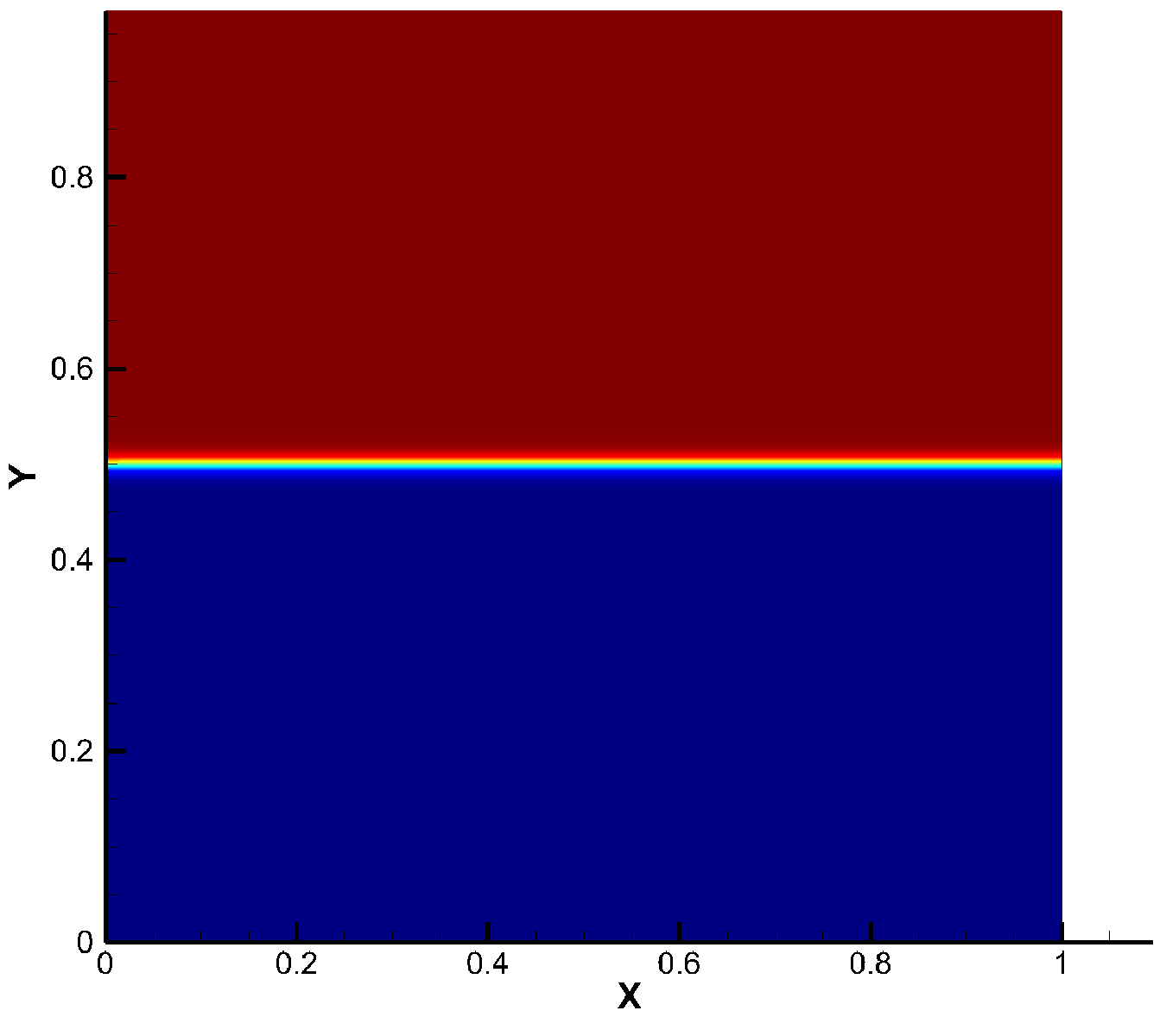}
		\end{minipage}
	}%
	\subfigure[t=2]{
		\begin{minipage}[t]{0.2\linewidth}
			\centering
			\includegraphics[width=\textwidth]{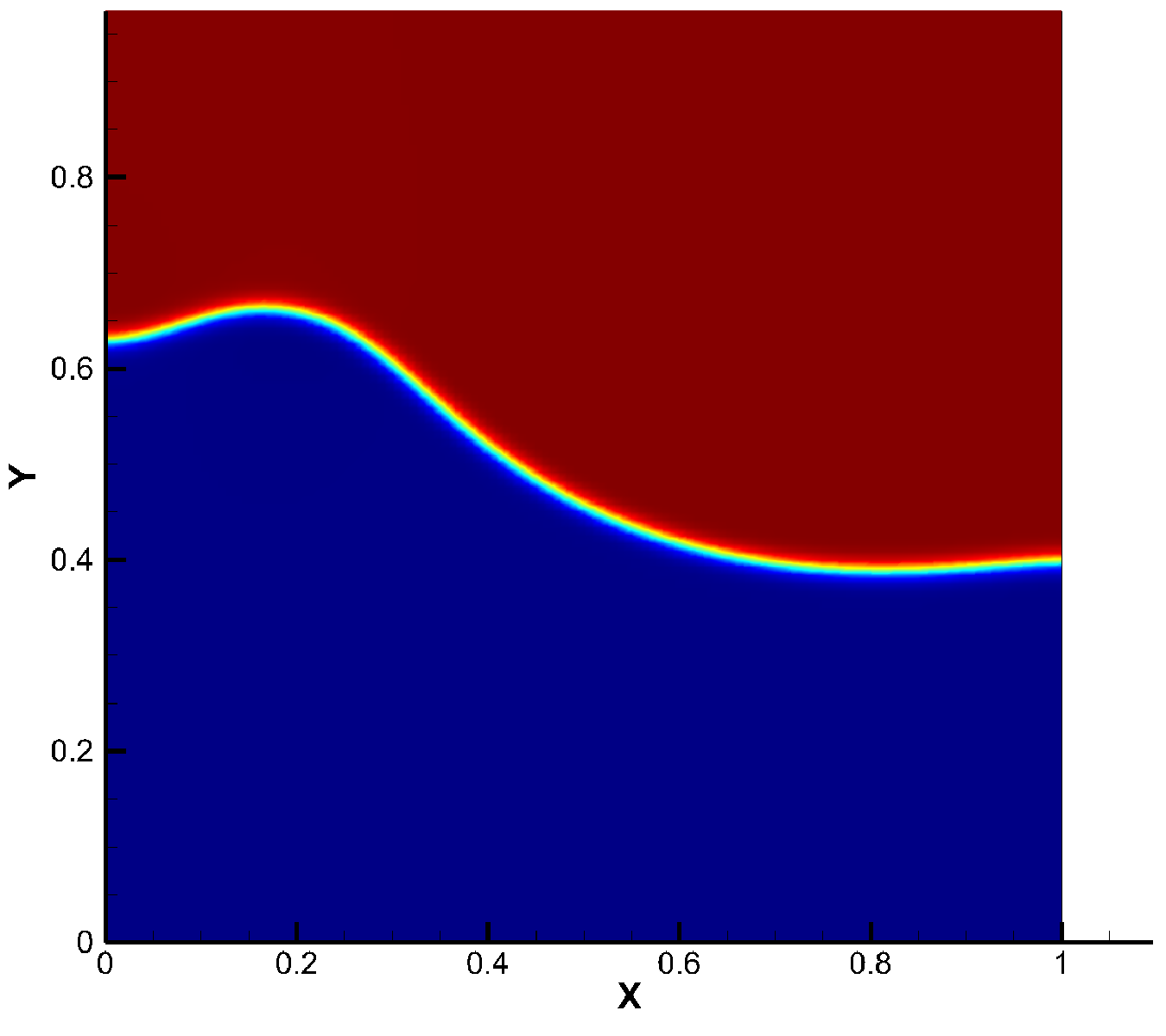}
		\end{minipage}
	}%
	\subfigure[t=3.1]{
		\begin{minipage}[t]{0.2\linewidth}
			\centering
			\includegraphics[width=\textwidth]{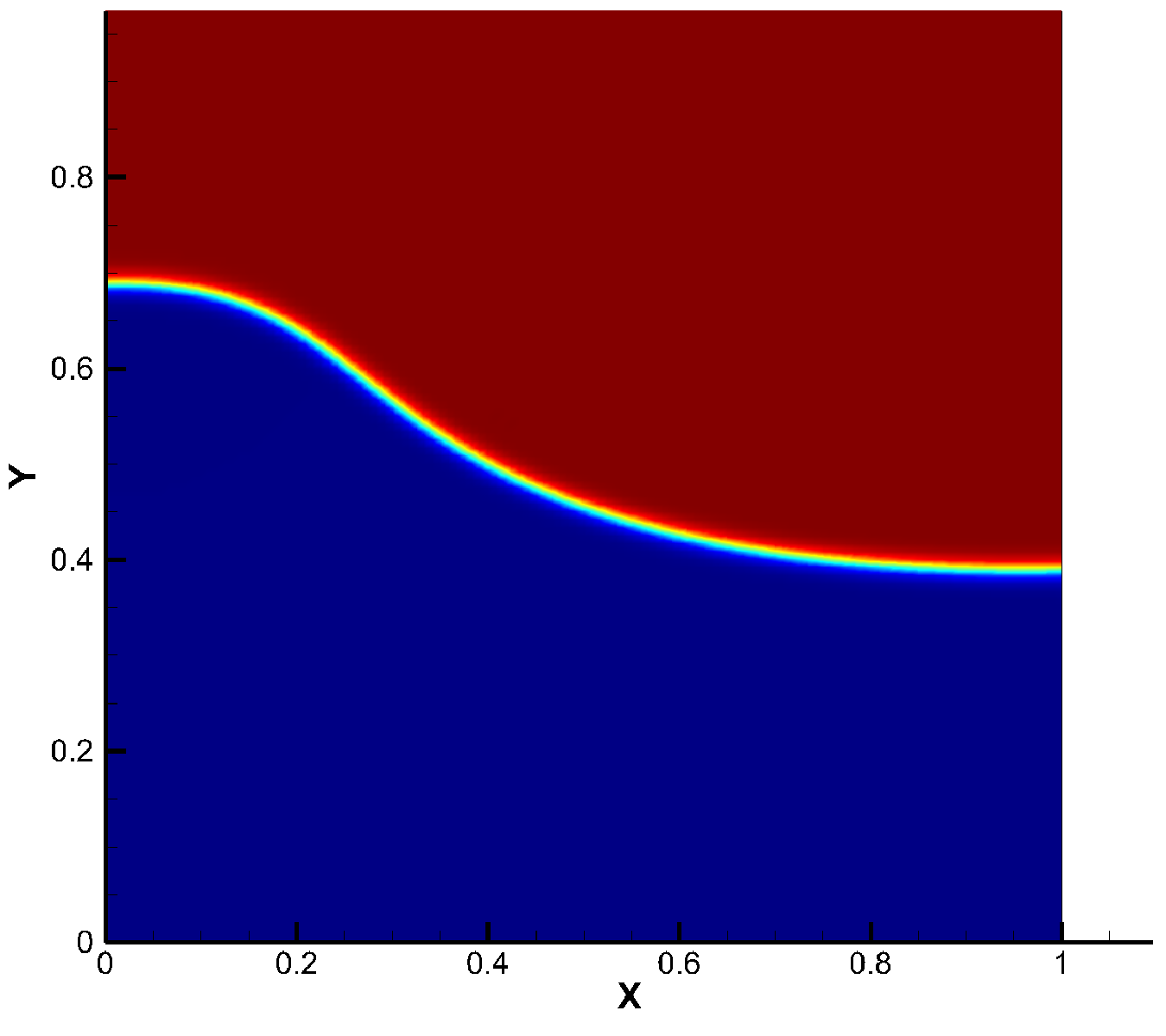}
		\end{minipage}
	}%
	\subfigure[t=5]{
		\begin{minipage}[t]{0.2\linewidth}
			\centering
			\includegraphics[width=\textwidth]{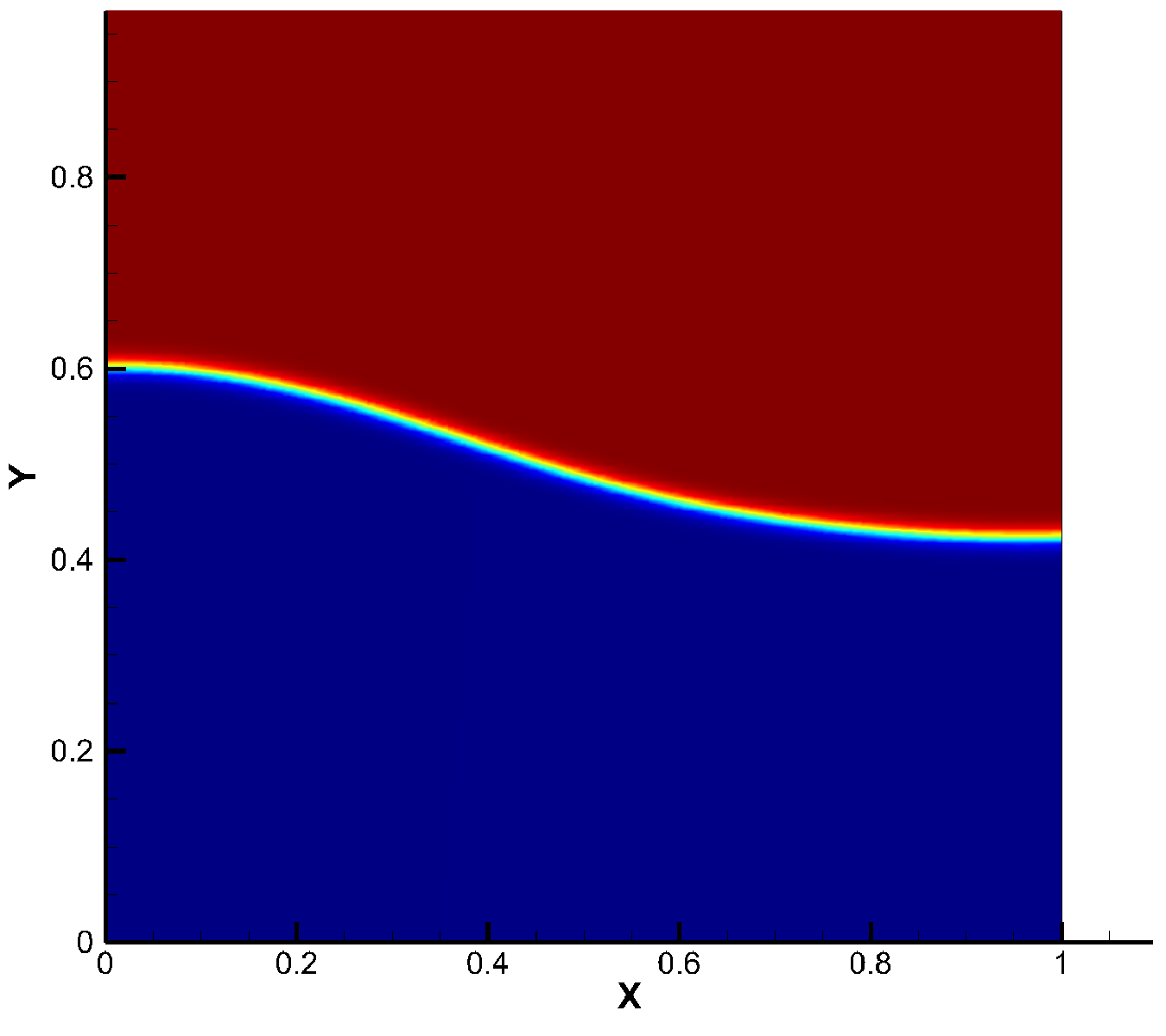}
		\end{minipage}
	}%
   \subfigure[t=9.1]{
		\begin{minipage}[t]{0.2\linewidth}
			\centering
			\includegraphics[width=\textwidth]{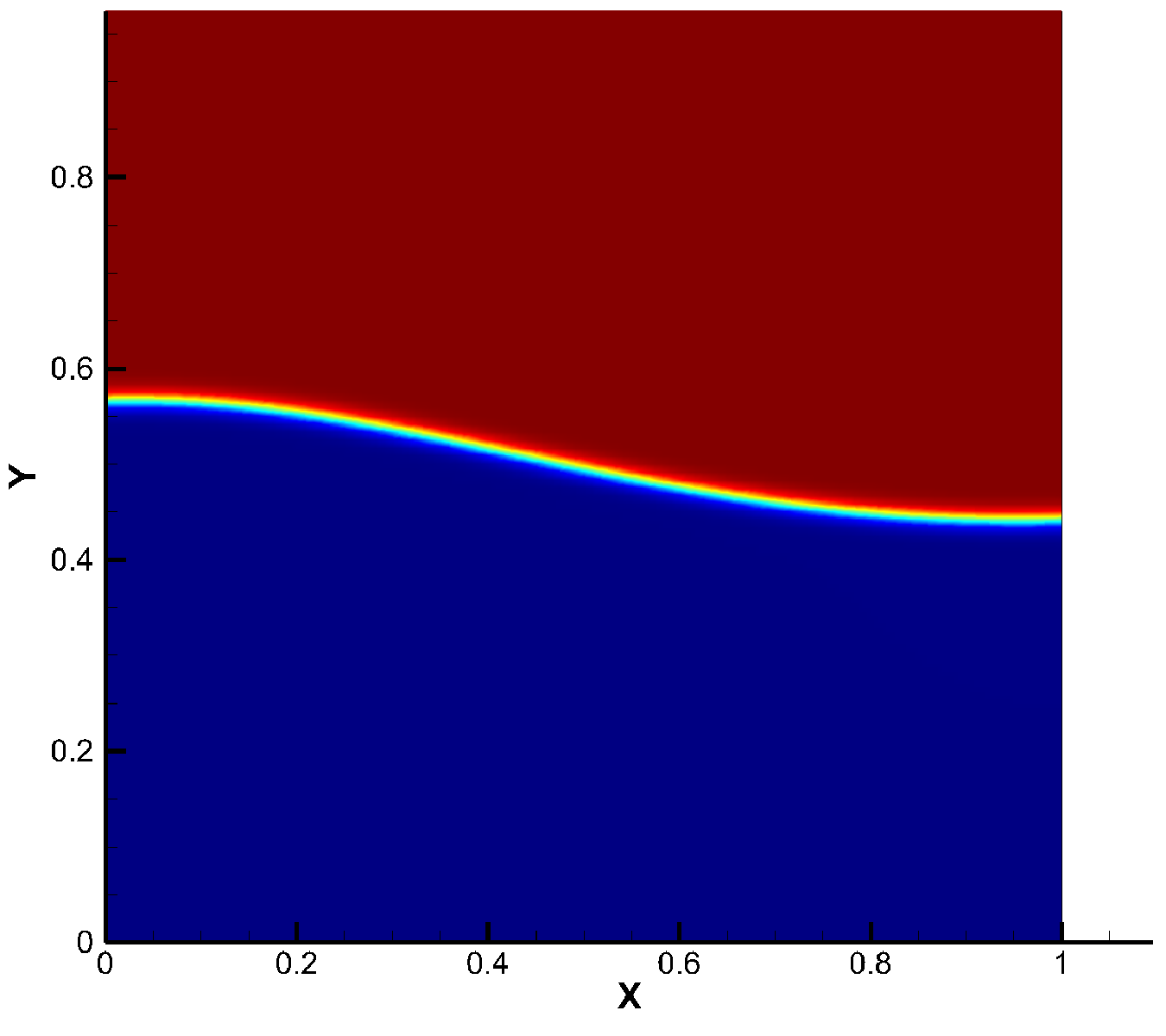}
		\end{minipage}
	}\vspace{5pt}
\subfigure[t=0.001]{
		\begin{minipage}[t]{0.2\linewidth}
			\centering
			\includegraphics[width=\textwidth]{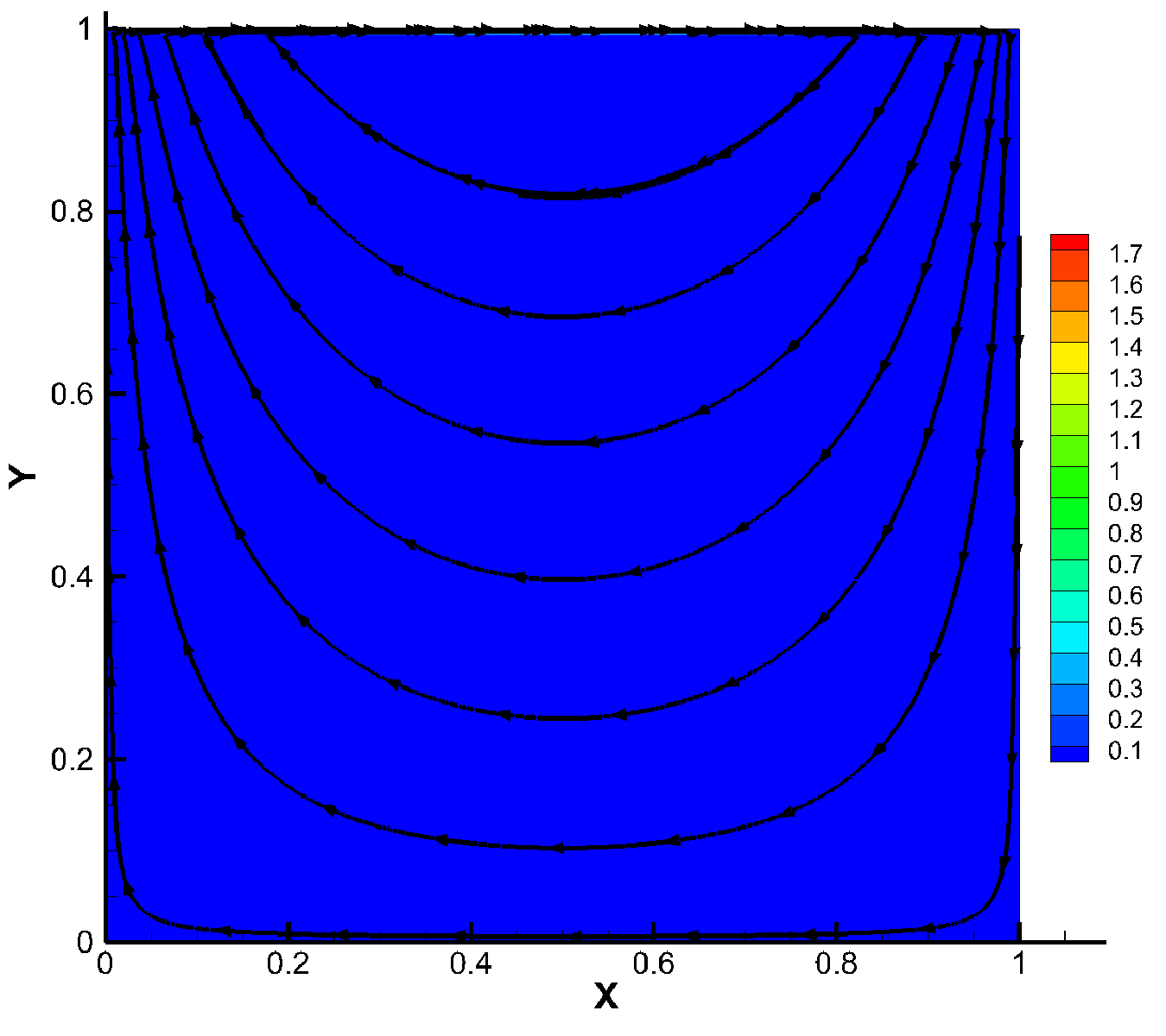}
		\end{minipage}
	}%
	\subfigure[t=2]{
		\begin{minipage}[t]{0.2\linewidth}
			\centering
			\includegraphics[width=\textwidth]{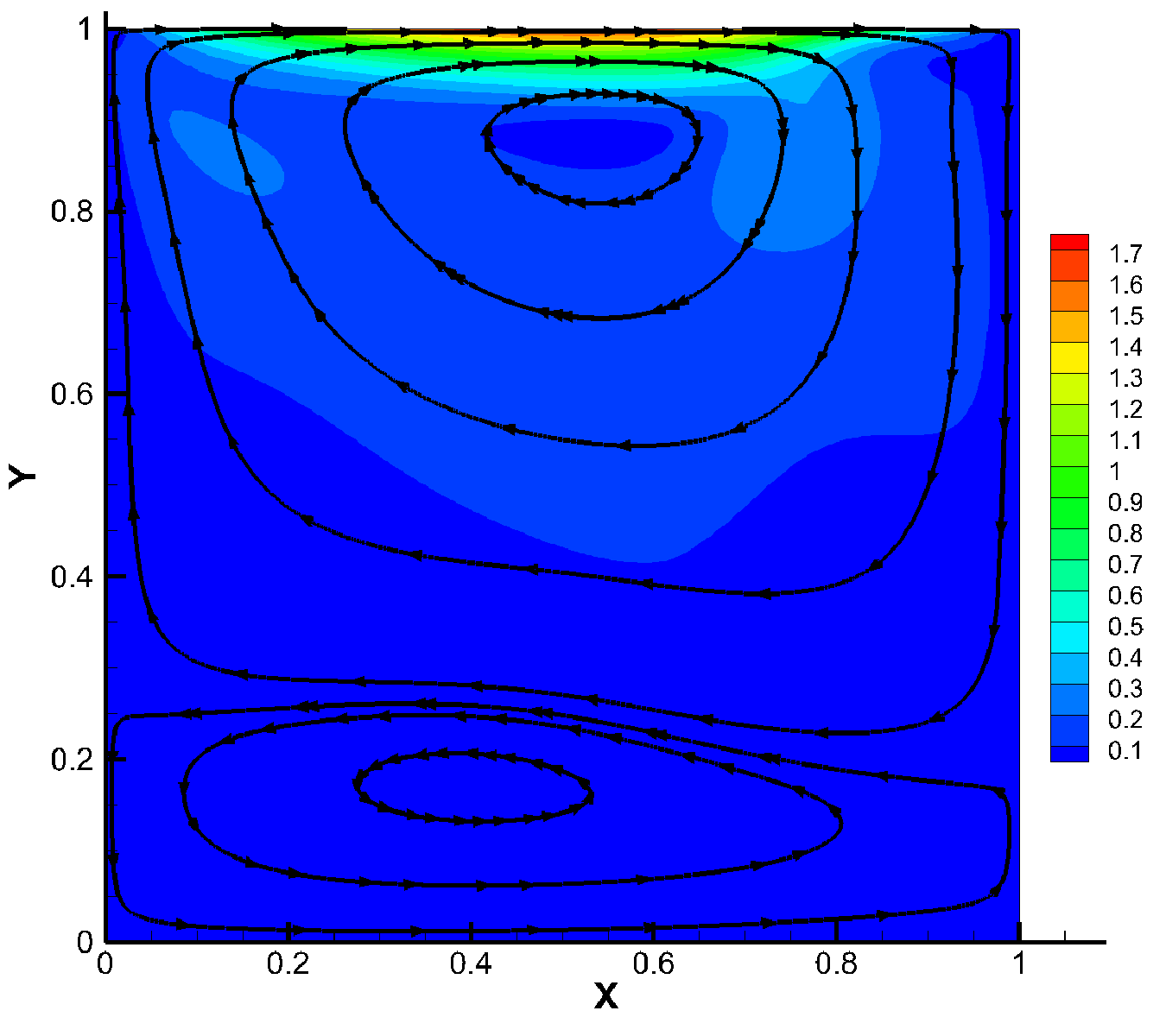}
		\end{minipage}
	}%
	\subfigure[t=3.1]{
		\begin{minipage}[t]{0.2\linewidth}
			\centering
			\includegraphics[width=\textwidth]{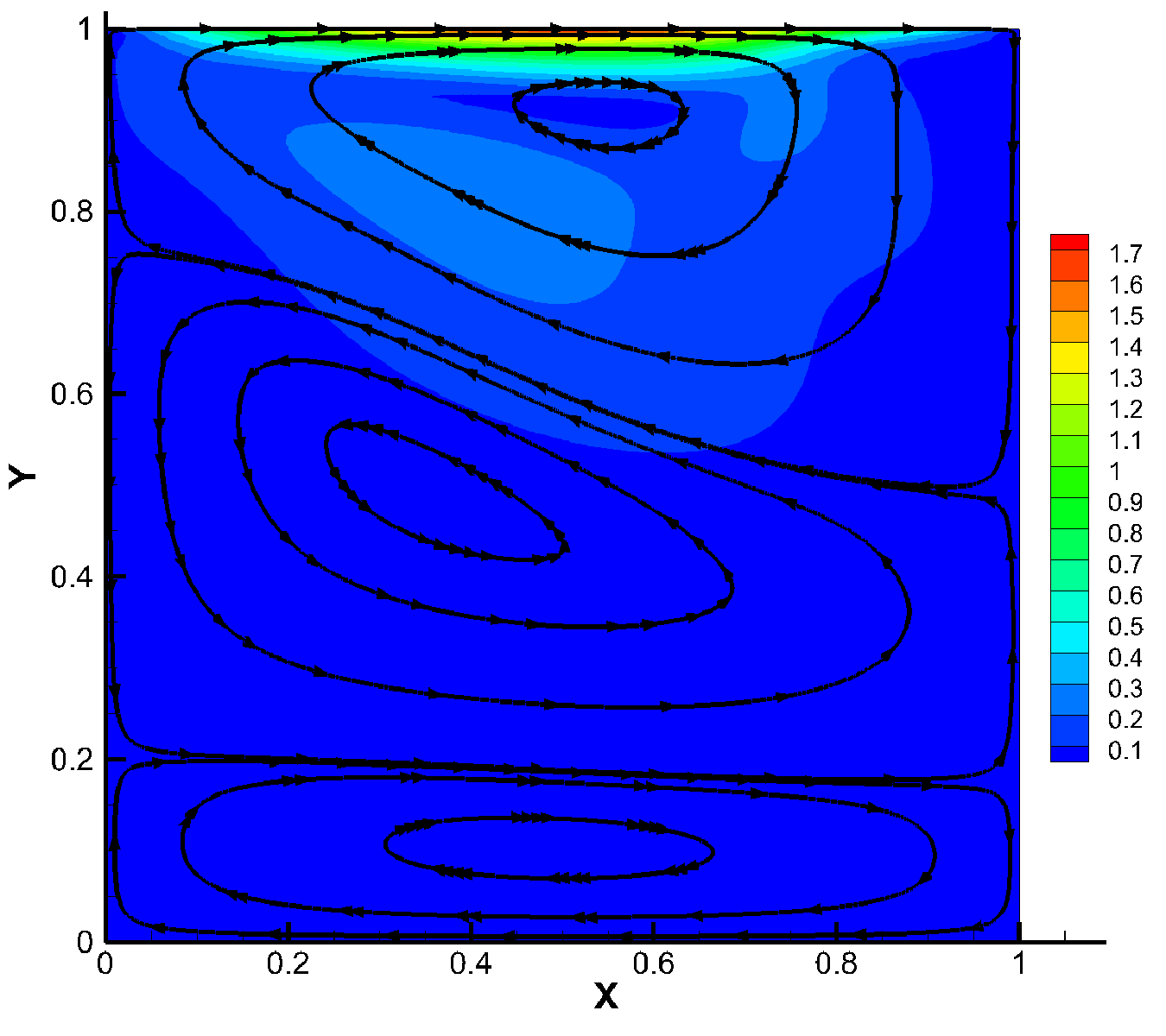}
		\end{minipage}
	}%
	\subfigure[t=5]{
		\begin{minipage}[t]{0.2\linewidth}
			\centering
			\includegraphics[width=\textwidth]{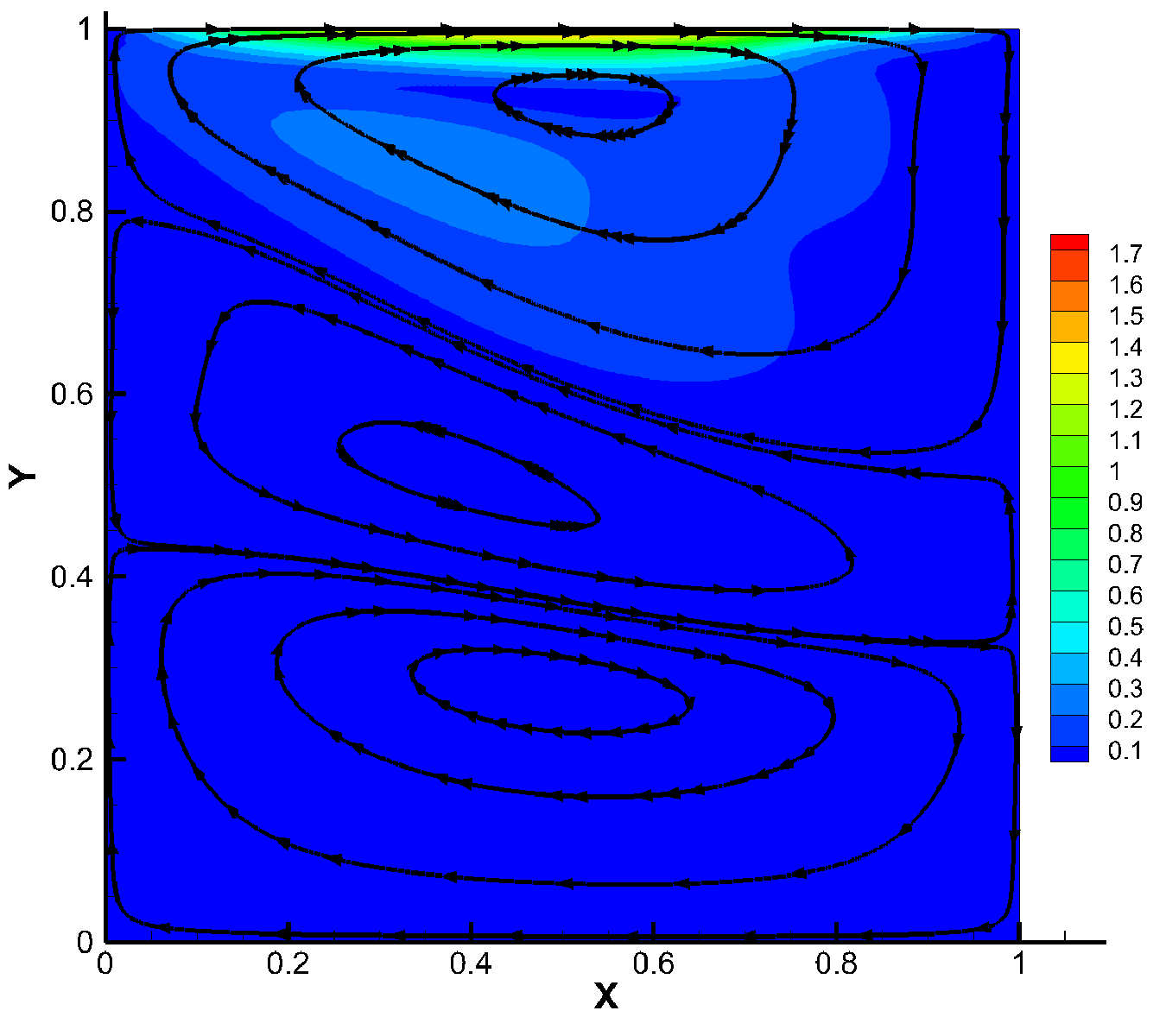}
		\end{minipage}
	}%
   \subfigure[t=9.1]{
		\begin{minipage}[t]{0.2\linewidth}
			\centering
			\includegraphics[width=\textwidth]{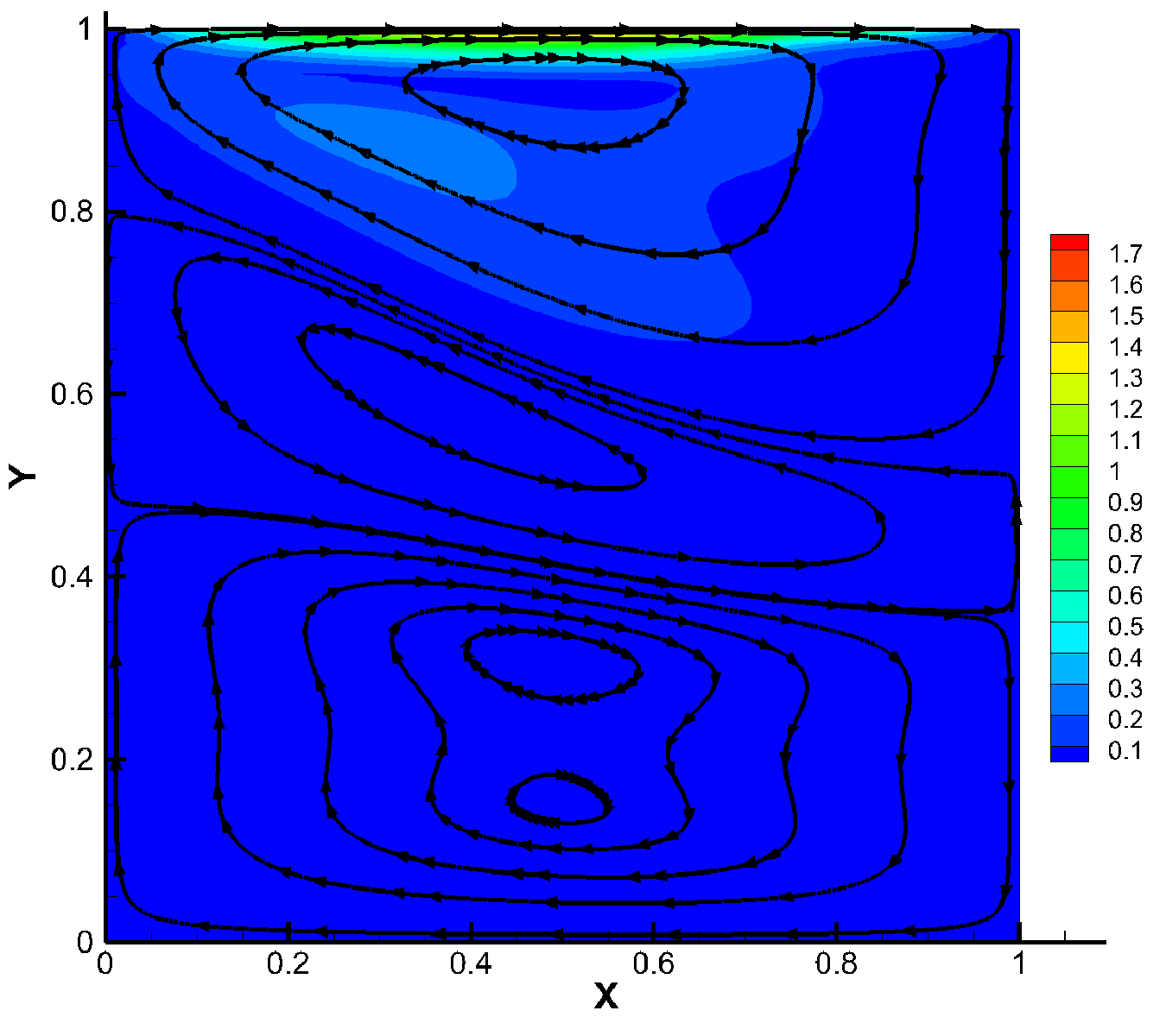}
		\end{minipage}
	}%
	\centering
	\caption{Snapshots of phase field (upper), velocity field (lower) dynamical evolution for lid driven cavity flow with $\mu$=0.6.}
	\label{lid-mu-06}
\end{figure}

\begin{figure}[h]
	\centering
\subfigure[t=0.001]{
		\begin{minipage}[t]{0.2\linewidth}
			\centering
			\includegraphics[width=\textwidth]{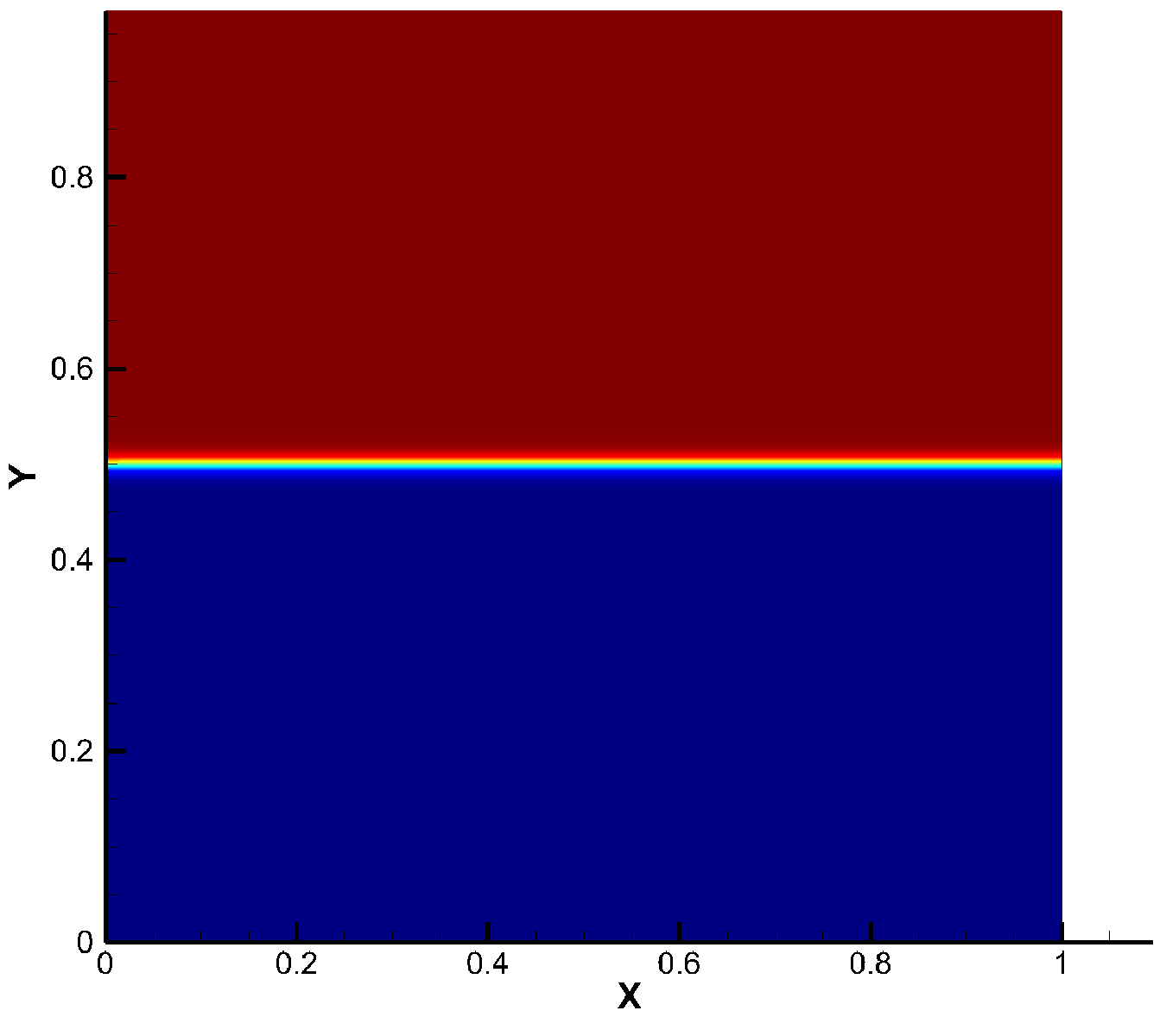}
		\end{minipage}
	}%
	\subfigure[t=2]{
		\begin{minipage}[t]{0.2\linewidth}
			\centering
			\includegraphics[width=\textwidth]{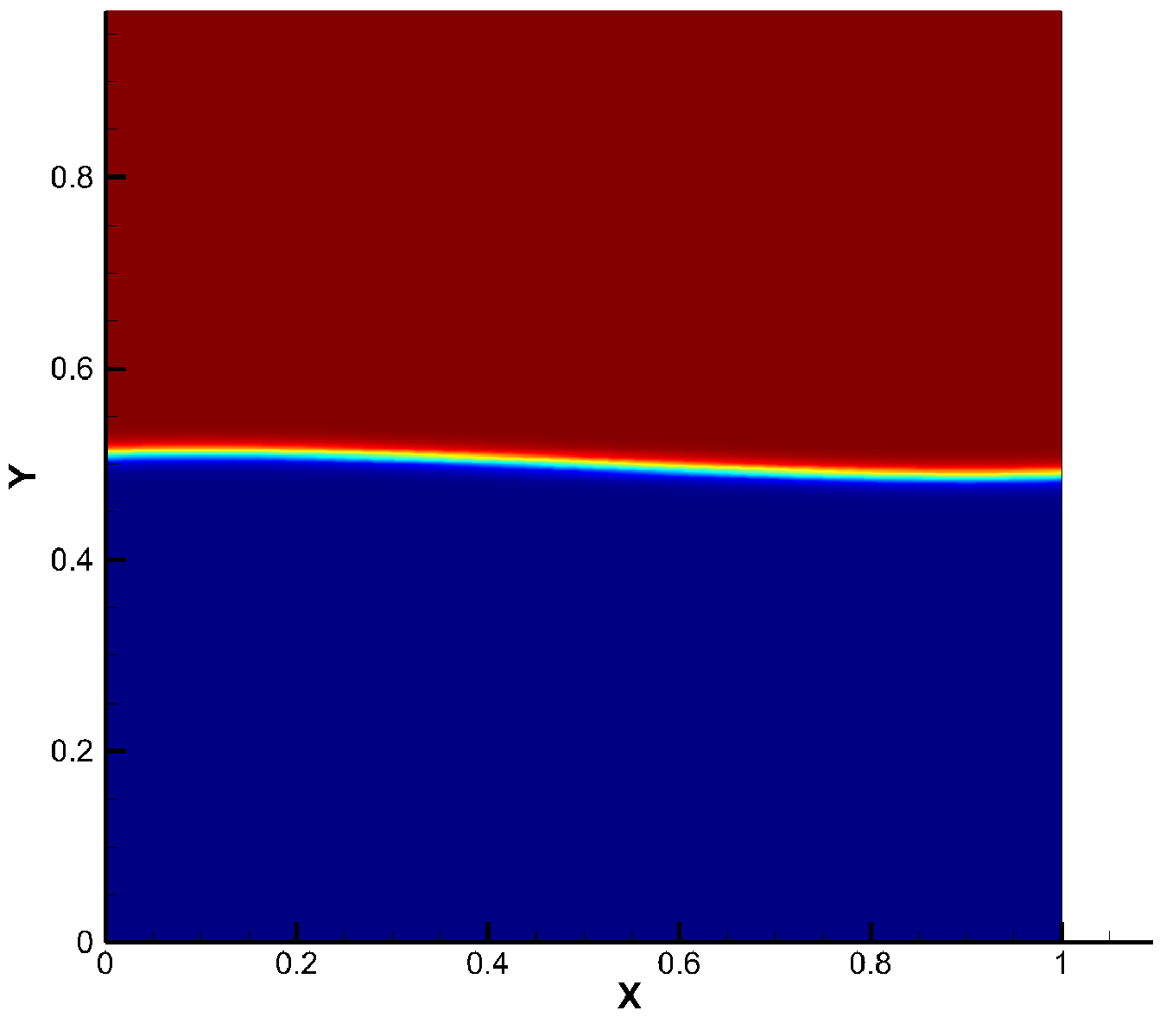}
		\end{minipage}
	}%
	\subfigure[t=3.1]{
		\begin{minipage}[t]{0.2\linewidth}
			\centering
			\includegraphics[width=\textwidth]{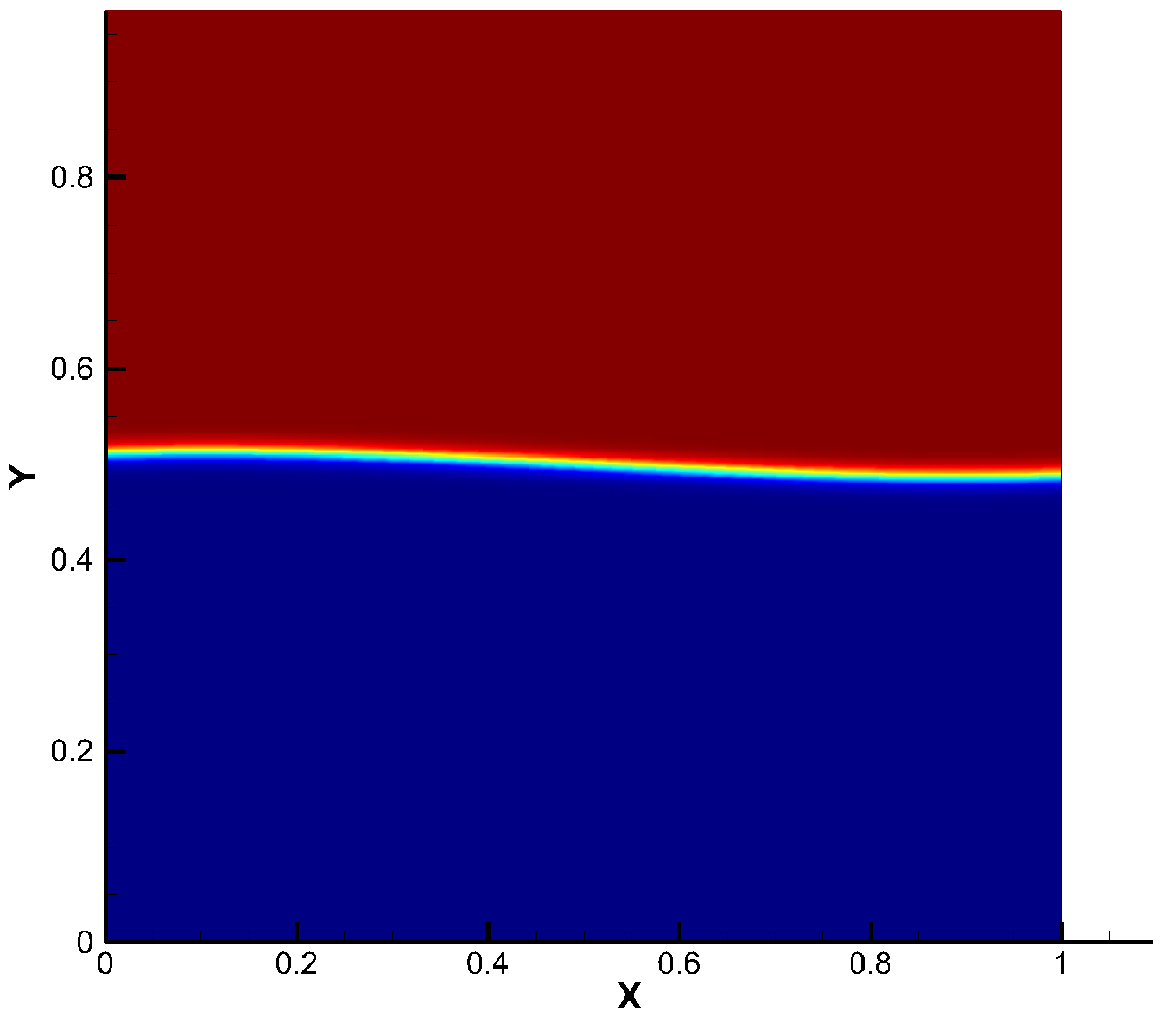}
		\end{minipage}
	}%
	\subfigure[t=5]{
		\begin{minipage}[t]{0.2\linewidth}
			\centering
			\includegraphics[width=\textwidth]{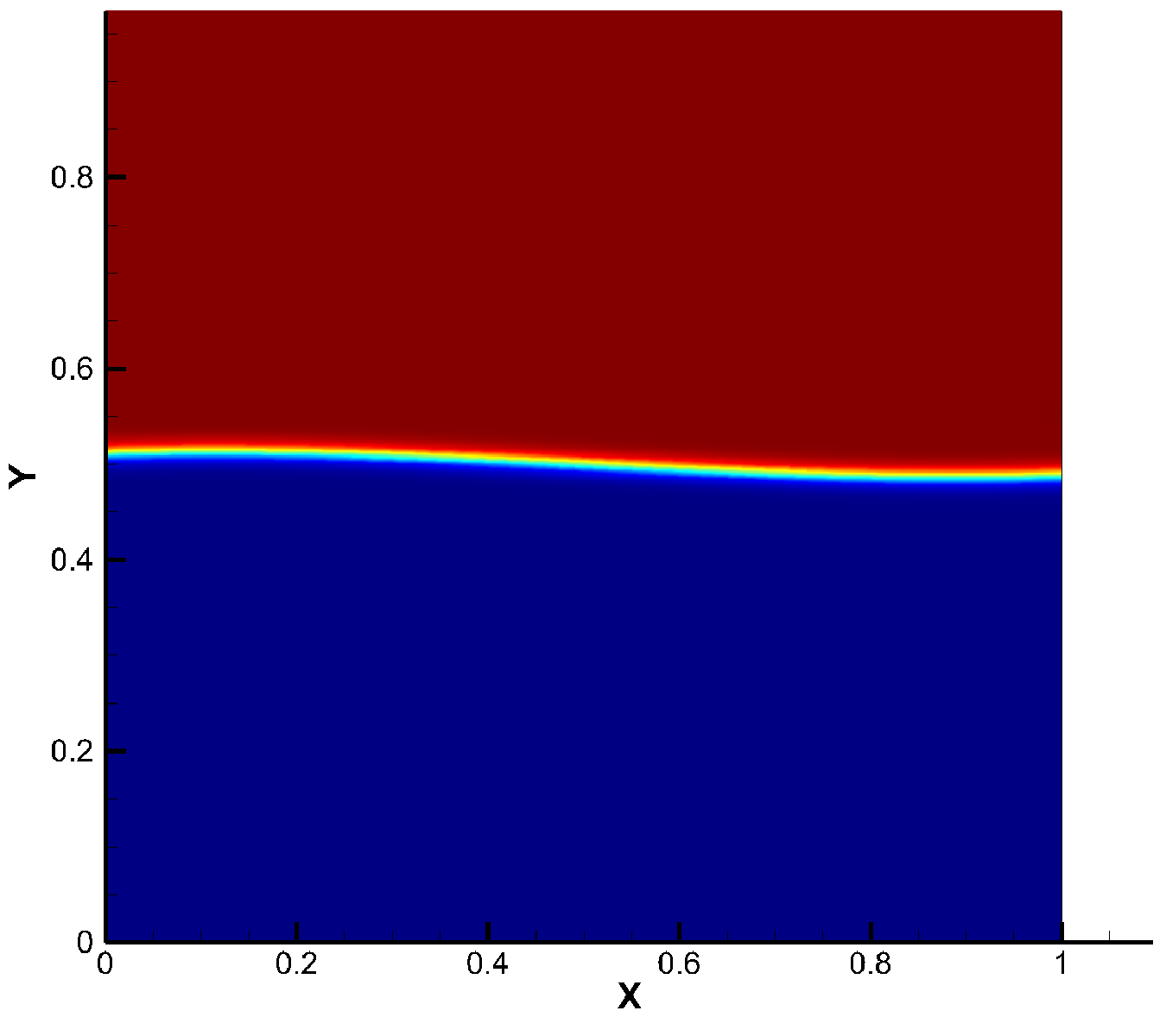}
		\end{minipage}
	}%
   \subfigure[t=9.1]{
		\begin{minipage}[t]{0.2\linewidth}
			\centering
			\includegraphics[width=\textwidth]{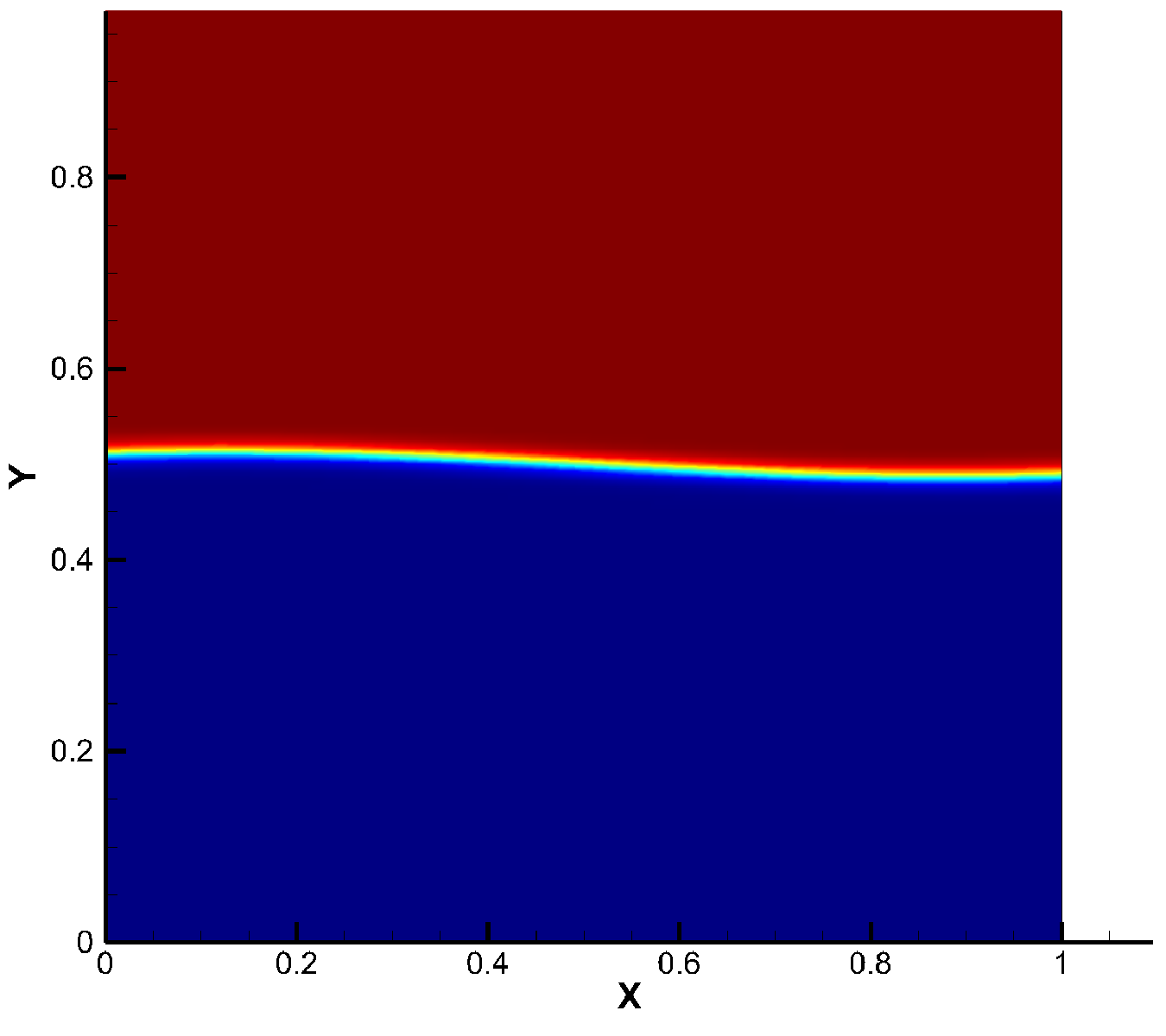}
		\end{minipage}
	}\vspace{5pt}
\subfigure[t=0.001]{
		\begin{minipage}[t]{0.2\linewidth}
			\centering
			\includegraphics[width=\textwidth]{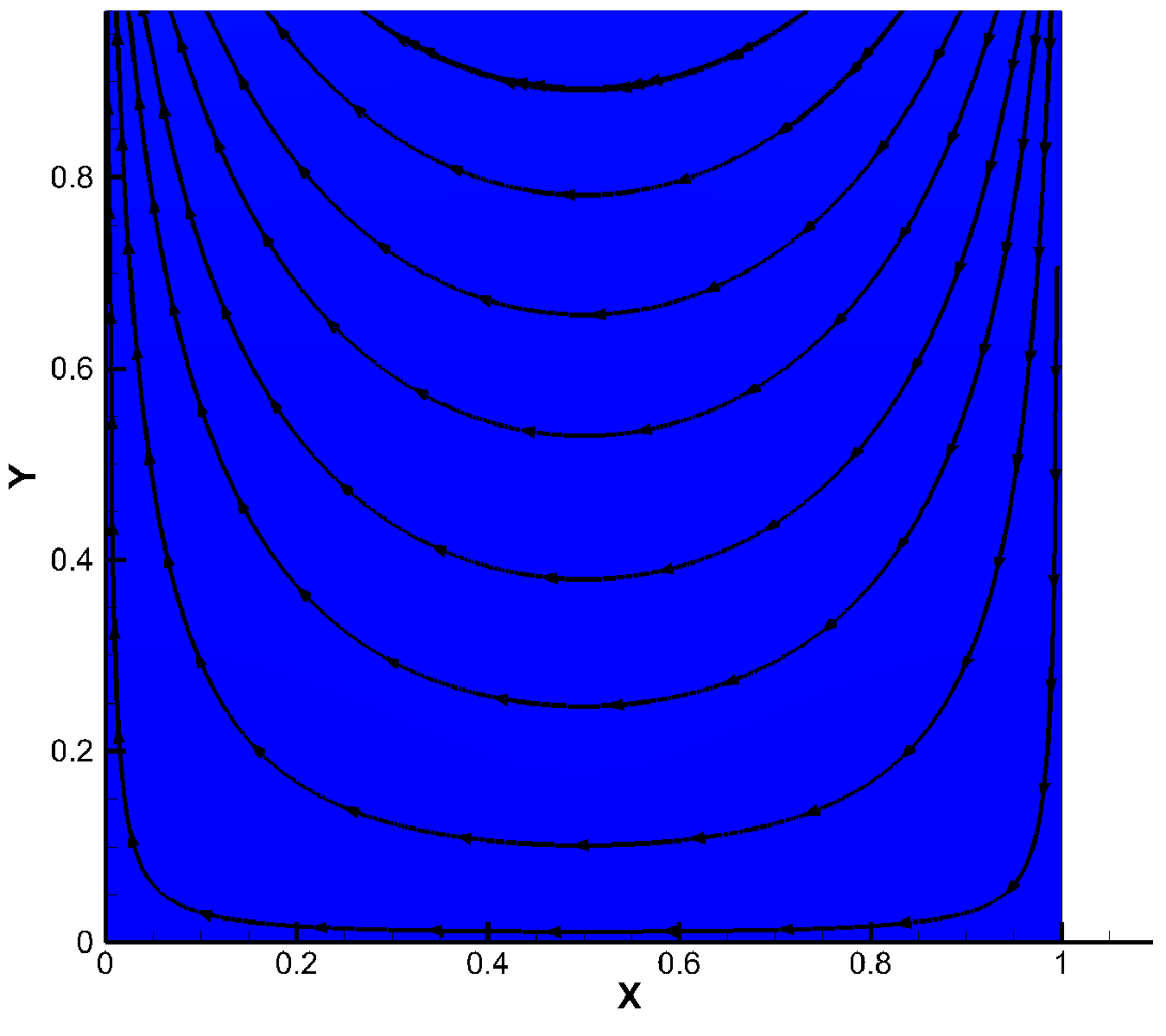}
		\end{minipage}
	}%
	\subfigure[t=2]{
		\begin{minipage}[t]{0.2\linewidth}
			\centering
			\includegraphics[width=\textwidth]{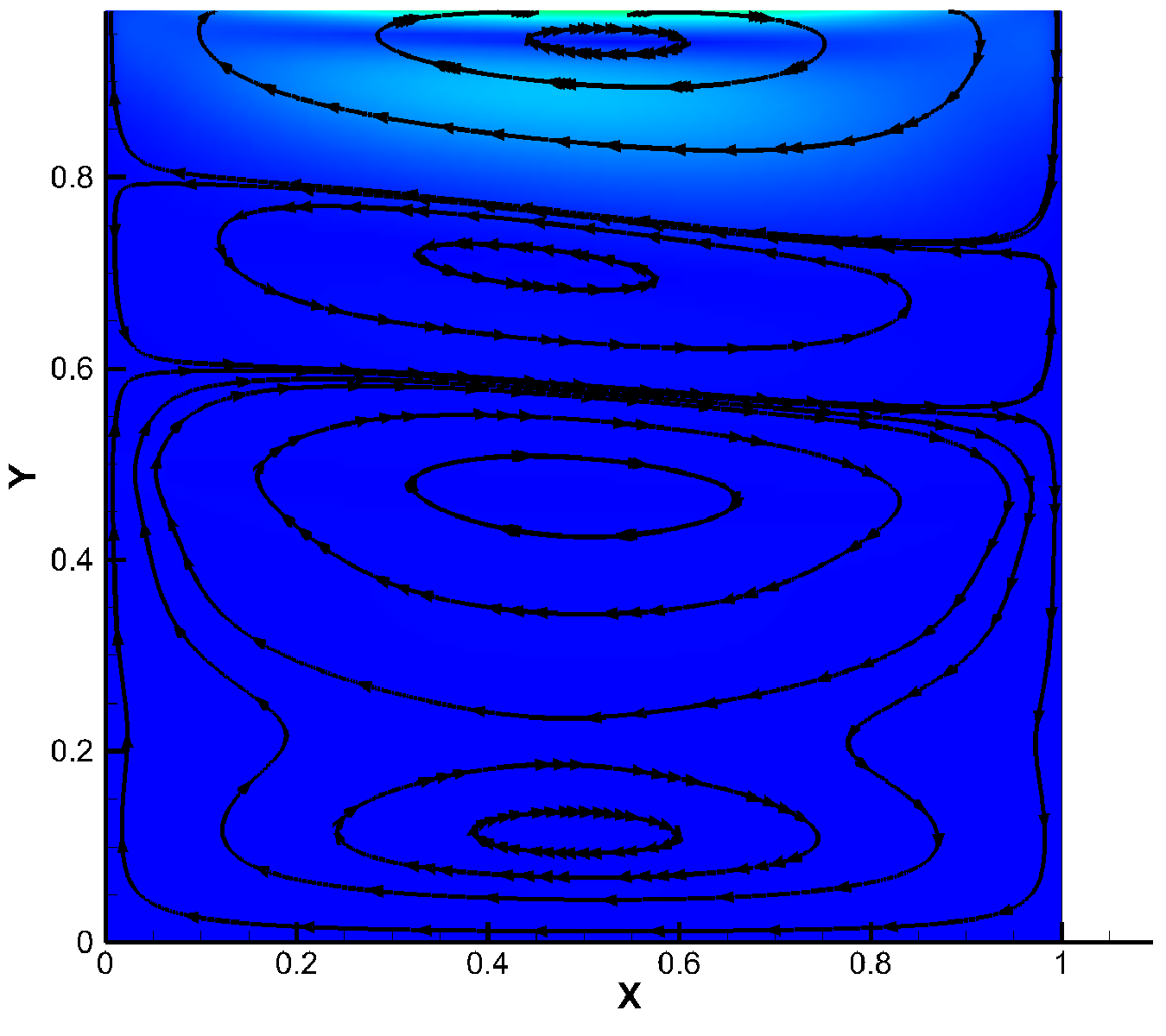}
		\end{minipage}
	}%
	\subfigure[t=3.1]{
		\begin{minipage}[t]{0.2\linewidth}
			\centering
			\includegraphics[width=\textwidth]{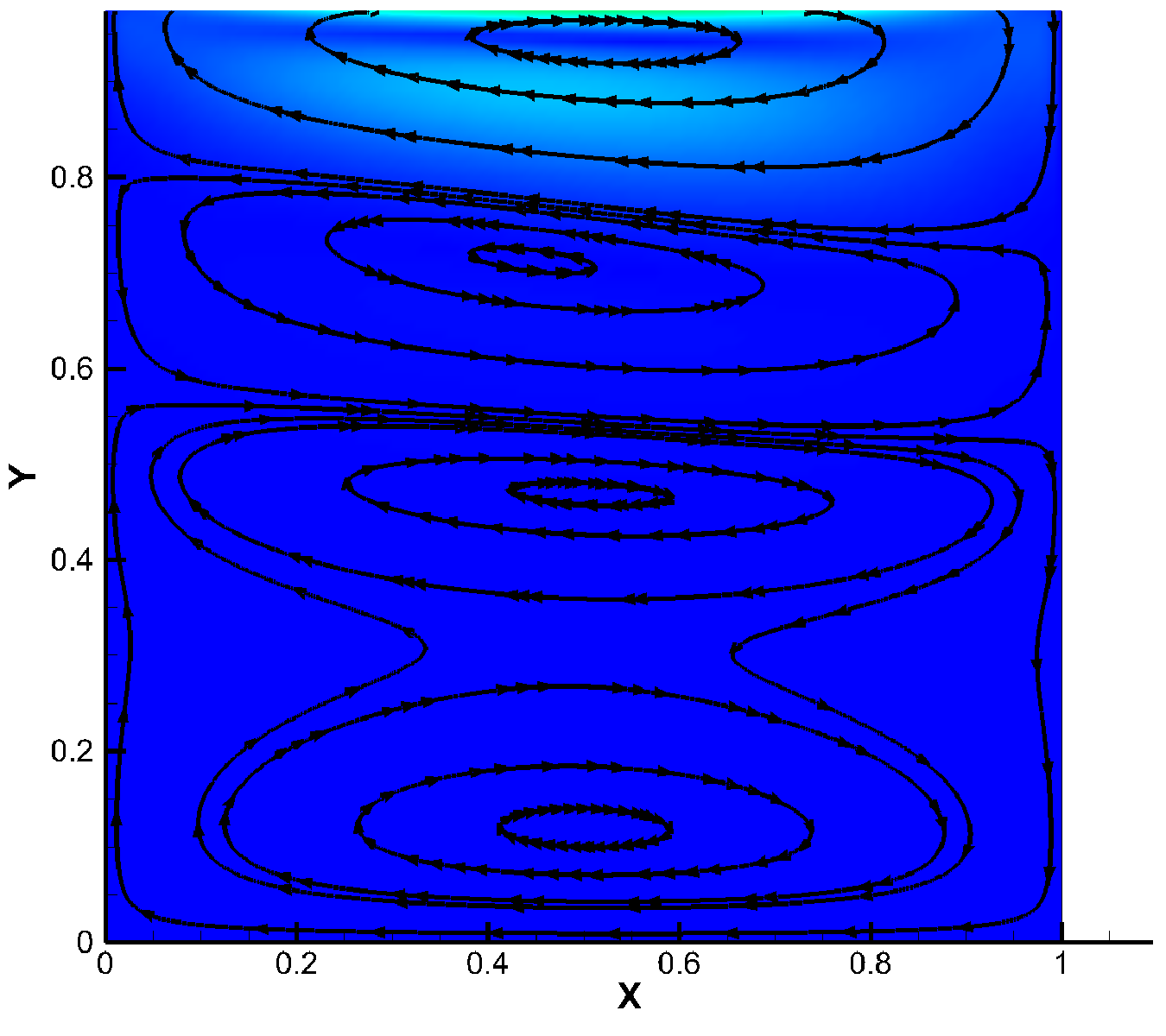}
		\end{minipage}
	}%
	\subfigure[t=5]{
		\begin{minipage}[t]{0.2\linewidth}
			\centering
			\includegraphics[width=\textwidth]{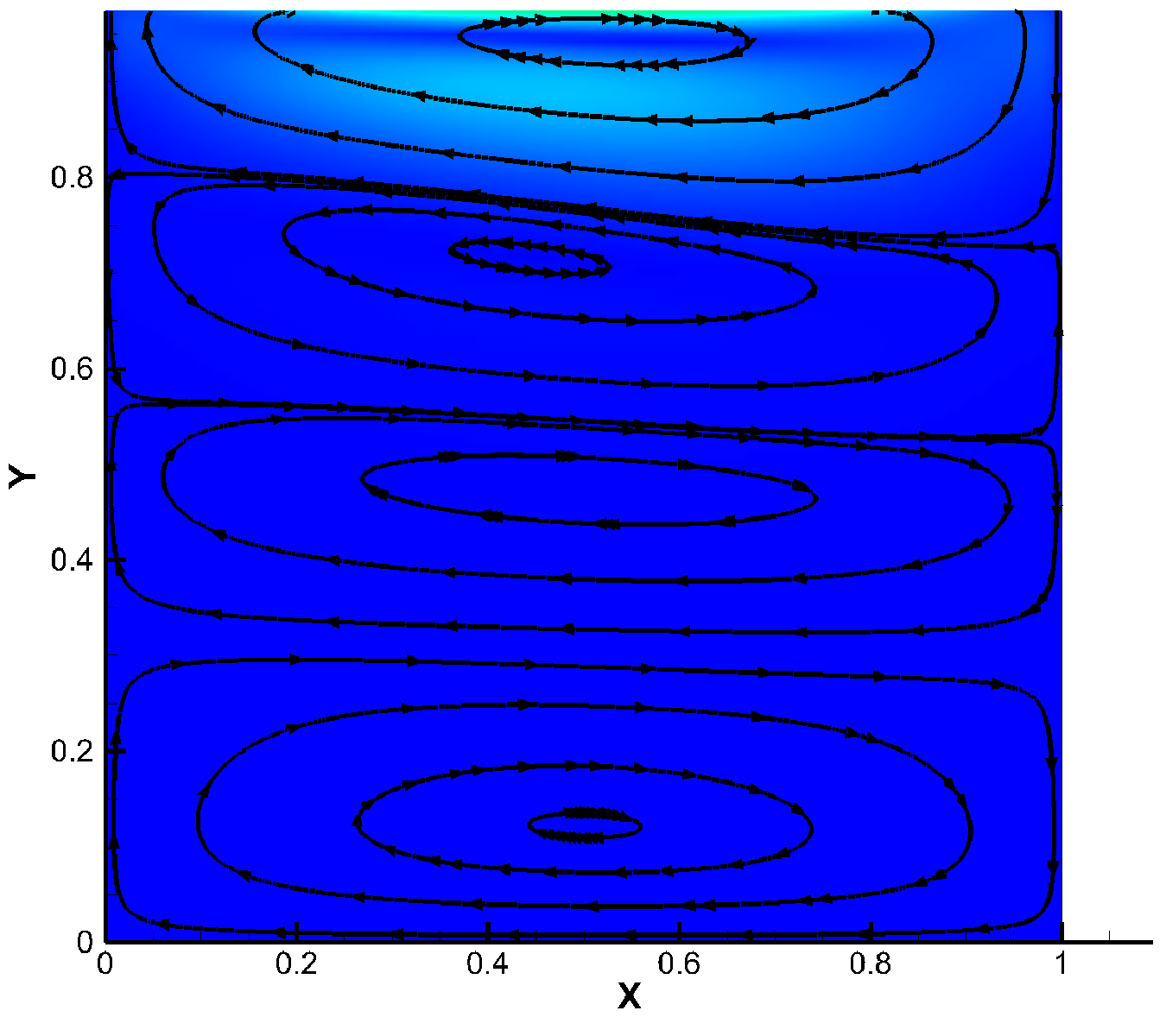}
		\end{minipage}
	}%
   \subfigure[t=9.1]{
		\begin{minipage}[t]{0.2\linewidth}
			\centering
			\includegraphics[width=\textwidth]{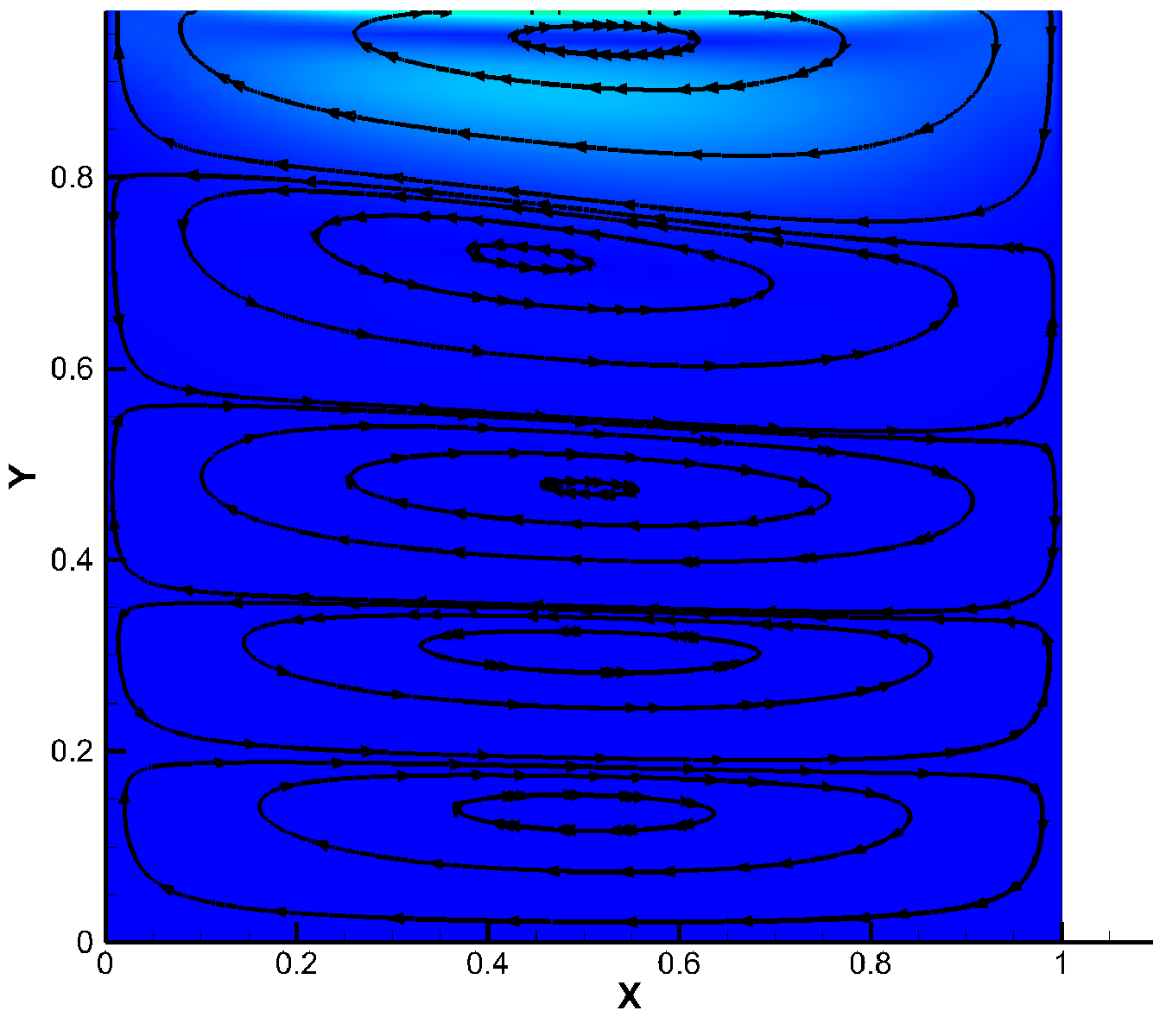}
		\end{minipage}
	}%
	\centering
	\caption{Snapshots of phase field (upper), velocity field (lower) dynamical evolution for lid driven cavity flow with $\mu$=0.1.}
	\label{lid-mu-01}
\end{figure}

In this subsection, we consider the well-known lid-driven cavity flow as a benchmark problem \cite{2008Finite, hintermuller2013adaptive} in a unit square $\Omega=[0,1]^{2}$.  We give the initial condition of the phase field as
\begin{equation*}
\phi_{0}=\tanh (100(y-0.5)).
\end{equation*}
The boundary conditions are set as  
\begin{equation*}
\frac{\partial\phi}{\partial \n}\bigg|_{\partial\Omega}=0, 
	\quad 
\frac{\partial\omega}{\partial \n}\bigg|_{\partial\Omega}=0, 
	\quad
\u\big|_{y=1}=(7x(x-1), 0)^\top, 
	\quad
\B\big|_{\partial\Omega}=  (1, 0)^\top,
\end{equation*}
and the velocity field $\u$ has a no-slip boundary condition on the other walls. Hereafter, we consider the unmatched mobility, viscosities, and electric conductivities of the two fluids, specifically, 
\begin{equation*}
M:=M(\phi)=\frac{M_{2}-M_{1}}{2}\phi+\frac{M_{2}+M_{1}}{2},\quad \nu:=\nu(\phi)=\frac{\nu_{2}-\nu_{1}}{2}\phi+\frac{\nu_{2}+\nu_{1}}{2},\quad \sigma:=\sigma(\phi)=\frac{\sigma_{2}-\sigma_{1}}{2}\phi+\frac{\sigma_{2}+\sigma_{1}}{2}.
\end{equation*}
We set the mesh size $h=1/120$, the time step $\Delta t=1/1000$, and the parameters as  
\begin{align}
&\gamma =1/120,	\quad M_1=M_2=0.12,	\quad \nu_{1}=1/1000,\quad \nu_{2}=1/100,\quad \lambda=1/1000,\quad \sigma_{1}=50, \quad \sigma_{2}=150. \label{3n2} 
\end{align}
Different values of $\mu$ imply different strengths of Lorentz forces. To clarify the effects of the magnetic induction field, we consider the parameter set defined in (\ref{3n2}) with  $\mu=2$, $0.6$, and $0.1$, which  are imposed in the numerical scheme  (\ref{scheme1})-(\ref{scheme5}). The results are as follows:
\begin{itemize}
 \item  The numerical results  for $\mu=2$  are displayed in  Figure \ref{lid-mu-2}. The applied boundary velocity pushes the free interface towards the upper region of the cavity, and as time progresses, a concave finger-like interface emerges in the cavity's left section.  A small velocity vortex forms in the lower right corner of the cavity, followed by the emergence of another in the lower left corner. Over time, both vortices gradually diminish in size. Similar numerical results are referred to \cite{2008Finite, boyer2002theoretical}.

 \item The numerical results  for $\mu=0.6$ and $0.1$  are presented in Figures \ref{lid-mu-06} and \ref{lid-mu-01}, respectively. The Lorentz forces become larger as $\mu$ decreases from 2 to 0.6 and then to 0.1. Compared with Figures \ref{lid-mu-06}-\ref{lid-mu-01}, we can observe that when the Lorentz force increases, the phase field evolution rate decreases. These simulations indicate that the large Lorentz force inhibits the stretching of the diffuse interface. Moreover, the velocity field is significantly influenced by the Lorentz force compared with Figures \ref{lid-mu-2}-\ref{lid-mu-01}.  The main vortex of velocity persists, and two smaller vortices will gradually develop on both sides at the lower part in Figure \ref{lid-mu-2}; The main velocity vortex gradually breaks down into several smaller vortices, and the vortices at the bottom gradually grow larger in Figure \ref{lid-mu-06}. The primary vortex in Figure \ref{lid-mu-01} gradually evolves into five small vortices that are essentially uniform.  The numerical examples are similar to the work in \cite{2022Highly}.

\end{itemize}

\subsection{2D/3D Kelvin-Helmholtz instability}
 
The Kelvin-Helmholtz (K-H) instability  is a common fluid instability caused by the velocity difference at the fluid interface \cite{Lee2015Two, 2023A}. Since the K-H instability has wide applications in natural and industrial fields, many researchers implemented numerical simulations of K-H instability  in recent years \cite{Yshin2018vortex}. In this simulation, we test the 2D/3D  K-H instability, where the parameter values are set to
\begin{equation*}
\gamma=1/100,	\quad 	M_{1}=M_{2}=1/100,	\quad 	\nu_{1}=\nu_{2}=1/1000,	\quad	\mu=1,	\quad \lambda=1/10000,	\quad \sigma_{1}=\sigma_{2}=1. 
\end{equation*}
We consider appropriate mesh sizes and time steps to effectively capture the dynamics of the interface in the computational domain $\Omega=[0, 1]^{d}$. The periodic boundary conditions for all variables are applied to the boundaries at $x=0$ and $x=1$.

\subsubsection{2D Kelvin-Helmholtz instability}

This example  illustrates  the dynamics of a sinusoidal perturbation at the interface between two fluids, characterized by a single mode of perturbation. We set the mesh size $h=1/150$, time step $\Delta t=1/1000$, and the following initial values:
\begin{eqnarray*}\label{K-H1}
\left\{
\begin{aligned}
\phi_{0}&=\tanh(\frac{y-0.5-0.01\sin(2\pi x)}{\sqrt{2}\gamma}),\\
\u_{0}&=\big(\tanh(\frac{y-0.5-0.01\sin(2\pi x)}{\sqrt{2}\gamma}), 0\big)^\top,\\
\B_{0}&=\big(1, 0\big)^\top.
\end{aligned}
\right.
\end{eqnarray*}
The boundary conditions for $\B$ at the top ($y=1$) and bottom ($y=0$) are given by $(-1, 0)^\top$, and  the  vertical component of $\u$ is $u_{2}=0$.
 
\begin{figure}[h]
	\centering
\subfigure[t=0.001]{
		\begin{minipage}[t]{0.25\linewidth}
			\centering
			\includegraphics[width=\textwidth]{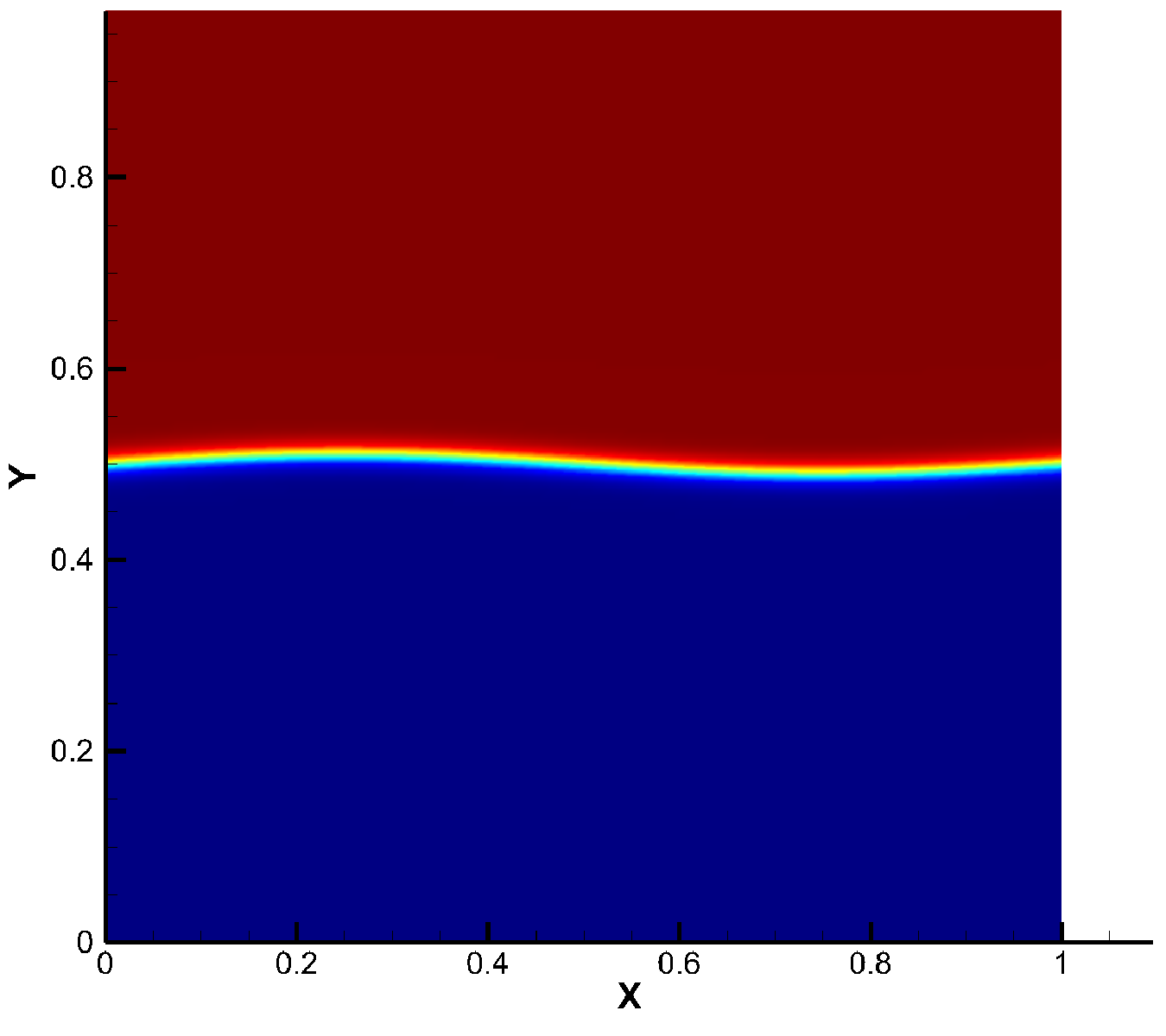}
		\end{minipage}
	}%
	\subfigure[t=0.6]{
		\begin{minipage}[t]{0.25\linewidth}
			\centering
			\includegraphics[width=\textwidth]{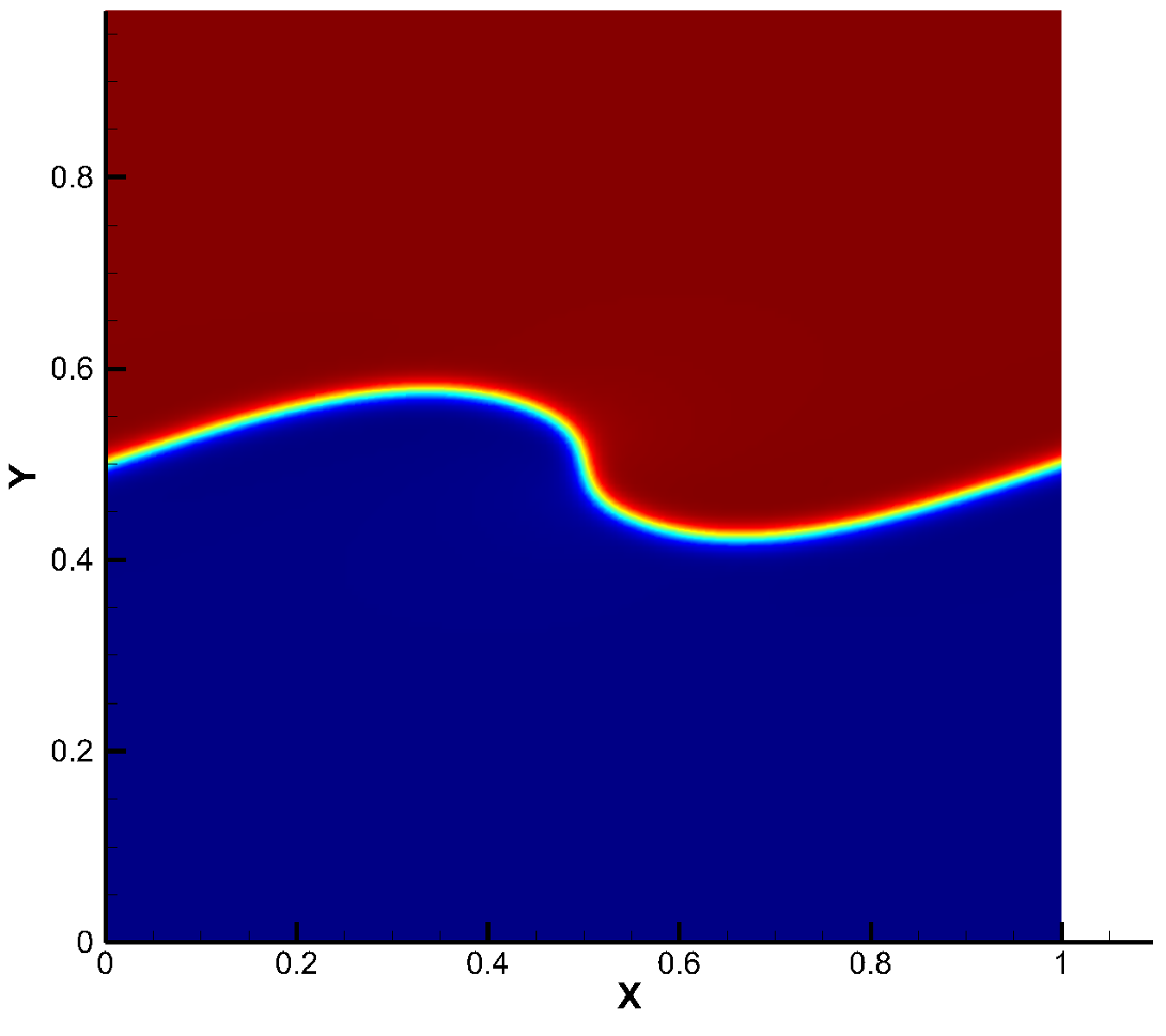}
		\end{minipage}
	}%
	\subfigure[t=0.85]{
		\begin{minipage}[t]{0.25\linewidth}
			\centering
			\includegraphics[width=\textwidth]{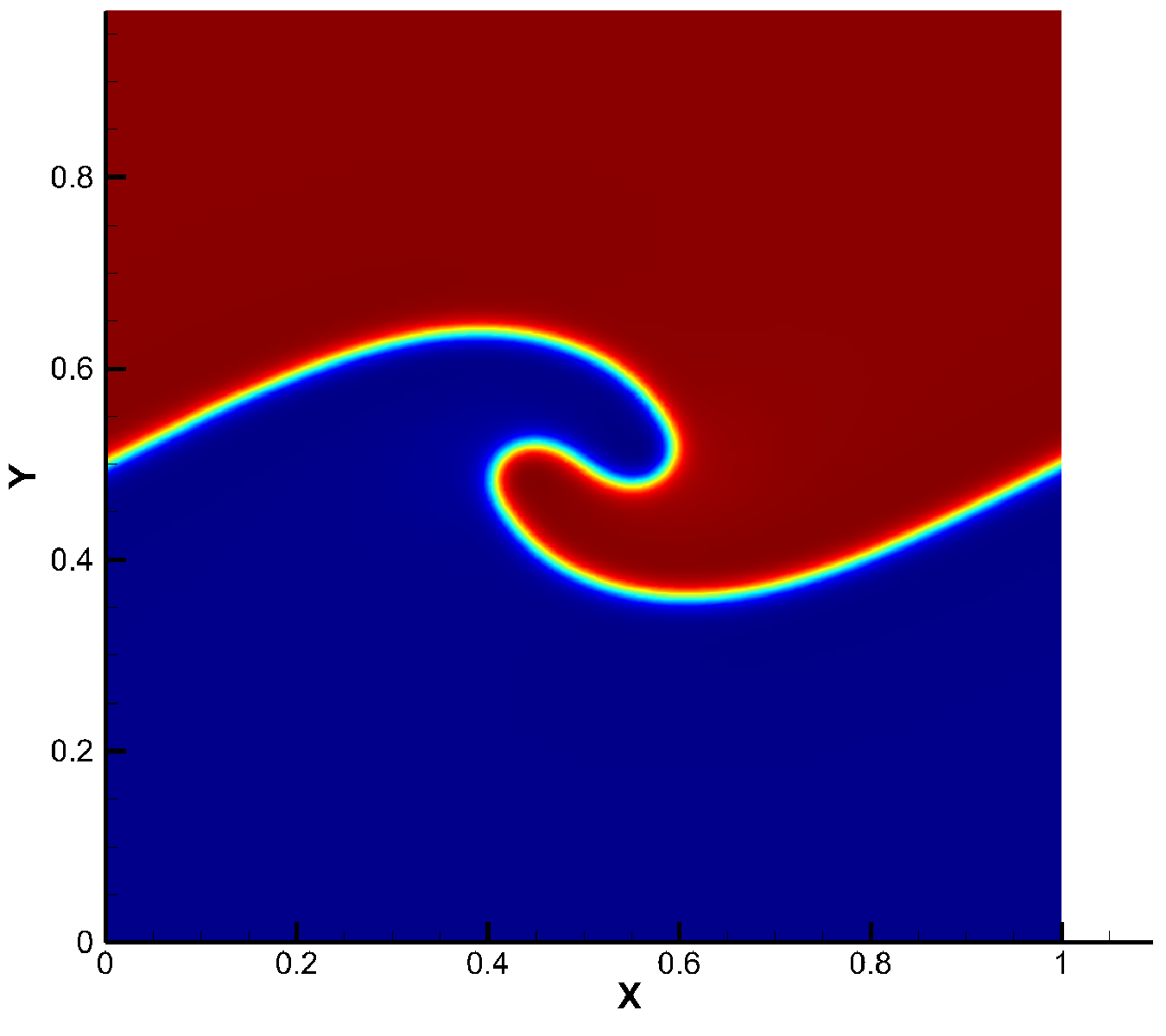}
		\end{minipage}
	}%
\subfigure[t=1]{
		\begin{minipage}[t]{0.25\linewidth}
			\centering
			\includegraphics[width=\textwidth]{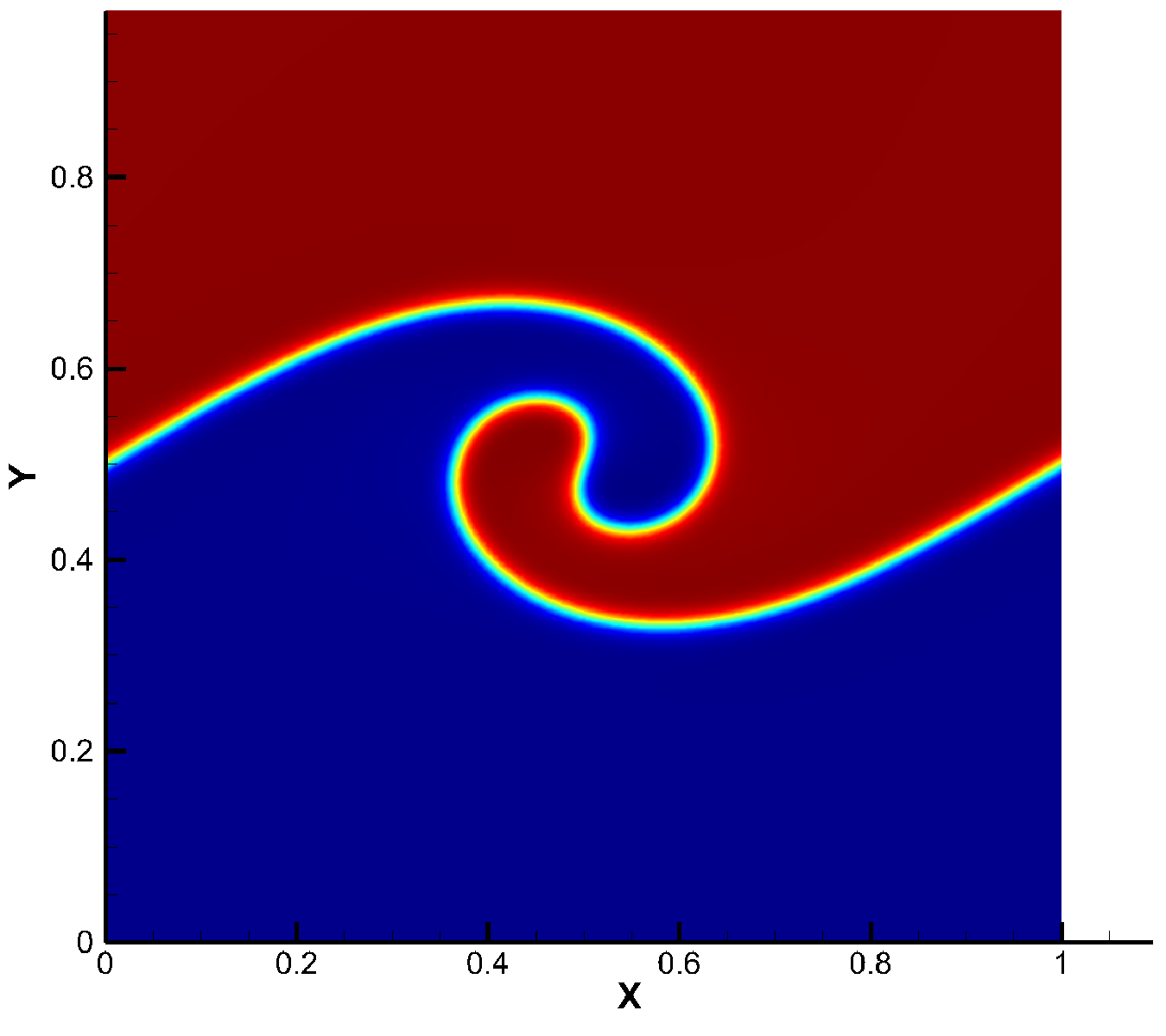}
		\end{minipage}
	}%
\\
	\subfigure[t=1.1]{
		\begin{minipage}[t]{0.25\linewidth}
			\centering
			\includegraphics[width=\textwidth]{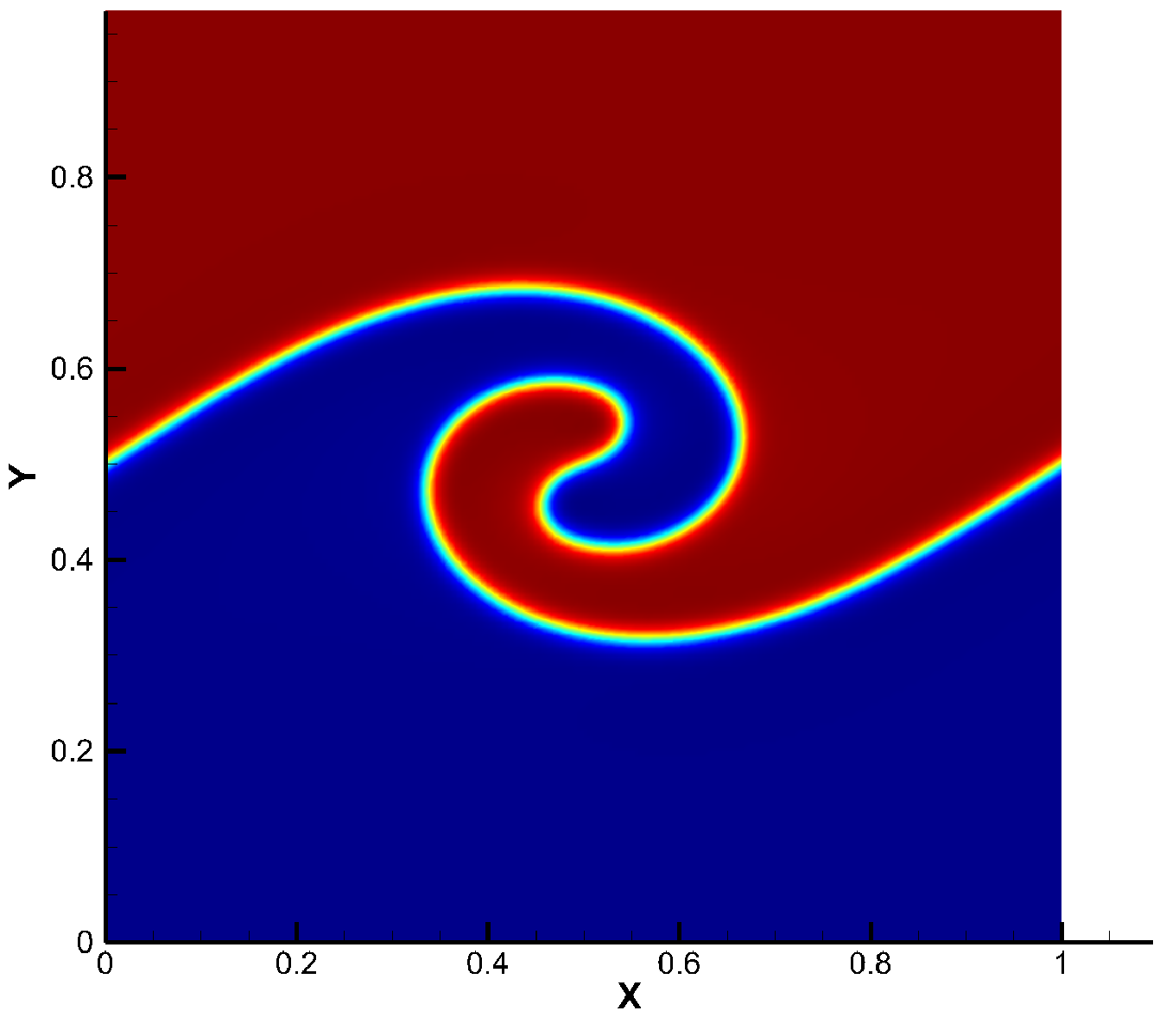}
		\end{minipage}
	}%
\subfigure[t=1.2]{
		\begin{minipage}[t]{0.25\linewidth}
			\centering
			\includegraphics[width=\textwidth]{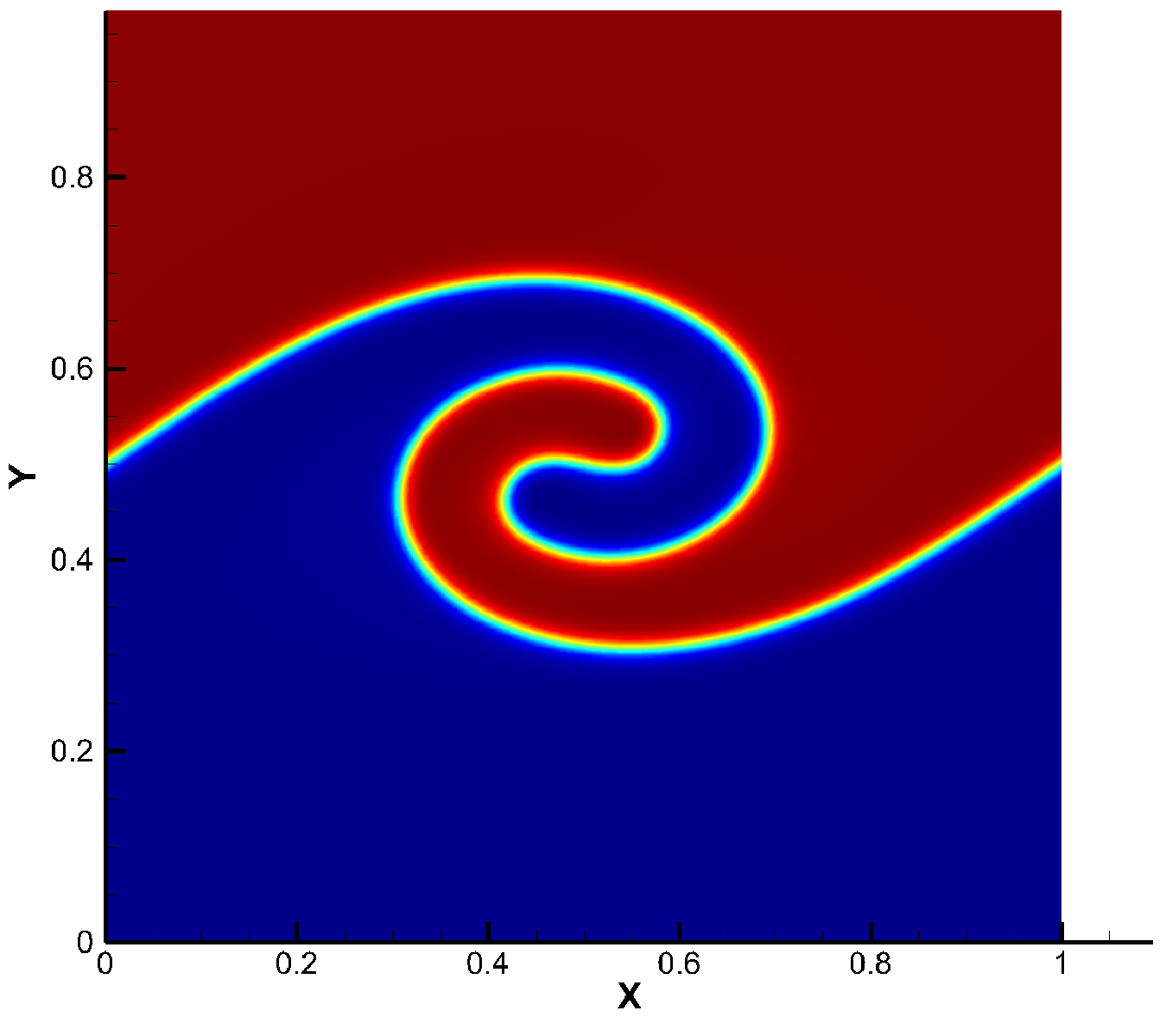}
		\end{minipage}
	}%
   \subfigure[t=1.4]{
		\begin{minipage}[t]{0.25\linewidth}
			\centering
			\includegraphics[width=\textwidth]{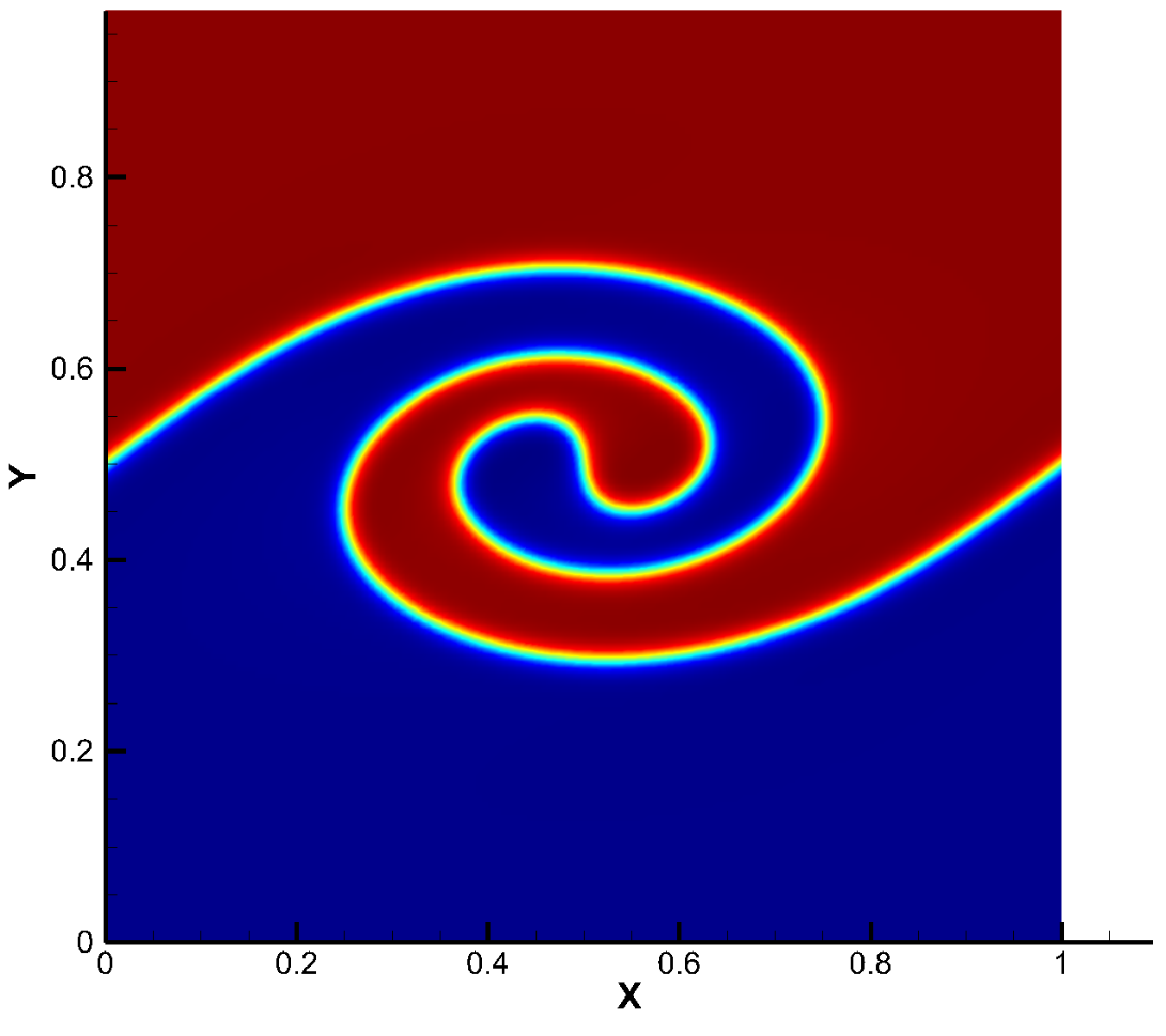}
		\end{minipage}
	}%
\subfigure[t=1.6]{
		\begin{minipage}[t]{0.25\linewidth}
			\centering
			\includegraphics[width=\textwidth]{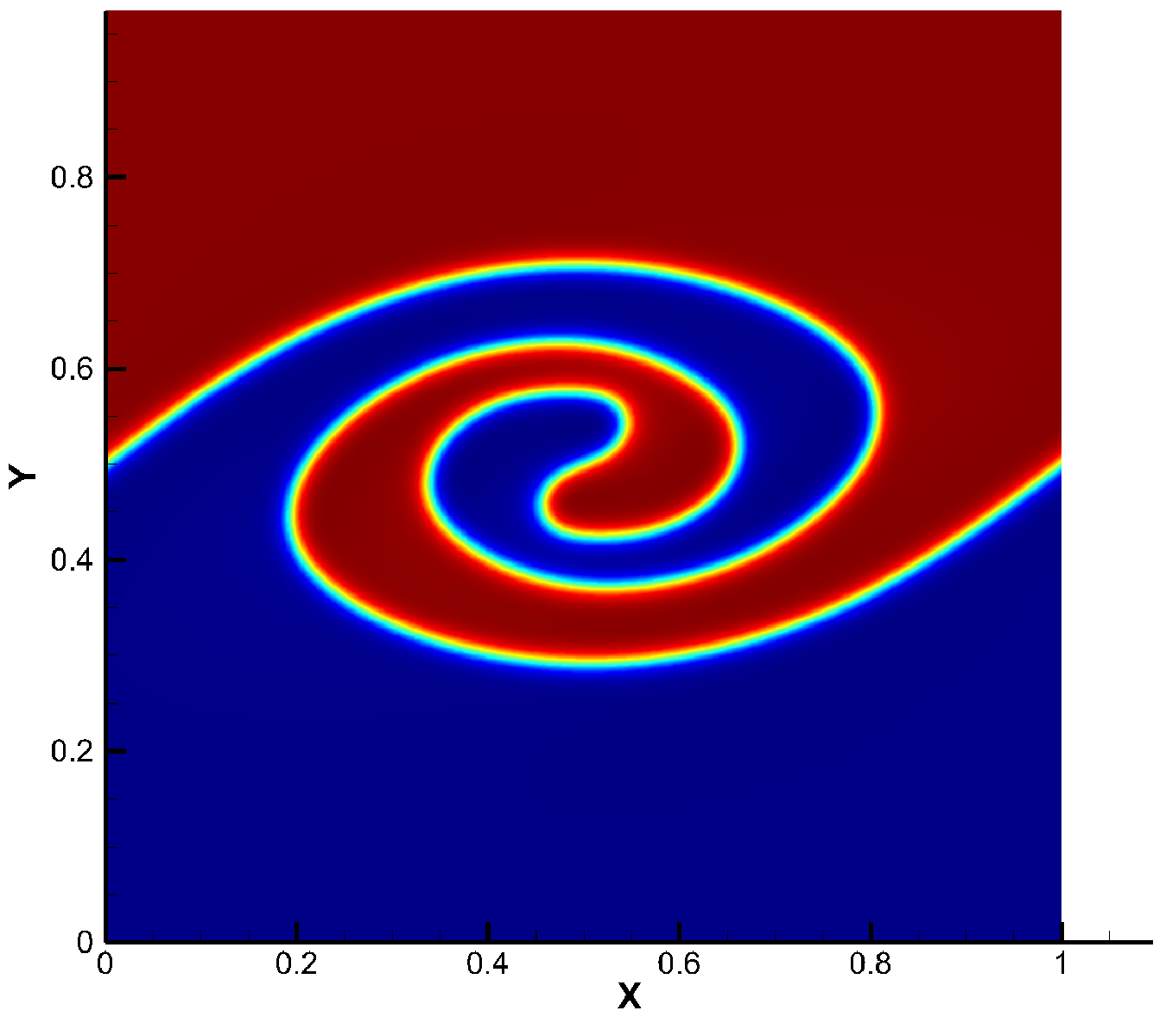}
		\end{minipage}
	}%
	\centering
	\caption{Temporal evolution of the phase field perturbed sinusoidal at t=0.001 (a), 0.6 (b), 0.85 (c), 1 (d), 1.1 (e), 1.2 (f), 1.4 (g), 1.6 (h) for case I.}
	\label{KH-phase}
\end{figure}

Figure \ref{KH-phase} shows the evolution of the phase field with a single-mode sinusoidal interface perturbation at different times. The interface undergoes a rolling up at the center of the domain at $t=0.6$. The rolling up of the interface forms a spiral shape at a later time, specifically showing  the characteristic features of K-H instability, as depicted in Figure \ref{KH-phase}.

\begin{figure}[h]
	\centering
\subfigure[t=0.001]{
		\begin{minipage}[t]{0.25\linewidth}
			\centering
			\includegraphics[width=\textwidth]{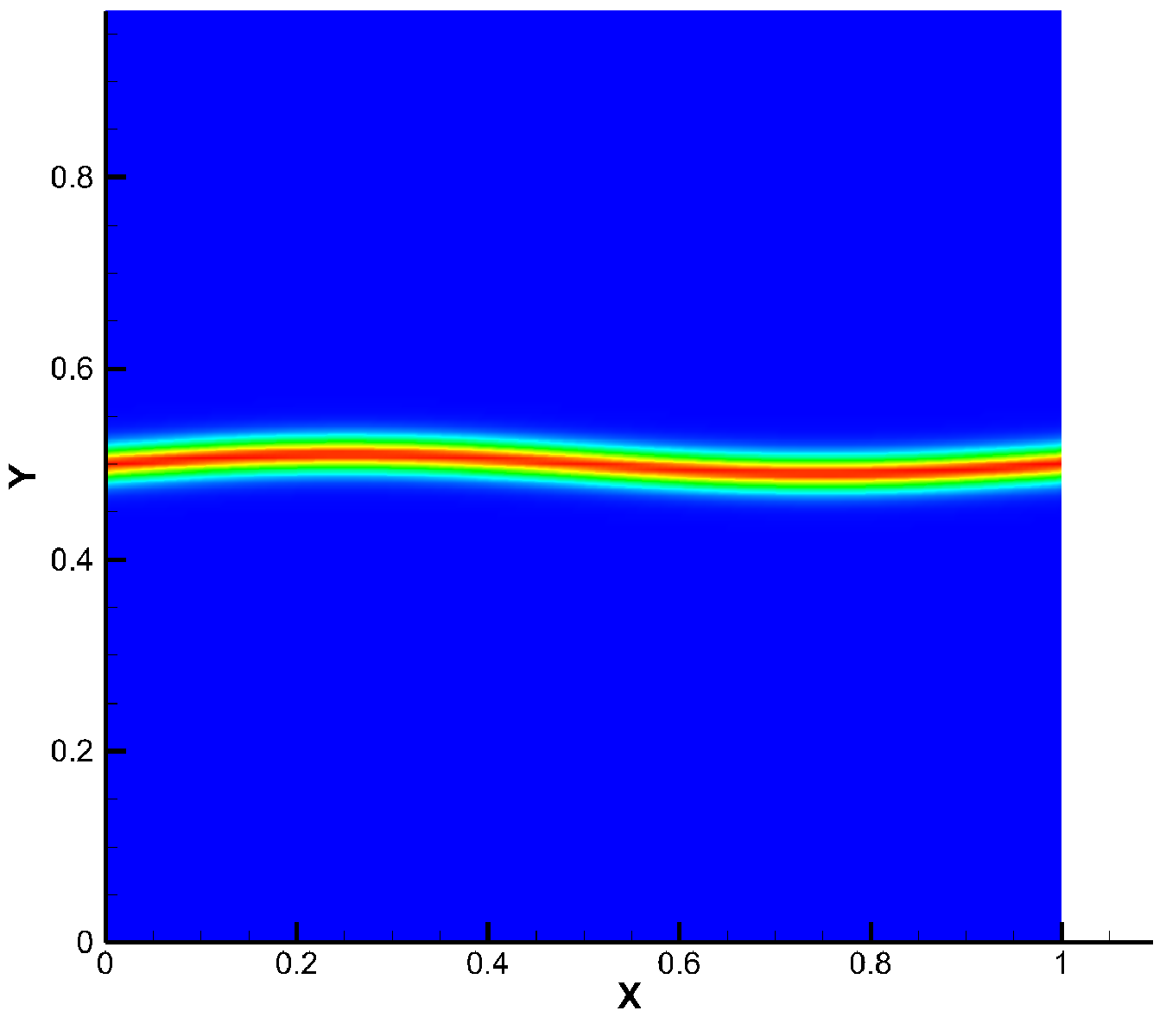}
		\end{minipage}
	}%
	\subfigure[t=0.6]{
		\begin{minipage}[t]{0.25\linewidth}
			\centering
			\includegraphics[width=\textwidth]{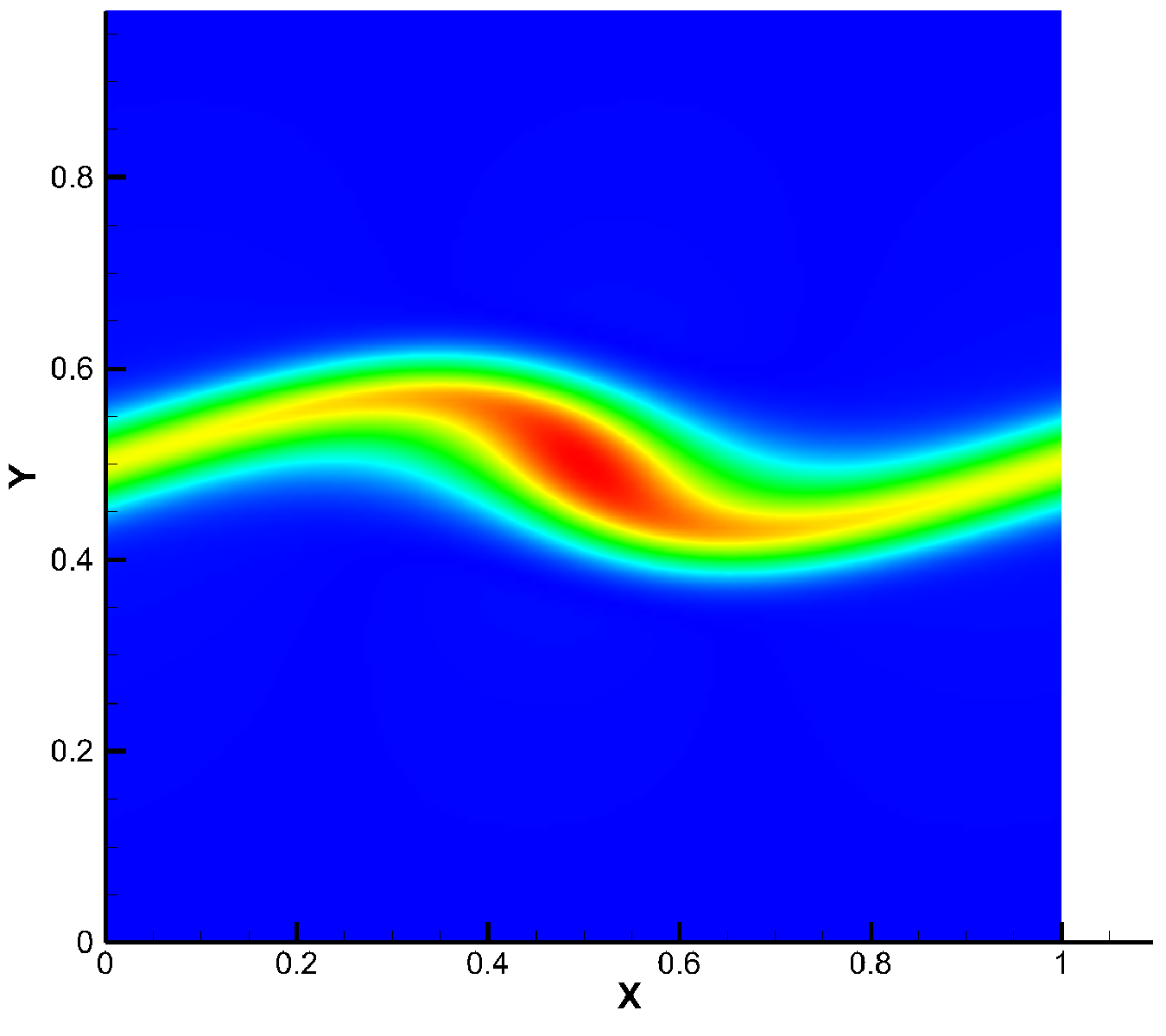}
		\end{minipage}
	}%
	\subfigure[t=0.85]{
		\begin{minipage}[t]{0.25\linewidth}
			\centering
			\includegraphics[width=\textwidth]{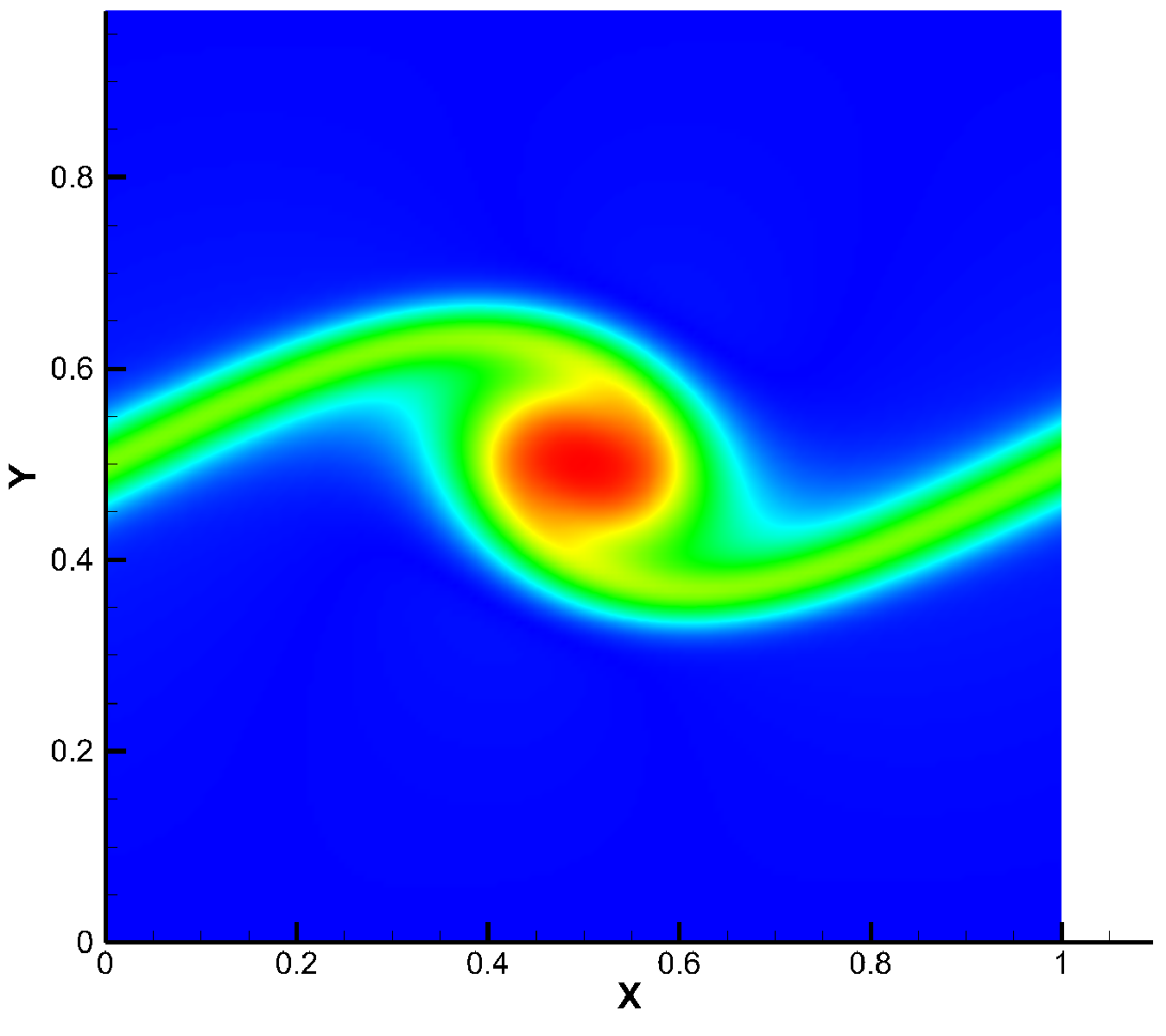}
		\end{minipage}
	}%
\subfigure[t=1]{
		\begin{minipage}[t]{0.25\linewidth}
			\centering
			\includegraphics[width=\textwidth]{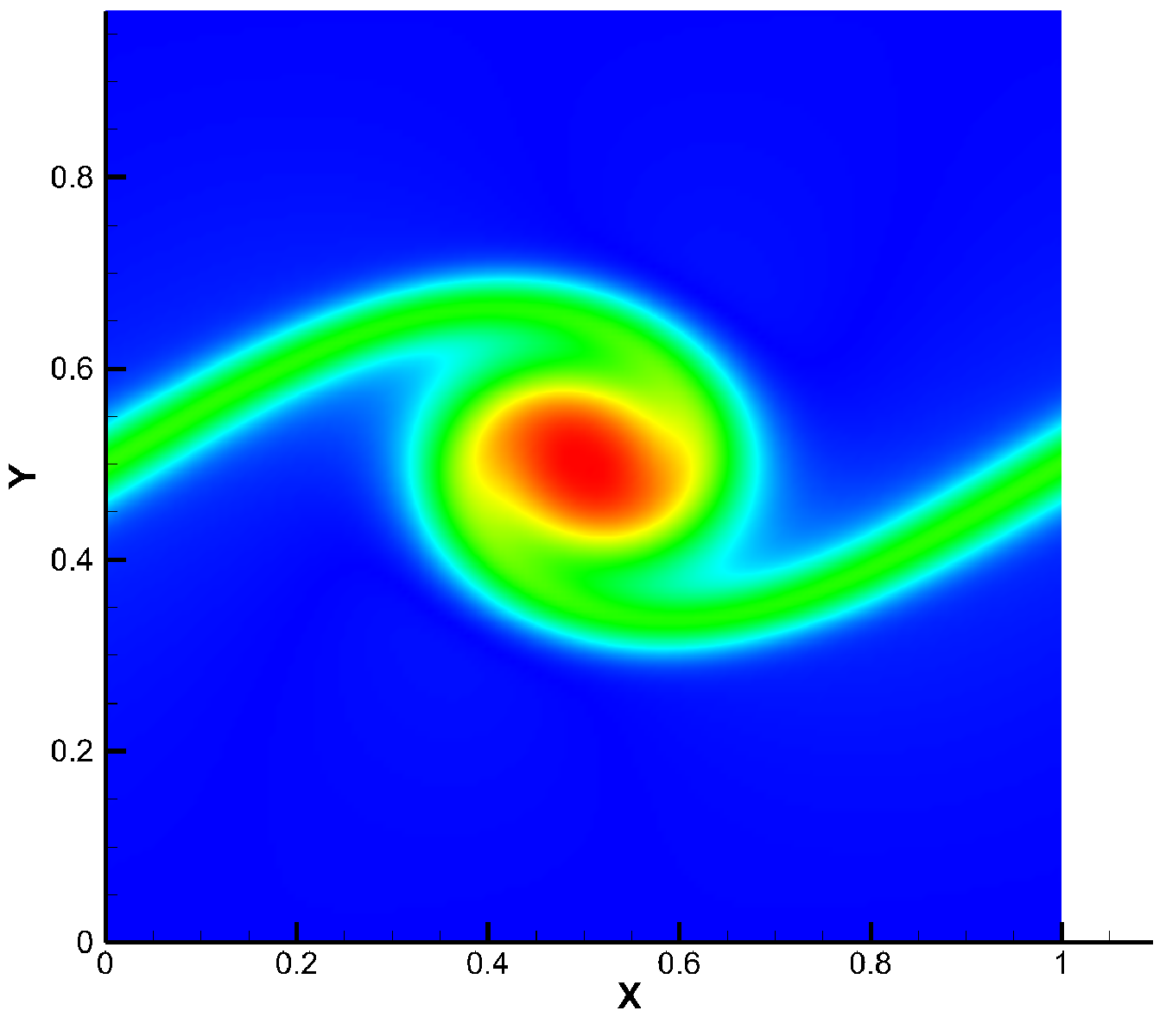}
		\end{minipage}
	}%
\\
	\subfigure[t=1.1]{
		\begin{minipage}[t]{0.25\linewidth}
			\centering
			\includegraphics[width=\textwidth]{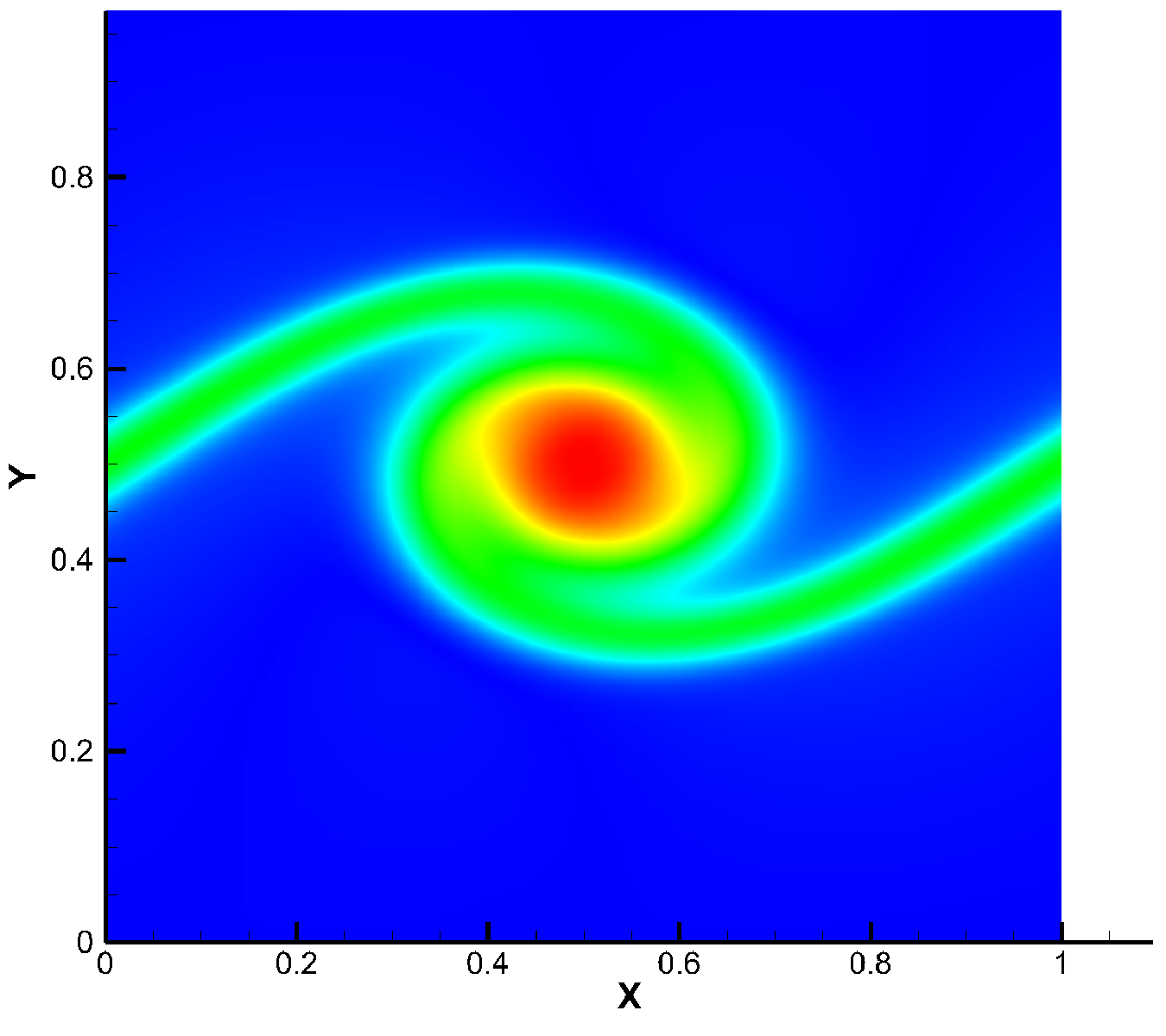}
		\end{minipage}
	}%
\subfigure[t=1.2]{
		\begin{minipage}[t]{0.25\linewidth}
			\centering
			\includegraphics[width=\textwidth]{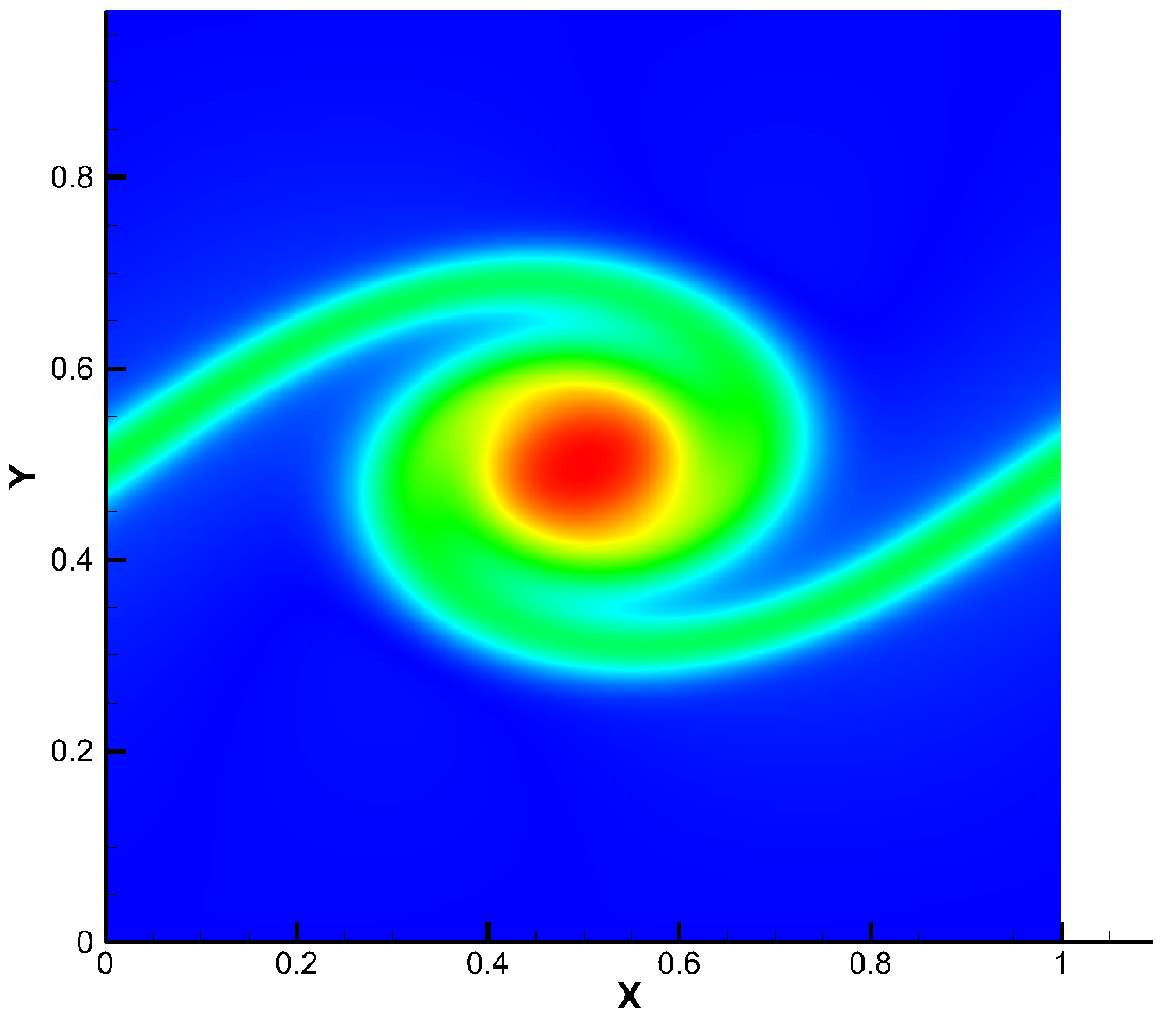}
		\end{minipage}
	}%
   \subfigure[t=1.4]{
		\begin{minipage}[t]{0.25\linewidth}
			\centering
			\includegraphics[width=\textwidth]{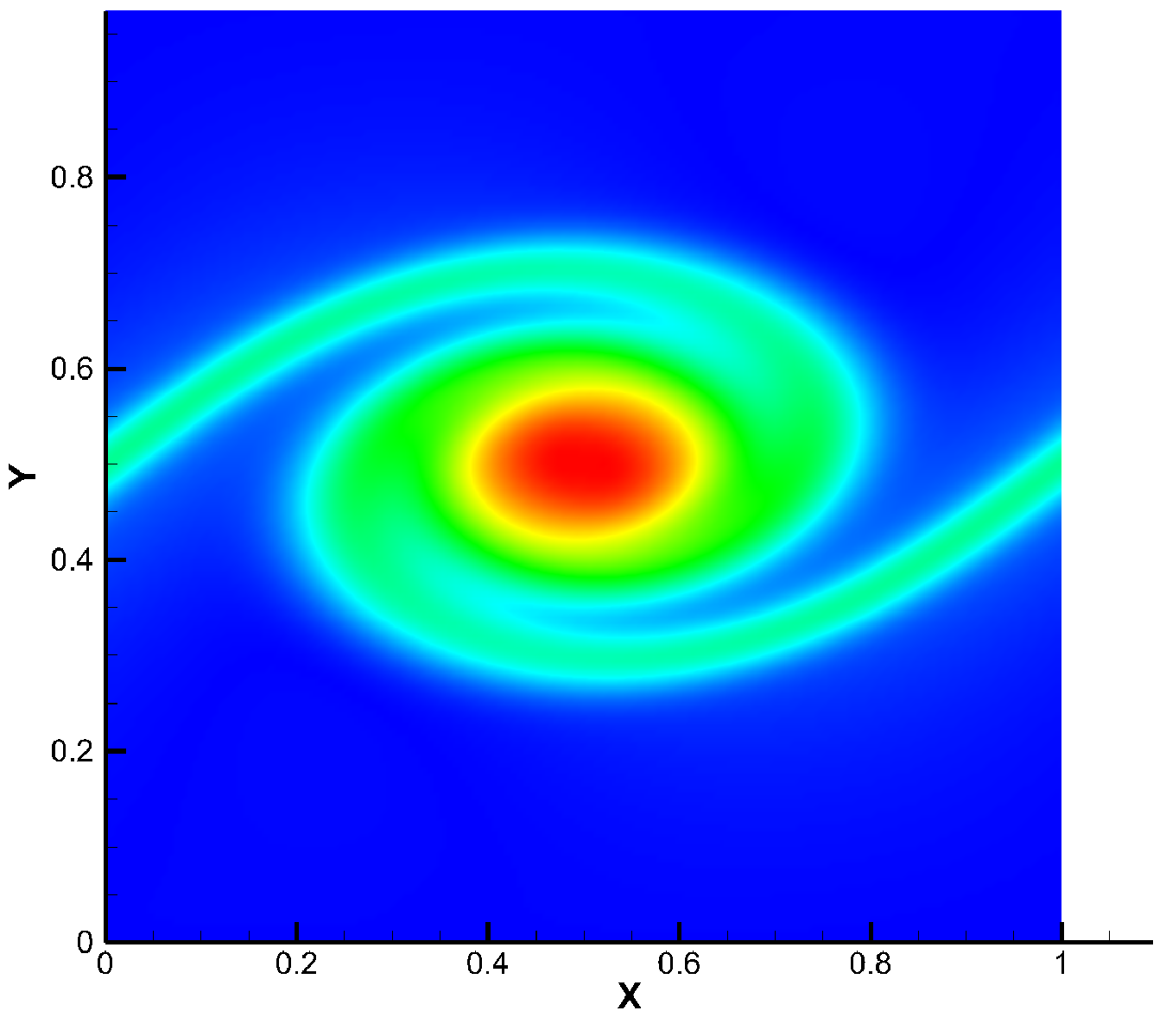}
		\end{minipage}
	}%
\subfigure[t=1.6]{
		\begin{minipage}[t]{0.25\linewidth}
			\centering
			\includegraphics[width=\textwidth]{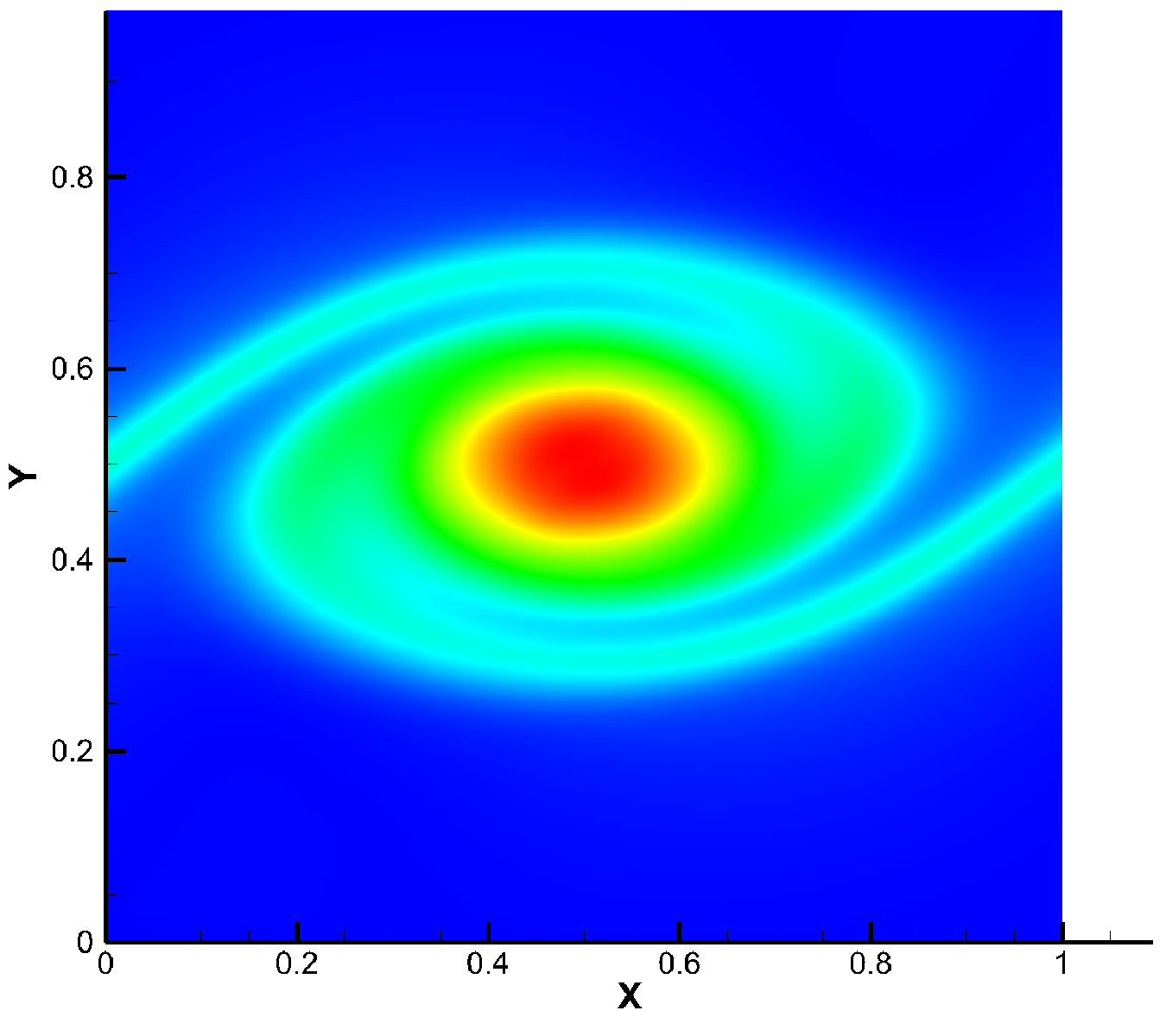}
		\end{minipage}
	}%
	\centering
	\caption{Single mode vorticity dynamics at time at t=0.001 (a), 0.6 (b), 0.85 (c), 1 (d), 1.1 (e), 1.2 (f), 1.4 (g), 1.6 (h) for case I.}
	\label{KH-vorticity}
\end{figure}

The snapshots of vorticity evolution are plotted in Figure \ref{KH-vorticity}. The fluids at the top and bottom flow in opposite directions, causing the vorticity to migrate towards the center of the region. As the vorticity  accumulates at the center, the interface starts to become more pronounced, and the amplitude of the instability increases. A roll-up phenomenon occurs, transforming the interface into a spiral that takes on a distinctive ``cat's eye'' configuration.

\subsubsection{3D Kelvin-Helmholtz instability}

In 3D simulation, we set the mesh size $h=1/16$, time step $\Delta t=1/1000$ and the following initial values:
\begin{eqnarray*}\label{K-H1}
\left\{
\begin{aligned}
\phi_{0}	&=  \tanh(\frac{z-0.5-0.01\sin(2\pi x)}{\sqrt{2}\gamma}),\\
\u_{0}	&=\big( \tanh(\frac{z-0.5-0.01\sin(2\pi x)}{\sqrt{2}\gamma}), 0, 0\big)^\top,\\
\B_{0}	&=\big(1, 0, 0\big)^\top.
\end{aligned}
\right.
\end{eqnarray*}
The magnetic induction field's boundary condition $(-1, 0, 0)^\top$ is imposed  on the faces where $y=0$, $y=1$, $z=0$, and $z=1$. On the upper boundary where $z=1$ and the lower boundary where $z=0$, the $u_{2}=u_{3}=0$ boundary condition is applied. In addition, the boundary conditions $u_{2}=0$ are equipped on the faces $y=0$ and $y=1$.

\begin{figure}[htbp]
	\centering
\subfigure[$t=0.001$]{
		\begin{minipage}[t]{0.25\linewidth}
			\centering
			\includegraphics[width=\textwidth]{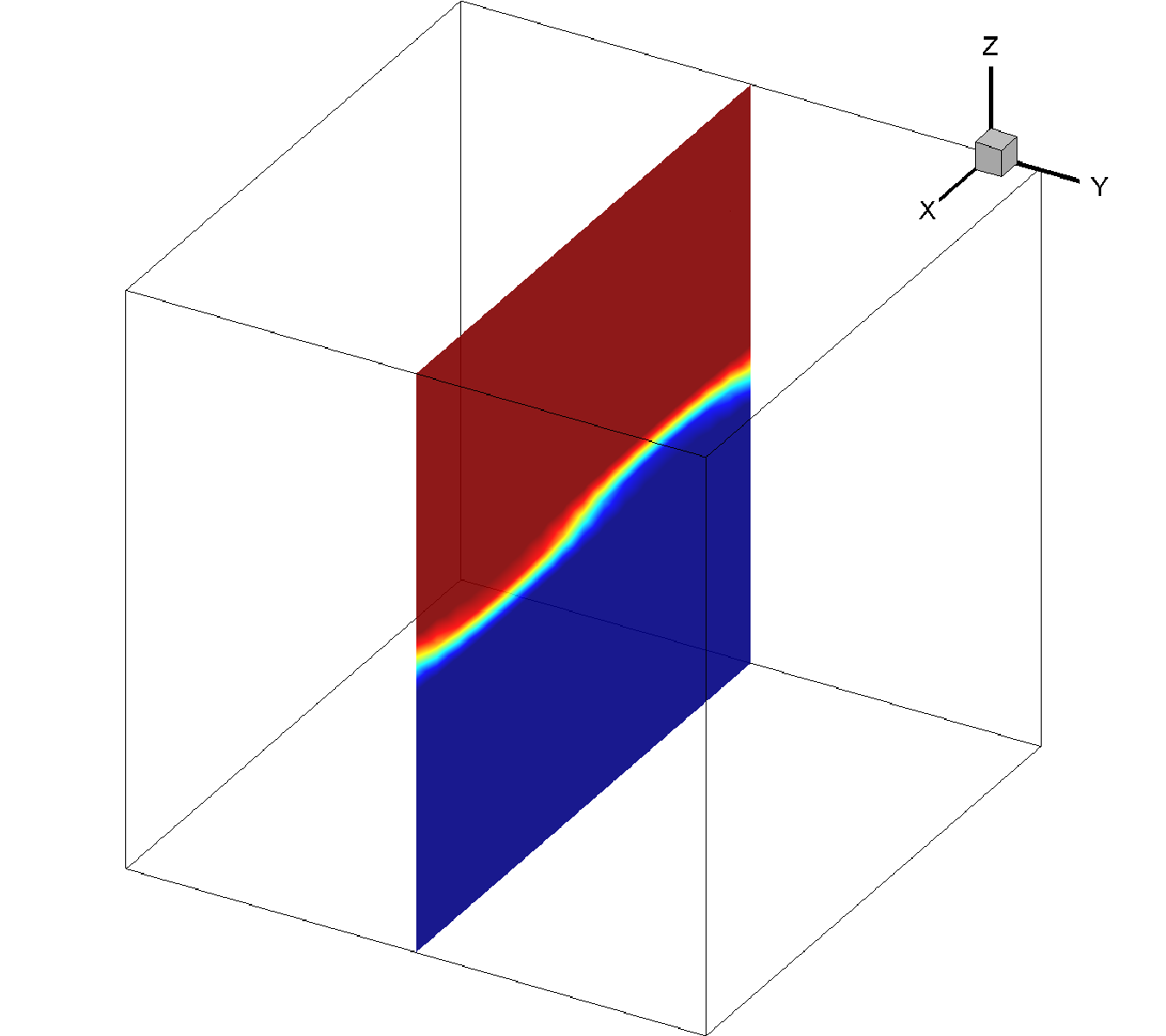}
		\end{minipage}
	}%
	\subfigure[t=0.6]{
		\begin{minipage}[t]{0.25\linewidth}
			\centering
			\includegraphics[width=\textwidth]{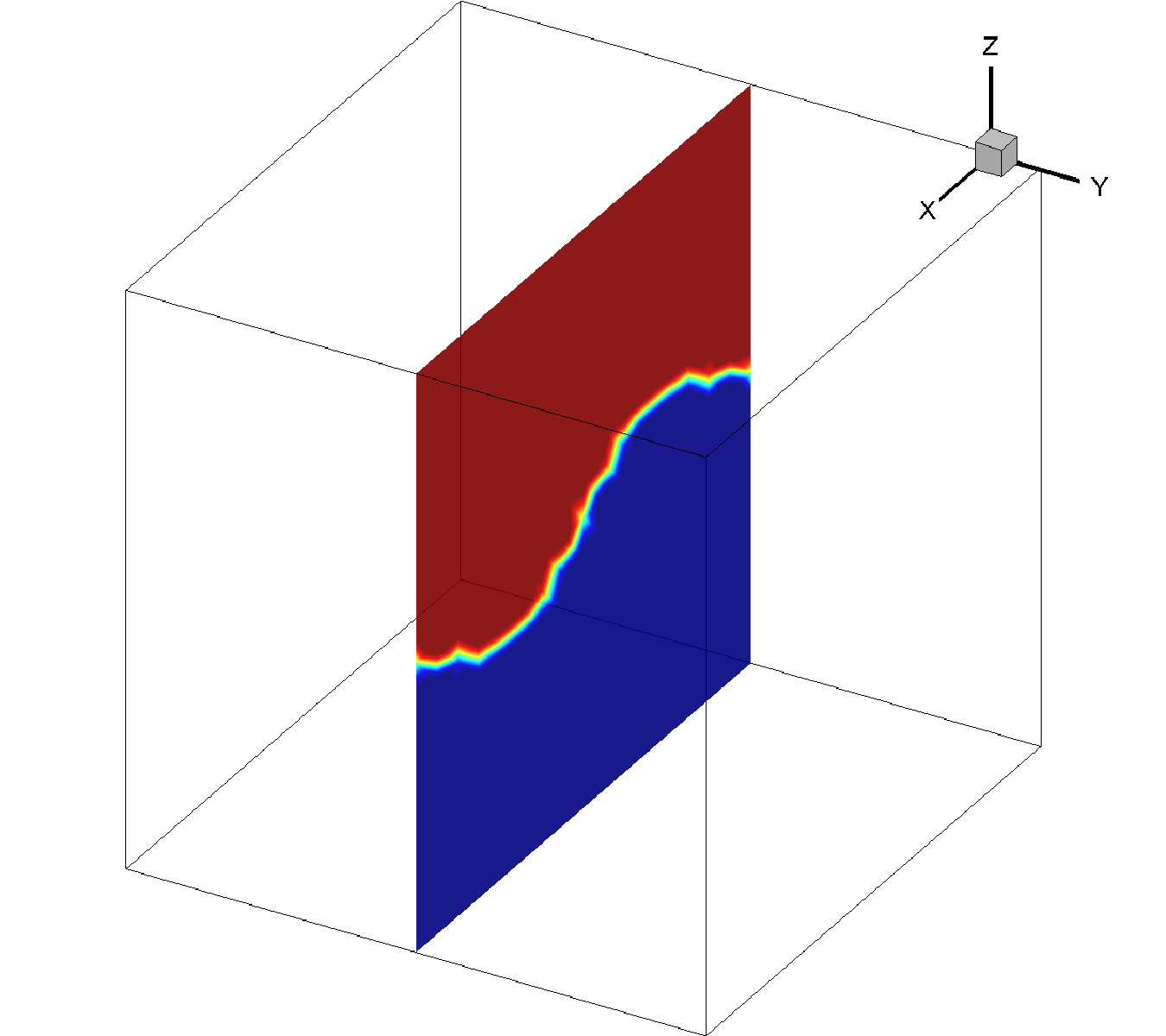}
		\end{minipage}
	}%
   \subfigure[t=1.3]{
		\begin{minipage}[t]{0.25\linewidth}
			\centering
			\includegraphics[width=\textwidth]{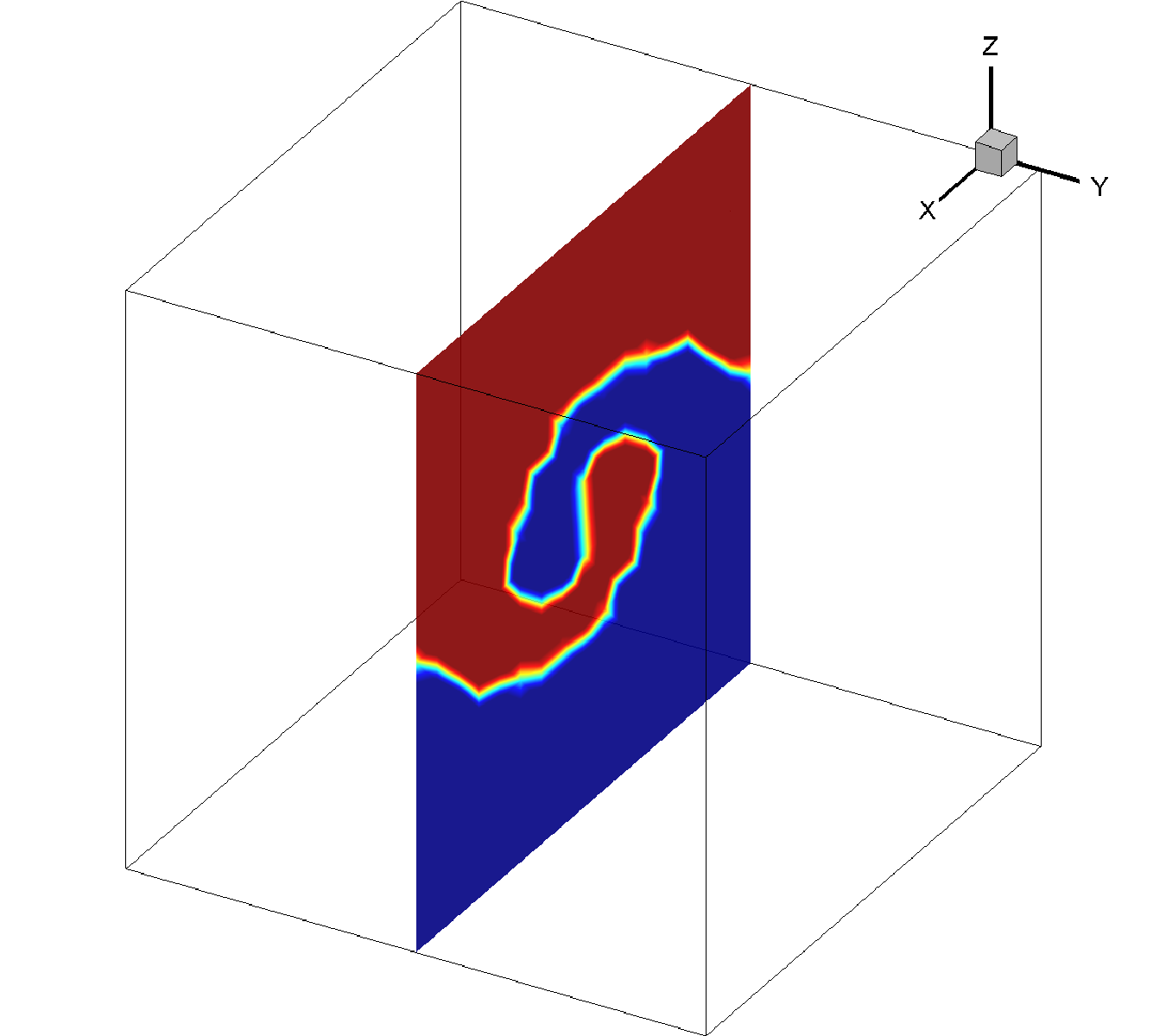}
		\end{minipage}
	}%
\subfigure[t=1.9]{
		\begin{minipage}[t]{0.25\linewidth}
			\centering
			\includegraphics[width=\textwidth]{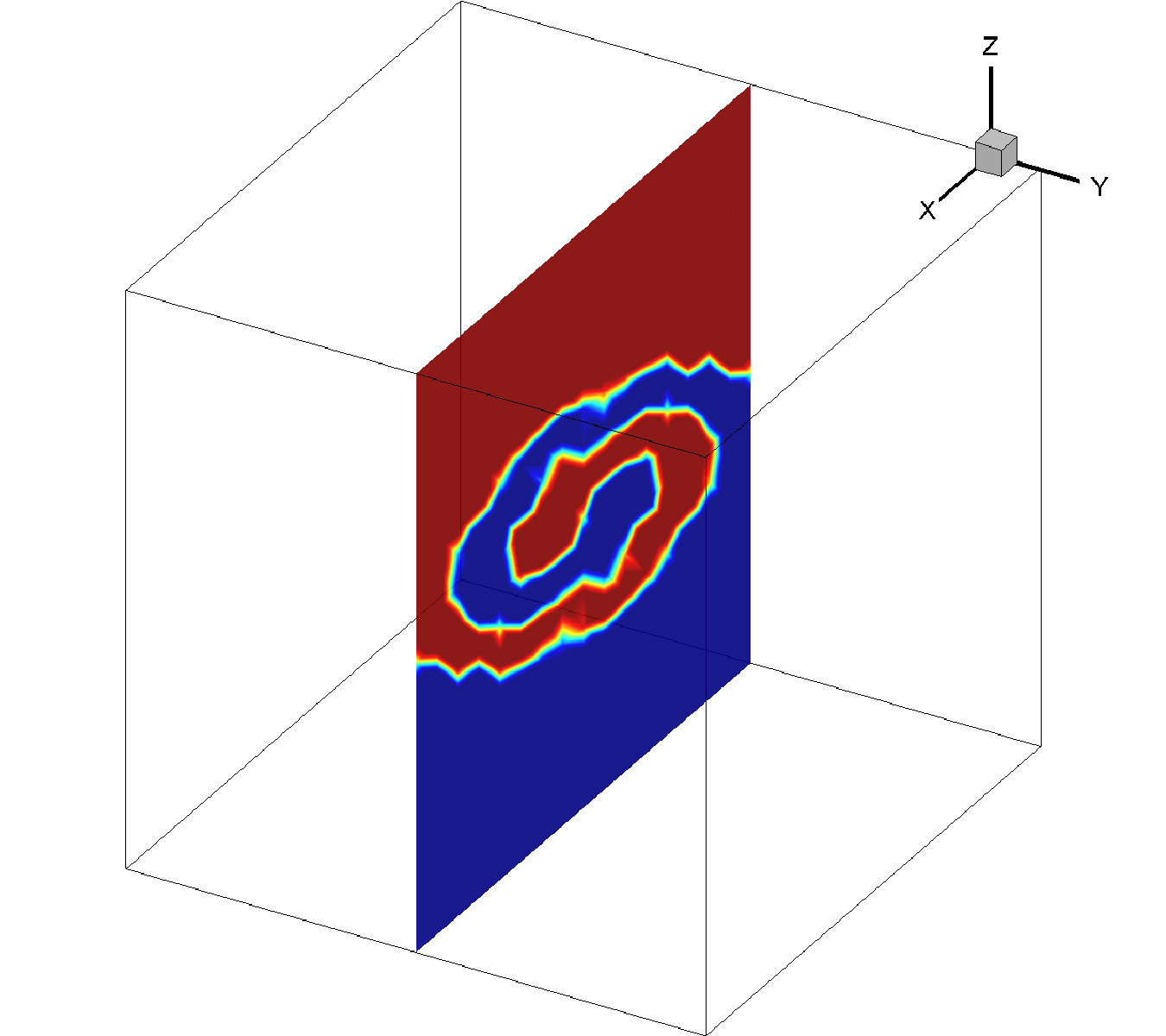}
		\end{minipage}
	}%
\\
\subfigure[t=0.001]{
		\begin{minipage}[t]{0.25\linewidth}
			\centering
			\includegraphics[width=\textwidth]{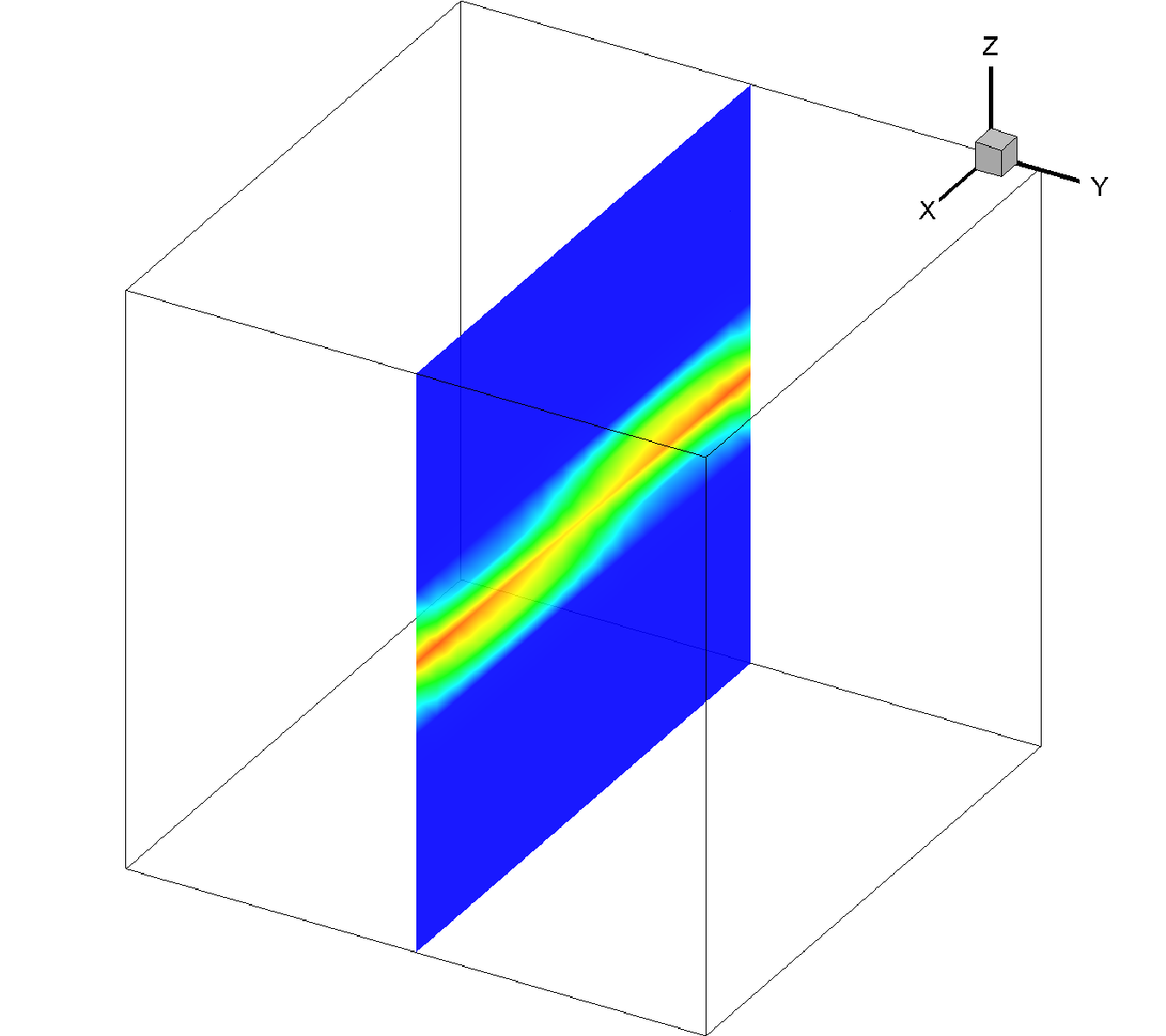}
		\end{minipage}
	}%
	\subfigure[t=0.6]{
		\begin{minipage}[t]{0.25\linewidth}
			\centering
			\includegraphics[width=\textwidth]{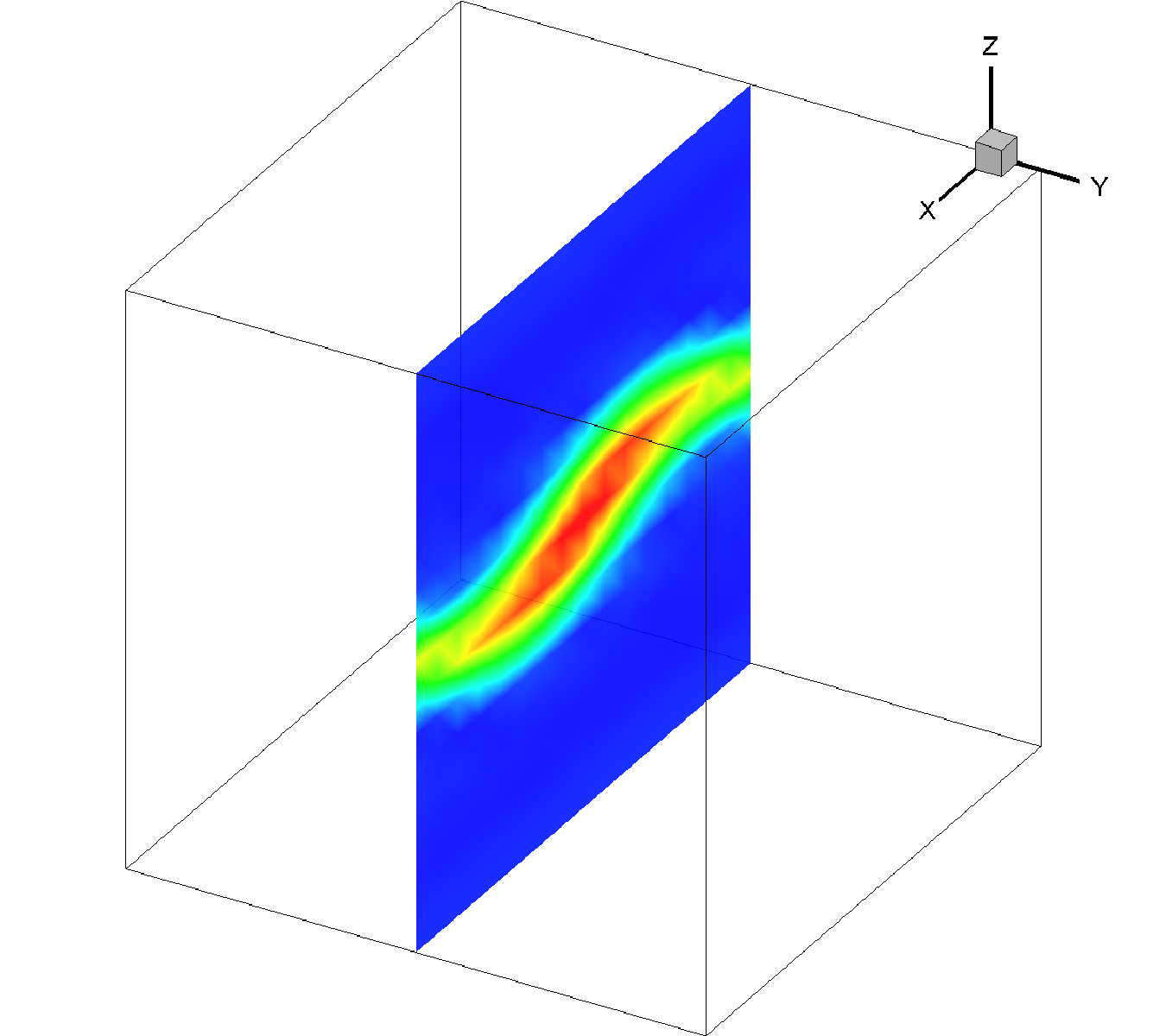}
		\end{minipage}
	}%
   \subfigure[t=1.3]{
		\begin{minipage}[t]{0.25\linewidth}
			\centering
			\includegraphics[width=\textwidth]{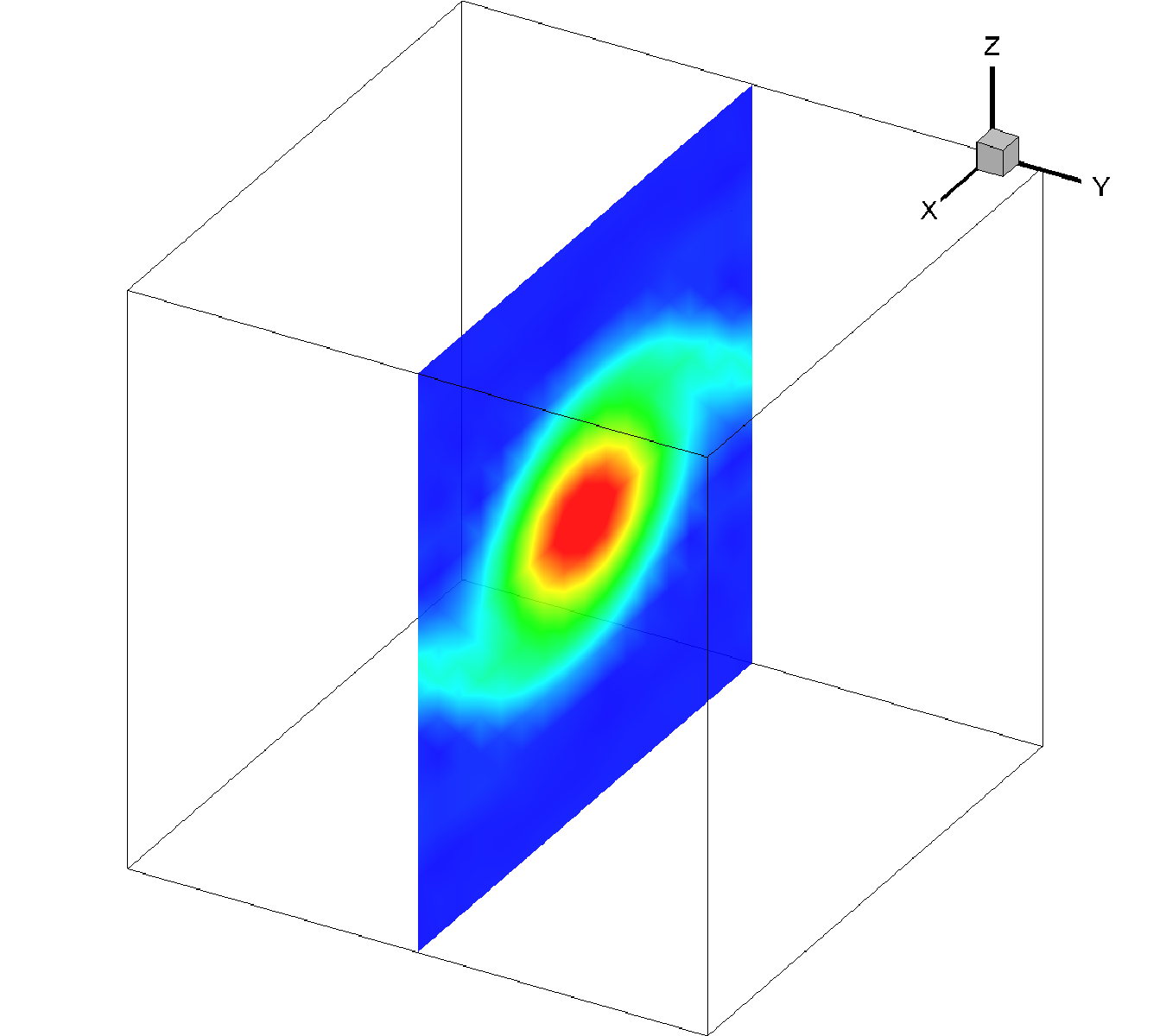}
		\end{minipage}
	}%
\subfigure[t=1.9]{
		\begin{minipage}[t]{0.25\linewidth}
			\centering
			\includegraphics[width=\textwidth]{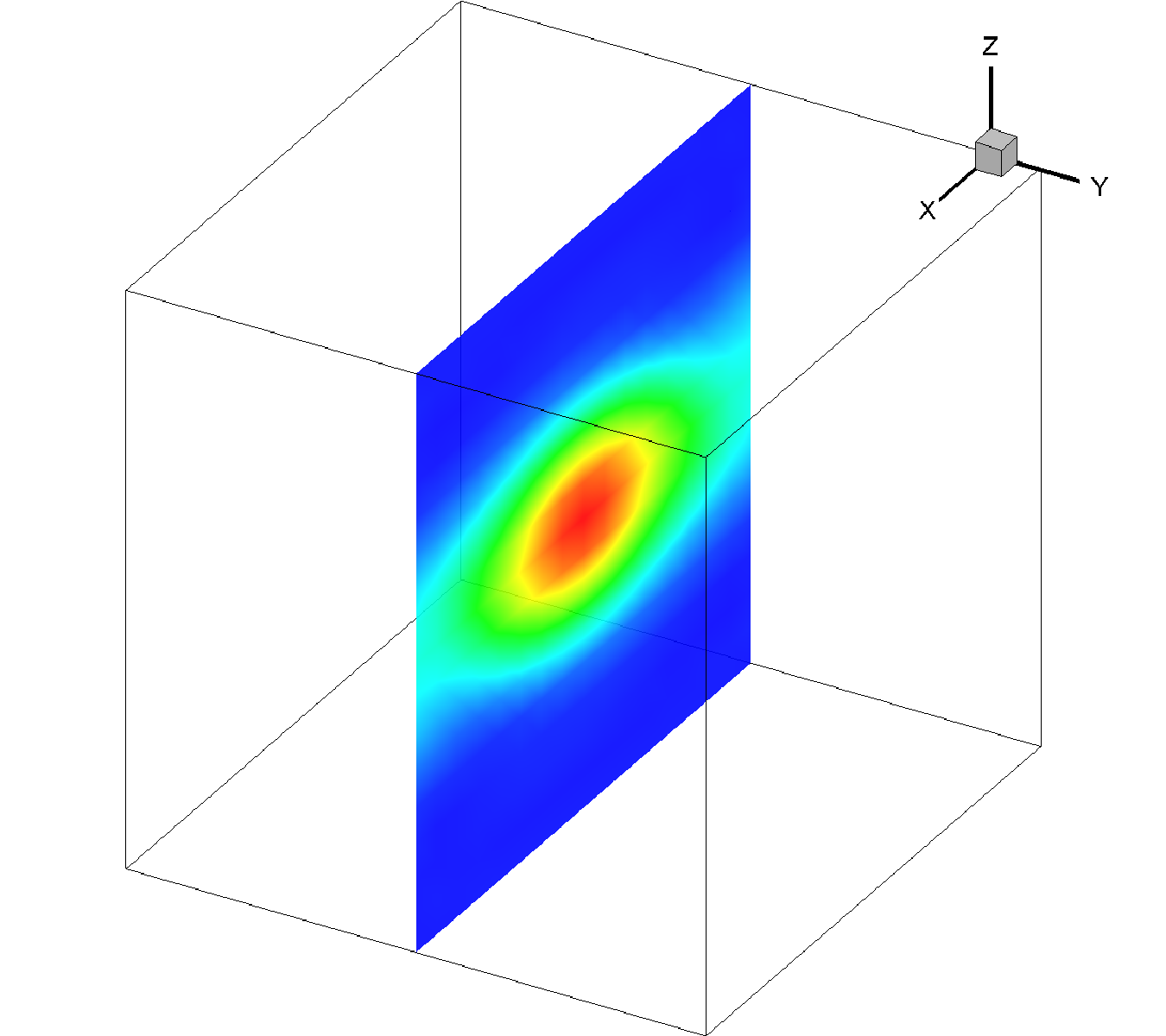}
		\end{minipage}
	}%
	\centering
	\caption{Temporal evolution of the phase field (upper), vorticity dynamics (lower) perturbed sinusoidal at t=0.001 (a), 0.6 (b),   1.3 (c), 1.9 (d) for case I 3D.}
	\label{KH-3D}
\end{figure}

We present the temporal evolution of the phase field (upper) and vorticity dynamics (lower), both of which are perturbed by a sinusoidal wave at $t=0.001, 0.6, 1.3, 1.9$ in Figure \ref{KH-3D}. The overall results are somewhat similar to the 2D cases, hence we explain more briefly here.  The phase field evolves over time, gradually bending and flipping in the interfacial region. The vorticity magnitude also evolves over time to exhibit a ``cat's eye'' pattern.

\section{Conclusion Remarks}\label{sec-conclusion}

In this paper, we develop the optimal error analysis for all variables  of a convex-splitting FEM for the  diffuse interface CH--MHD model in a convex domain.  When the magnetic induction field is discretized by the N\'{e}d\'{e}lec edge elements and the other variables are  discretized by the Lagrange elements, the optimal $\L^{2}$- and $\H^{1}$-norm error analysis is obtained with the help of the Ritz, Stokes and Maxwell  quasi-projections.   
 
Considering that the present results are only applicable for first-order numerical scheme, the optimal $\L^{2}$- and $\H^{1}$-norm error analysis for second-order scheme will be investigated in the future work.

\bibliographystyle{plain}
\bibliography{reference}







\end{document}